\documentclass[a4paper, 12pt]{article}
\usepackage{makeidx}
\usepackage{enumerate, theorem}
\usepackage{amsmath, amsfonts, amssymb}
\usepackage[height=22.5cm, width=15cm]{geometry}
\usepackage[all]{xy}
\usepackage{mathrsfs, yfonts}
\usepackage{graphicx}
\DeclareMathOperator{\C}{\mathbb{C}}
\newcommand{\parag}[1]{\paragraph{\sc{#1.}} }

\newtheorem{thm}{Theorem}[subsection]
\newtheorem{defn}[thm]{Definition}
\newtheorem{cor}[thm]{Corollary}
\newtheorem{prop}[thm]{Proposition}
\newtheorem{lemma}[thm]{Lemma}

\begin{document}

\title{On partial differential operators which annihilate the roots of the universal equation of degree $k$}

 
 \bigskip

 \author{Daniel Barlet\footnote{Institut Elie Cartan, G\'eom\`{e}trie,\newline
Universit\'e de Lorraine, CNRS UMR 7502   and  Institut Universitaire de France.}.\\
\quad \\
 $\hfill $  \qquad \qquad \qquad \qquad \qquad \qquad \qquad \qquad {\it A Mil\`{e}ne qui a support\'{e} courageusement}\\ 
 $\hfill $ \qquad \qquad \qquad \qquad \qquad \qquad \qquad \qquad {\it deux confinements avec cette \'equation}}
 
 \bigskip

\maketitle

\parag{Abstract} The aim of this paper is to study in details the regular holonomic $D-$module introduced in \cite{[B.19]} whose local solutions outside the polar  hyper-surface $\{\Delta(\sigma).\sigma_k = 0 \}$ are given by the local system generated by the local branches of the  multivalued function which is the root  of the universal  degree $k$ equation  $z^k + \sum_{h=1}^k (-1)^h.\sigma_h.z^{k-h} = 0 $. Note that it is surprising that this regular holonomic  $D-$module is given by the quotient of $D$ by a left ideal which has very simple explicit generators despite the fact it necessary encodes the analogous systems for any root  of the universal  degree $l$ equation for each  $l \leq k$.\\
 Our main result is to relate this $D-$module with the minimal extension of the  irreducible local system associated to the difference of two branches of the multivalued function defined above. Then we obtain again a very simple explicit description of this minimal extension in term  of the generators  of its left ideal in the Weyl algebra.\\
 As an application we show how these results allow to compute the Taylor expansion of the root near $-1$ of the equation $z^k + \sum_{h=-1}^k (-1)^h.\sigma_h.z^{k-h} - (-1)^k = 0 $.
 
 \parag{AMS Classification} 32 C 38 - 35 Q 15 - 35 C 10 -32 B 06
 
  \newpage

 \tableofcontents

\section{Introduction} 

There are several ways to define an interesting  function. Of course the simplest one is to give its value at each point by an explicite finite formula or as a sum of an infinite series (converging somewhere at least). Another way is to give a functional equation which characterizes it. A third approach is to give a partial differential system which has our function $f$ as its unique solution (up to normalisation). \\
For instance, the function $f(z) = e^z$ may be define as 
\begin{enumerate}
\item $f(z) = \sum_{n=0}^\infty \frac{z^n}{n!}$.
\item $f(z+z') = f(z).f(z') $ with  $f(0) = 1$ and $f(1) = e$.
\item $\frac{\partial f}{\partial z} = f $ and $f(0) = 1$.
\end{enumerate}
In general, to increase our understanding of such a function it is useful to dispose of at least two kinds of the characterization as above. For instance, in the basic example of $e^z$ the description 3. gives easily the formula 1 and also the functional equation 2. Note that the third approach will often leads to a description of the first kind via the Taylor expansion at least when we dispose of a regular holonomic system defining $f$ which is enough simple and suitably described in order to allow an inductive explicite computation of the coefficients of the Taylor expansion. But this means that we are at least able to well describe essentially all partial differential operators which annihilate $f$.\\
We shall consider, in this paper, the case of the multivalued function $z(\sigma) $ on $\C^k$, with  $\sigma := (\sigma_1, \dots, \sigma_k)$, which is defined as the root of the universal monic polynomial of degree $k$:
$$P_\sigma(z) := \sum_{h=0}^k (-1)^h.\sigma_h.z^{k-h} \quad {\rm with \ the \ convention} \quad \sigma_0 \equiv 1 .$$
It is well known that the description of this function with the first approach is quite difficult (at least for $k \geq 5$). The definition given above of this multivalued function may be seen as a description of the second kind.\\
The aim of this paper is to give a description of the third kind which  characterizes  this multivalued function. More precisely we describe completely the structure of the regular holonomic $D_N-$module\footnote{$N$ will be $\C^k$ with coordinates $\sigma_1, \dots, \sigma_k$.} $D_N\big/\mathcal{J}_\lambda$ where $\mathcal{J}_\lambda$ is the left ideal in $D_N$ which annihilates $z^\lambda(\sigma)$.\\
 The case where $\lambda$ is in $\mathbb{Z}$ is of special interest (for $\lambda = 0$  the left ideal $\mathcal{J}_0$ will be defined in a natural way inside the annihilator of the function $1$).\\
As an application of the structure theorem for the $D-$module $\mathcal{N}_1 := D_N\big/\mathcal{J}_1$  we compute the Taylor series at the point $\sigma^0 = (0, 0, \dots, -1)$ of the holomorphic function of $\sigma_1, \dots, \sigma_k$  which gives the root of the polynomial 
$$z^k + \sum_{h=1}^k (-1)^h.\sigma_h.z^{k-h} - (-1)^k = 0$$
 which is near $-1$. The fact that the associated $D_N-$module corresponding to $\lambda = 1$ is not simple make this computation quite complicate when we use the $D_N-$module $\mathcal{N}_1$ itself. \\
But with the remark that $z(\sigma) - \sigma_1/k$ is a solution of the simple part $\mathcal{N}^\square_1$ of this $D-$module deduced from the structure theorem \ref{fin 1}  we obtain a complete explicite computation of the Taylor series at $\sigma^0$  of the root which is near $-1$.\\
Of course this method is valid to compute (with some more numerical complications, but without theoretical difficulty) the Taylor expansion of any (uni-valued) holomorphic branch of the multivalued function $z(\sigma)$ near any point $\sigma^0 \in N$. 
\newpage

\section{The $D-$modules $\mathcal{W}$ and $\mathcal{M}$}

\parag{Notations}
We fix in the sequel and integer $k \geq 2$. Let $\C[\sigma]\langle \partial \rangle$ be the Weyl algebra in the variables $\sigma_1, \dots, \sigma_k$. We shall note $N = \C^k$ which is the target of the quotient map
 $$q : M:=  \C^k \to \C^k\big/\mathfrak{S}_k  = N \simeq \C^k$$
  by the natural action of the permutation group $\mathfrak{S}_k$ on $\C^k$. We shall  note $T_q$ its tangent map.\\
  Then $D_N$ denotes the sheaf of holomorphic differential operators on $N$ and we shall use the same notations for modules on $\C[\sigma]\langle \partial \rangle$ and for the corresponding sheaves of $D_N-$modules on $N$.\\
  For basic results on $D-$modules the reader may consult \cite{[Bjork]} or \cite{[Borel]}.\\

\subsection{The $D-$module $\mathcal{W}$}

 In this section we shall consider  the $D_N-$module $\mathcal{W} := D_N\big/\mathcal{A}$ where $\mathcal{A}$ is the left ideal sheaf in $D_N$ generated by
\begin{equation}
A_{p, q} := \partial_p\partial_q - \partial_{p+1}\partial_{q-1} \quad {\rm for} \quad p \in [1, k-1] \quad {\rm and} \quad q \in [2, k] 
\end{equation}

\parag{Notations} Let $\mathcal{W}(m)$ be the sub$-\mathcal{O}_N-$module in $\mathcal{W}$ of the classes induced by germs in $D_N(m)$ the partial differential operators of order at most equal to $m$. \\
As we have $\mathcal{A}(m) := D_N(m) \cap \mathcal{A} = \sum_{p, q} D_N(m-2).A_{p, q}$  we have $\mathcal{W}(m) = D_N(m)\big/\mathcal{A}(m) $ which injects in $\mathcal{W}$ and 
$$\mathcal{W} = \cup_{m \geq 0} \mathcal{W}_m.$$
Note that $\mathcal{A}(1) = 0$ so $\mathcal{W}(1) = D_N(1)$.\\

It is clear that the characteristic variety of the $D_N-$module $\mathcal{W}$ is equal to $N \times S(k)$ in the cotangent bundle $ T_N^* \simeq N \times \C^k$ of $N$,  where $S(k)$ is the algebraic cone in $\C^k$ defined by the equations
\begin{equation}
 \eta_p.\eta_q - \eta_{p+1}.\eta_{q-1} = 0 \quad \forall p \in [1, k-1] \quad {\rm and} \quad  \forall q \in [2, k].
 \end{equation}
We describe this two-dimensional cone and the corresponding ideal in the appendix (see section 6). We shall use in the present section the following results  which are proved in the appendix (prop. 6.0.5 and cor. 6.0.6).

\begin{prop}\label{local description}
Let $L_1 := \{\eta_1 = 0\} \cap S(k)$ and $L_k := \{\eta_k = 0 \} \cap S(k)$. Then $L_1$ is the line directed by the vector $(0, \dots, 0, 1)$ and $L_k$ the line directed by the vector $(1, 0, \dots, 0)$.
The maps $\varphi_1 : S(k)\setminus L_1 \to \C^*\times \C$ and $\varphi_k : S(k) \setminus L_k \to \C^*\times \C$ which are defined by the formulas
\begin{equation}
\varphi_1(\eta) := (\eta_1, -\eta_2/\eta_1) \quad {\rm and} \quad \varphi_k(\eta) := (\eta_k, -\eta_{k-1}/\eta_k)
\end{equation}
are isomorphisms. So $S(k)\setminus \{0\}$ is smooth and connected surface.
\end{prop}

\begin{cor}\label{irreduct.}
The ideal  of $\C[\eta]$ defined by the equations in $(2)$ is prime. Moreover $S(k)$ is a normal surface.
\end{cor}

\parag{Notation} For $\alpha$ in $\mathbb{N}^k$  define $q := \vert \alpha\vert := \sum_{h=1}^k \alpha_h$ and $r := w(\alpha) = \sum_{h=1}^k h.\alpha_h$.

\begin{defn}\label{bi-homog.}
Let  $P$ be a germ of section  of $D_N$. We say that $P$ is {\bf bi-homogeneous of type $(q, r)$} if we may write 
$$ P = \sum_{\vert \alpha\vert = q, w(\alpha) = r} a_\alpha.\partial^\alpha $$
where $a_\alpha$ are  germs of holomorphic functions in $N$.
\end{defn}

It is clear that  any germ $P$  of section of $D_N$ has a unique decomposition
 $$P = \sum_{q, r} P_{q, r}$$
 where $P_{q, r}$ is a  bi-homogeneous germ of section of $D_N$  of type $(q, r)$. Note that this sum is finite because for a given order $q$ the corresponding type $(q, r)$ has non zero representative only when $r$ is in $[q, k.q]$.

\begin{lemma}\label{tres simple}
Let $P$ be a germ section of $D_N$ and write the decomposition of $P$ in  its bi-homogeneous components  as $P = \sum_{r = 0}^N  P_{q, r} $. Then $P$ is a germ of section in $\mathcal{A}$ if and only if for each type $(q, r)$ $P_{q, r}$ is a germ of section in $ \mathcal{A}$.
\end{lemma}

\parag{proof} It is clear that $P$ is in $\mathcal{A}$ when each $P_{q, r}$ is in $ \mathcal{A}$. Conversely, assume that $P $ is in $ \mathcal{A}$. Then we may write 
$$ P = \sum_{(i, j)} B_{i, j}.A_{i, j} \quad {\rm with}  \quad i \in [1, k-1] \quad {\rm and} \quad j \in [2, k] $$
and where $B_{i, j}$ are germs of sections of $D_N$.\\
Write $B_{i, j} = \sum_{q, r} (B_{i, j})_{q, r} $ the decomposition of $B_{i, j}$ in  its bi-homogeneous components; this gives 
$$ P = \sum_{q, r}\sum_{p \geq 0} \big(\sum_{i+j=p} (B_{i, j})_{q, r}.A_{i, j}\big) $$
where $\sum_{i+j=p} (B_{i, j})_{q, r}.A_{i, j}$ is bi-homogeneous of type $(q+2, r+p)$ for each $(i, j)$ such that  $i+j=p$. This implies that $P_{q, r}$ is equal to the sum $ \sum_{i+j+s = r} (B_{i, j})_{q-2, s}.A_{i, j}$. So each $P_{q, r}$ is a germ of section in $\mathcal{A}$.$\hfill \blacksquare$\\

\begin{lemma}\label{base}
 The class  of  $\partial^\alpha$ in $\mathcal{W}$ only depends on $q := \vert \alpha\vert$ and $r := w(\alpha)$.It will be denoted $y_{q, r} $. Moreover, if $\mathcal{W}_q$ is the sub$-\mathcal{O}_N-$module of $\mathcal{W}$ generated by the $y_{q, r} $ for $r \in [q, k.q]$, $\mathcal{W}_q$ is a free $\mathcal{O}_N-$module of rank $k.q - q + 1$  with basis  $y_{q, q}, y_{q, q+1}, \dots, y_{q, k.q}$ and, as $\mathcal{O}_N-$module, we have the direct decompositions: 
\begin{equation}
 \mathcal{W}(m) = \oplus_{q=0}^m \mathcal{W}_q \quad {\rm and} \quad  \mathcal{W} = \oplus_{q \in \mathbb{N}} \ \mathcal{W}_q .\\
 \end{equation}
 \end{lemma}
 
Remark that the action of $D_N$ on $\mathcal{W}$ is defined by  
\begin{equation}
\partial_j(y_{q, r}) = y_{q+1, r+j} \quad \forall j \in [1, k] \quad \forall q \in \mathbb{N} \quad {\rm and} \quad \forall r \in [q, k.q]
\end{equation}
 and that $\partial_j\mathcal{W}_q \subset \mathcal{W}_{q+1}$.

\parag{proof} The fact that the class induced by $\partial^\alpha$ in $\mathcal{W}$ depends only on $\vert \alpha \vert$ and $w(\alpha)$  is a direct consequence of the fact that the class induced by $x^\alpha$ in $\C[x_1, \dots, x_k]\big/IS(k)$ only depends on $q := \vert \alpha\vert$ and $r := w(\alpha)$ (see proposition \ref{le retour} in the appendix). Then it is clear that $y_{q, q}, y_{q, q+1}, \dots, y_{q, k.q}$ is a $\mathcal{O}_N-$basis of $\mathcal{W}_q$ looking at the symbols and using the appendix (section 6) over  the sheaf of $\C-$algebras $\mathcal{O}_N$.$\hfill \blacksquare$\\

The global polynomial solutions of  the $D_N-$module $\mathcal{W}$ are described by our next lemma.

\begin{defn}\label{les mqr} For each $q \in \mathbb{N}$ and each $r \in [q, k.q]$ define the polynomial
\begin{equation}
m_{q, r}(\sigma) := \sum_{\vert\alpha\vert = q, w(\alpha) = r} \frac{\sigma^\alpha}{\alpha !}
\end{equation}
\end{defn}

\begin{lemma}\label{sol. W}
Any $m_{q, r} \in \C[\sigma_1, \dots, \sigma_k]$  is annihilated by the left ideal $\mathcal{A}$ in $D_N$ and if a polynomial  $P \in \C[\sigma_1, \dots, \sigma_k]$ is annihilated by $\mathcal{A}$, $P$ is, in a unique way, a $\C-$linear combination of the $m_{q, r}$ for $q \geq 0$ and $r \in [q, k.q]$ which gives the bi-homogeneous decomposition of $P(\partial_1, \dots, \partial_k) \in \C[\partial_1, \dots, \partial_k]$ (see lemma \ref{tres simple}).
\end{lemma}

\parag{Proof} First we shall verify that each polynomial $m^q_r$ is annihilated by each $A_{i, j}$ for all $ i \in [1, k-1]$ and all $j \in [2, k]$.  We have for each $(i, j) \in [1, k]^2$:

$$\partial_i\partial_j (m_{q, r})(\sigma) = \sum_{\vert \beta\vert = q-2, w(\beta) = r-(i+j)} \frac{\sigma^\beta}{\beta !}  = m_{q-2, r-(i+j)(\sigma)}$$
because  $\alpha_i.\alpha_j = 0$ implies $\partial_i\partial_j \sigma^\alpha = 0 $. The right hand-side above only depends on $i+j$ for $q$ and $r$ fixed. This is enough to conclude our verification.\\
Note also that the uniqueness is obvious because of the uniqueness of the Taylor expansion of a polynomial.\\
Let now $P := \sum_{\alpha} c_\alpha.\frac{\sigma^\alpha}{\alpha !} $ a polynomial in $\C[\sigma]$ which is annihilated by the left ideal $\mathcal{A}$ in $D_N$. We want to show that if $\alpha $ and $ \beta$ satisfy $\vert \alpha\vert = \vert \beta\vert$ and $w(\alpha) = w(\beta)$  we have $c_\alpha = c_\beta$. It is enough to prove this equality when there exist $i \in [1, k-1], j \in [2, k]$ and $\gamma \in \mathbb{N}^k$ such that $\sigma^\alpha = \sigma_i.\sigma_j.\sigma^\gamma$ and $\sigma^\beta = \sigma_{i+1}.\sigma_{j-1}.\sigma^\gamma$ by definition of the equivalence relation\footnote{which is noted $\alpha \sharp \beta$ in the section 6.} given by $\vert \alpha\vert = \vert \beta\vert$ and $w(\alpha) = w(\beta)$. In this case the coefficient of $\sigma^\gamma/\gamma !$ in $\partial_i\partial_j P$ is $c_\alpha$ and in $\partial_{i+1}\partial_{j-1} P$ is $c_\beta$. So they are equal.$\hfill \blacksquare$\\
 
 It is easy to see that an entire holomorphic function $F : N \to \C$ is solution of $\mathcal{W}$ if and only if its Taylor series at the origin may be written, for some choice of  $c_{q, r} \in \C$: 
 $$ F(\sigma) = \sum_{q, r}  c_{q, r}.m_{q, r}(\sigma) .$$
 In the same way, a holomorphic germ $f : (N, \sigma^0) \to (\C, z^0)$ is solution of $\mathcal{W}$ if and only if its Taylor series may be written in the form
 $$ f(\sigma^0 + \sigma) =   \sum_{q, r}  c_{q, r}.m_{q, r}(\sigma) $$
 with $c_{0, 0} = z^0$.

\subsection{The $D_N-$module $\mathcal{M}$}

\begin{defn}\label{ideal}
Let $m \in [2, k]$ be an integer and define the second order differential operators in the Weyl algebra $\C[\sigma]\langle \partial \rangle$
\begin{equation}
\mathcal{T}^m := \partial_1\partial_{m-1} + \partial_m.E.  \quad {\rm where} \quad E := \sum_{h=1}^k \sigma_h.\partial_h
\end{equation}
Then define the left ideal $\mathcal{I}$ in $D_N$ as 
\begin{equation}
\mathcal{I} := \mathcal{A} + \sum_{m=2}^k D_N.\mathcal{T}^m
\end{equation}
and let $\mathcal{M}$ be  the $D_N-$module 
 \begin{equation}
\mathcal{M} := D_N\big/\mathcal{I}
\end{equation}
\end{defn}

We shall now recall and precise some results of \cite{[B.19]}.

Let $Z$ be the complex (algebraic) subspace in $N \times \C^k$ (with coordinates $\sigma_1, \dots, \sigma_k, \eta_1, \dots, \eta_k$)  defined by the ideal of $(2, 2)-$minors of the matrix
\begin{equation}
 \begin{pmatrix} \eta_1& -l_{\sigma}(\eta) \\ \eta_2 & \eta_1 \\ . & . \\ . & . \\ \eta_k & \eta_{k-1} \end{pmatrix} 
 \end{equation}
 where $l_\sigma(\eta) := \sum_{h=1}^k \sigma_h.\eta_h$. We shall note $I_Z$ the ideal of $\mathcal{O}_{N \times \C^k}$ generated by these minors and by $p_*I_Z$ its direct image by the projection
 $p : N\times \C^k \to N$. For each integer $q \geq 0$ the sub-sheaf  $p_*I_Z(q)$ of sections of $p_*I_Z$ which are homogeneous of degree $q$ along the fibers of $p$ is a coherent $\mathcal{O}_N-$module.

 \begin{prop}\label{charact. 1}
 The complex subspace  $Z$ is reduced, globally  irreducible and $Z$ is the characteristic cycle of the $D_N-$module $\mathcal{M}$.
\end{prop}

\parag{Proof} The fact that $Z$ is globally irreducible is already proved in \cite{[B.19]} prop. 4.2.6  as $Z$ is conic over $N$. This implies that $Z$  is reduced as a complex sub-space:\\
Assume that $I_Z $ is not equal to the reduced ideal $I_{\vert Z\vert}$ of the complex analytic subset $\vert Z\vert$ in $N \times \C^k$. By homogeneity in the variables $\eta_1, \dots, \eta_k$
there exists $q \geq 0$ such that the quotient $\mathcal{Q}(q) := I_{\vert Z\vert} (q)\big/ I_Z(q)$ is not $\{0\}$ and then the coherent sheaf $p_*(\mathcal{Q}(q))$ is not $\{0\}$ on $N$. But this contradicts
the fact that any global section on $N$ of $p_*I_{\vert Z\vert}(q)$ is a global section on $N$ of $p_*I_Z(q)$ which is the content of the proposition 4.2.6 in {\it loc. cit.} \\
To complete the proof that $Z$ is the characteristic cycle of $\mathcal{M}$ it is enough to see that the symbol of any germ $P$ of section in $\mathcal{I}$ vanishes on $\vert Z\vert$. This is obvious by definition of $I_Z$.$\hfill \blacksquare$\\

The following proposition, which is a local version of the theorem 5.1.1 in \cite{[B.19]}, will be useful:

 \begin{prop}\label{charact. 0}
 Let $\mathcal{I}_+$ the left ideal in $D_N$  of germs of differential operators $P$ such that $P(N_m) = 0$ for each Newton polynomial $N_m, m \in \mathbb{N}$. Then $\mathcal{I}_+ = \mathcal{I}$
 \end{prop}
 
 \parag{Proof} The proposition 4.1.2 in \cite{[B.19]} already proves the inclusion $\mathcal{I} \subset \mathcal{I}_+$. To prove the other inclusion we shall argue by contradiction. So assume that at some point $\sigma^0$ in $N$  we have $\mathcal{I}_{+, \sigma^0} \setminus \mathcal{I}_{\sigma^0} \not= \emptyset$ and let $P$ be in $\mathcal{I}_{+, \sigma^0} \setminus \mathcal{I}_{\sigma^0}$ with minimal order say $q$. Thanks to the proposition 4.2.8 in {\it loc. cit.} we know\footnote{This proposition proved that the symbol of a non zero germ of section of $\mathcal{I}_+$ vanishes on $\vert Z\vert$.} that the symbol $s(P)$ is in $p_*(I_Z)_{\sigma^0}$ thanks to the equality $I_Z = I_{\vert Z\vert}$ proved above. So there exists a germ $P_1$   in $ \mathcal{I}_{\sigma^0} \setminus \{0\}$ with symbol $s(P_1) = s(P)$. Then the order of $P - P_1$ is strictly less than $q$. But then $P - P_1$ is in $\mathcal{I}_{+, \sigma^0}$ with order strictly less than $q$ and then it is in $\mathcal{I}_{\sigma^0} $. Contradiction. So $\mathcal{I} = \mathcal{I}_+$.$\hfill \blacksquare$\\

\parag{Notations}
We note $D_N(m)$ the sub-sheaf of differential operators of order at most equal to $m$ and  $\mathcal{I}(m)$ the sub$-\mathcal{O}_N-$module generated in $\mathcal{I}$ by classes induced by differential operators of order at most equal to $m$.\\
Then we note $\mathcal{M}(m) := D_N(m)\big/\mathcal{I}(m)$.\\
For any non zero germ of section $P$ of $D_N$ we note $s(P)$ its symbol in $\mathcal{O}_N[\eta_1, \dots, \eta_k]$. For $P = 0$, let $s(P)$ be $0$.\\
Recall that we note $p : N\times \C^k \to N$ the projection.

\begin{lemma}\label{limite}
We have $\lim_{m \to \infty} \mathcal{M}(m)  \simeq  \mathcal{M} $, where the maps $\mathcal{M}(m) \to \mathcal{M}(m+1)$ are induced by the obvious inclusions
$$ D_N(m) \hookrightarrow D_N(m+1) \quad {\rm and} \quad \mathcal{I}(m) \hookrightarrow  \mathcal{I}(m+1) .$$
\end{lemma}

\parag{proof} Beware that the maps $\mathcal{M}(m) \to \mathcal{M}(m+1)$ are not "a priori" injective.\\
There is an obvious map $\lim_{m \to \infty} \mathcal{M}(m) \to \mathcal{M}$ which is clearly surjective. The point is to prove injectivity. Let $P$ be a non zero germ at some $\sigma \in N$ of order $m$
such that its image in $\mathcal{M}_\sigma$ is $0$. Then, by definition, there exists germs $B_h,  h \in [2, k]$ and $C_{p, q}, (p, q) \in [1, k]^2$ in $D_{N, \sigma}$ such that
$$ P = \sum_{h=2}^k B_h.\mathcal{T}^h + \sum_{p, q} C_{p, q}.A_{p, q} .$$
Let $r$ be the maximal order of the germs $B_h$ and $C_{p,q}$. Then the equality above shows that $P$ is in $\mathcal{I}(r+2)$. So the image of $P$ in $\lim_{m \to \infty} \mathcal{M}(m)$ is zero, as it is already $0$ in $\mathcal{M}(r+2)$.$\hfill \blacksquare$\\

\begin{lemma}\label{remplace} 
Let $P$ be a non zero germ of  section  of the sheaf $\mathcal{I}$. Assume that $P$ has order at most $1$. Then $P = 0$.
\end{lemma}

\parag{proof} 
 Let $P = a_0 + \sum_{h=1}^k a_h.\partial_h $. Recall that for each $h \in [1, k]$ and each $m \in \mathbb{N}$ we have (see proposition 5.2.1 in \cite{[B.19]}):
\begin{equation}
\partial_hN_m = (-1)^{h-1}.m.DN_{m-h} 
\end{equation}
where the polynomials
$$DN_m := \sum_{P_\sigma(x_j) = 0}  \frac{x_j^{m+k-1}}{P'_{\sigma}(x_j)} $$
vanish for $m \in [-k+1, -1]$ and $DN_0 = 1$.\\
Then  the equality $\mathcal{I}_+ = \mathcal{I}$ proved in the  proposition \ref{charact. 0} implies that for each integer $m $ we have
$$ a_0.N_m + \sum_{h=1}^k a_m.(-1)^{h-1}.m.DN_{m-h} = 0 ,\quad \forall m \in \mathbb{N} .$$
For $m = 0$ this gives $a_0 = 0$; if we have $a_0 = a_1 = \dots = a_p = 0$ for some $p \in [0, k-1]$  then $P[N_{p+1}] = 0$ gives $\sum_{h = p+1}^k a_h.(-1)^h.(p+1).DN_{p+1-h} =  a_{p+1}.(p+1).DN_0 = 0$ and then $a_{p+1} = 0$. 
So $P = 0$.$\hfill \blacksquare$\\

\begin{lemma}\label{computation}
Let $q \geq 2$ be an integer, $\alpha \in \mathbb{N}^k$ such that $\vert \alpha \vert = q-2$  and let $m$ be an integer in the interval $[2, k]$. The class induced by $\partial^\alpha.\mathcal{T}^m$ in $\mathcal{W}$ only depends on the integers $q$ and $r := w(\alpha) + m$. This class is given by the formula (with the convention $\sigma_0 \equiv 1$)
\begin{equation}
[\partial^\alpha.\mathcal{T}^m] = \sum_{h=0}^k \sigma_h.y_{q, r+h} + (q-1).y_{q, r},
\end{equation}
where $y_{q, r}$ is the class induced by $\partial^\gamma$ in $\mathcal{W}$ for any $\gamma \in \mathbb{N}^k$ such that $\vert \gamma \vert = q$ and $w(\gamma) = r$ (see lemma \ref{base}).\\
Let  $\lambda$ be a complex number and  let $\beta \in \mathbb{N}^k$ such that $\vert \beta\vert = q-1$ and $w(\beta) = r$. The class induced by $\partial^\beta.(U_0 - \lambda)$ in $\mathcal{W}$ only depends on $\lambda$ and on the integers $q$ and $r$. This class is given by
\begin{equation}
[\partial^\beta.(U_0 - \lambda)] = \sum_{h=1}^k  h.\sigma_h.y_{q, r+h} + (r - \lambda).y_{q-1, r} .
\end{equation}
Also the class induced by $\partial^\beta.U_{-1}$ in $\mathcal{W}$, again for $\vert\beta\vert = q-1$ and $w(\beta) = r$,  only depends on the integers $q$ and $r$. This class is given by
\begin{equation}
[\partial^\beta U_{-1} ]= \sum_{h=0}^k (k-h).\sigma_h.y_{q, r+h+1}  + (k.(q-1) - r).y_{q-1, r +1}
\end{equation}
where, for $r = k.(q-1)$, the last term in $(14)$ is equal to $0$ by convention.
\end{lemma}

\parag{proof} By definition $\mathcal{T}^m = \partial_1\partial_{m-1} + \sum_{h=1}^k \sigma_h.\partial_h\partial_m + \partial_m $  which implies 
$$ \partial^\alpha\mathcal{T}^m = \partial^\alpha \partial_1\partial_{m-1} + \sum_{h=1}^k \sigma_h.\partial_h\partial_m\partial^\alpha + (q-1).\partial^\alpha\partial_m $$
as we have  $\partial^\alpha.\sigma_h.\partial_h = \sigma_h.\partial_h\partial^\alpha + \alpha_h.\partial^\alpha$ for any $\alpha \in \mathbb{N}^k$ and any $h \in [1, k]$. Now formula $(12)$ follows from the lemma \ref{base}, proving our first assertion.\\
As   $U_0 := \sum_{h=1}^k h.\sigma_h.\partial_h$  we have  
$$ \partial^\beta.(U_0 - \lambda) = \sum_{h=1}^k h.\partial^\beta\sigma_h\partial_h - \lambda.\partial^\beta = \sum_{h=1}^k h.\sigma_h.\partial_h\partial^\beta + (w(\beta) - \lambda).\partial^\beta $$
which gives the formula $(13)$ using the  lemma \ref{base}, and this proves our second assertion.
The third one is analogous using the fact that $U_{-1} = \sum_{h=0}^{k-1} (k-h).\sigma_h.\partial_{h+1}$ with the convention $\sigma_0 \equiv 1$ and the equalities
$$ \sum_{h=0}^{k-1} (k-h).\beta_h = k.((q-1) - \beta_k) - (w(\beta) - k.\beta_k) = k.(q-1) - r $$
with the convention $\beta_0 = 0$.$\hfill \blacksquare$\\

\parag{Notations} \begin{enumerate}
\item Let $V_q \subset \mathcal{W}_q$ be the $\mathcal{O}_N-$sub-module with basis $y_{q, r}$ for $r \in [k.(q-1)+1, k.q]$. Remark that $V_0 = \mathcal{W}_0 = \mathcal{W}(0)= \mathcal{O}_N $ and $V_1 = \mathcal{W}_1 = \oplus_{h=1}^k \mathcal{O}_N.\partial_h$.
\item Let $L_q : \mathcal{W}_q \to \mathcal{M}(q)$ be the map induced by restriction to $\mathcal{W}_q$ of the quotient map $\mathcal{W}(q) \to \mathcal{M}(q)$ and $l_q : V_q \to \mathcal{M}(q)$ its restriction to $V_q$.
\end{enumerate}

 \begin{lemma}\label{missing}
 Fix an integer $q \geq 0$. Then for any $Y \in \mathcal{W}_q$ there exists $X \in V_q$ such that  $L_q(Y - X)$  is in $\mathcal{M}(q-1)$, with the convention $\mathcal{M}(-1) = \{0\}$.
 \end{lemma}
 
 \parag{proof} Remark that for $Y = y_{q, r}$ with $r \in [k.(q-1)+1, k.q]$ we may choose $X = Y$. So it is enough to prove the lemma for $Y$ in the sub-module with basis $y_{q, r}$ with $r \in [q, k.(q-1)]$.\\
 Note that for $q = 0$ and for $q = 1$ there is  nothing more to prove.\\
 For each $q \geq 2$ and $r \in [q, k.(q-1)]$ there exists $m \in [2, k]$ such that $r-m$ is in $[q-2, k.(q-2)]$ because the addition map $(s, m) \to s+m$  is surjective\footnote{For $r \in [q, k.(q-2)+2]$ take $s = r-2$ and $m = 2$, for $r = k.(q-2)+j$ with $j \in [2, k]$ take $s = r-j$ and $m = j$.} from $[q-2, k.(q-2)]\times [2, k]$ to $[q, k.(q-1)]$. So, there exists $\alpha \in \mathbb{N}^k$ such that $\vert \alpha\vert = q-2$ and $w(\alpha) = r-m$. Then
 $$ \partial^\alpha\mathcal{T}^m = y_{q, r} + (q-1).y_{q-1, r}+ \sum_{h=1}^k \sigma_h.y_{q, r+h} $$
 and the class induced by $y_{q, r}$ in $\mathcal{M}(q)$ is, modulo $\mathcal{M}(q-1)$, in the sub$-\mathcal{O}_N-$module of $\mathcal{M}(q)$ induced by the images of classes of $y_{q, r'}$ with $r' > r$. By a descending induction on $r \in [q, k.(q-1)]$ we see that, modulo $\mathcal{M}(q-1)$, the image of  $\mathcal{W}_q$ by $L_q$ is equal to $L_q(V_q)$. This implies our statement by induction on $q$.$\hfill \blacksquare$\\
 

\begin{prop}\label{description}
  For any $q \in \mathbb{N}$ there is  a natural isomorphism of  $\mathcal{O}_N-$modules
 \begin{equation}
\Lambda_q := \oplus_{p=0}^q l_p :  \oplus_{p=0}^{q}\  V_p  \longrightarrow  \mathcal{M}(q)
\end{equation}
which is compatible with the natural map $\oplus_{p=0}^{q}\  V_p \hookrightarrow \oplus_{p=0}^{q+1}\  V_p$ and  the natural map $\mathcal{M}(q) \rightarrow \mathcal{M}(q+1)$.
 \end{prop}
 
 \parag{Proof} For $q = 0$ we have $\mathcal{M}(0) = V_0 = \mathcal{O}_N.y_{0, 0}$ where $y_{0, 0} = 1$. So $\Lambda_0$ is an isomorphism. For $q = 1$, the lemma \ref{remplace} shows that the map $\Lambda_1$ is injective. As it is surjective  (we have $V_0 = \mathcal{W}_0$ and  $V_1 = \mathcal{W}_1$) the assertion is clear.\\
 Assume that we have proved that $\Lambda_{q-1}$ is an isomorphism of $\mathcal{O}_N-$modules for some $q \geq 2$ . We shall prove  that $\Lambda_q$ is also an isomorphism.\\
 Consider $Y := \sum_{p=0}^q  Y_p$ with $Y_p \in V_p$ for each $p \in [0, q]$, which is in the kernel of $\Lambda_q$. If $Y_q = 0$ the induction hypothesis allows to conclude that $Y = 0$.\\
 So assume that $Y_q \not= 0$. As $Y_q $ is induced\footnote{Note that for each $r \in [k.(q-1)+1, k.q], r = k.(q-1)+j$, then $y_{q, r}$ is induced in $W_q$ by $\partial_k^{q-1}.\partial_j$.} by a differential operator of the form $ \sum_{j=1}^k  b_j.\partial_k^{q-1}\partial_j$, with $b_j \in \mathcal{O}_N$ for $j \in [1, k]$, we may choose a differential operator $P \in \mathcal{I}$ of order $q$ which induces $Y$ such that its symbol is equal to $\eta_k^{q-1}.\sum_{j=1}^k  b_j.\eta_j$. This symbol vanishes on  $Z$, and, as $\eta_k$ does not vanish on any non empty open set on $Z$, we conclude that $\sum_{j=1}^k  b_j.\eta_j$ vanishes on $Z$. The injectivity of $\Lambda_1$ implies that $b_1 = \dots = b_k = 0$ showing that $Y_q = 0$ and this contradicts our hypothesis. So $\Lambda_q$ is injective.\\
 We have already noticed that the lemma \ref{missing} implies the surjectivity of $\Lambda_q$ for $q \geq 2$. So the proof is complete.$\hfill \blacksquare $
 
 \begin{cor}\label{non torsion 0}
 The $D_N-$module $\mathcal{M}$ has no $\mathcal{O}_N-$torsion.
 \end{cor}

\parag{Proof} This an  easy consequence of the previous proposition  giving that each $\mathcal{M}(q)$ is a free $\mathcal{O}_N-$module  because for any $\sigma \in N$ a non zero torsion germ in $\mathcal{M}_\sigma$  has to come from a non zero torsion element in $\mathcal{M}(q)_\sigma$ for some $q$ large enough (may be much more larger than the order of the germ in $D_{N, \sigma}$ inducing this class in $\mathcal{M}_\sigma$) thanks to the lemma \ref{limite}.$\hfill \blacksquare$
 
 .
 \subsection{On  quotients of $\mathcal{M}$}
 
 \subsubsection{A first result}
 
 We shall use the description of the characteristic variety of $\mathcal{M}$ to examine the holonomic quotients of $\mathcal{M}$  supported by an irreducible complex subset of $N$.
 
 \begin{prop}\label{co-torsion 1}
 Let $Q$ be a holonomic quotient of $\mathcal{M}$ which is supported by an analytic subset $S$ of $N$ with empty interior in $N$. Then $S$ is a hyper-surface and $S$ is contained in $\{\sigma_k = 0\}\cup \{\Delta(\sigma) = 0\}$.
 \end{prop}
 
 \parag{proof} Let $S_0$ be an irreducible component of $S$, the support of a holonomic quotient $Q$ of $\mathcal{M}$. Let $d \geq 1$ be the co-dimension of $S_0$. Then near the generic point in $S_0$ the  co-normal sheaf of $S_0$ is a rank $d$ vector bundle over $S_0$ which is contained in $Z$. As the fibres of $Z$ over $N$ have pure dimension $1$ we have $d \leq 1$ and  then $d = 1$ and $S_0$ is a hyper-surface in $N$. Then $S$ is also a hyper-surface in $N$.\\
 Let now $S_0$ be an irreducible component of $S$ which is not contained in $\{\Delta(\sigma) = 0 \}$. Then near the generic point in $S_0$ the quotient map $q : M \to N$ is an \'etale cover and this shows that $\mathcal{M}$ locally is isomorphic to the quotient of $D_{\C^k}$ by the let ideal with generators $\frac{\partial^2}{\partial z_i\partial z_j}$ for $i \not= j$ in $[1, k]$. So the characteristic variety of $\mathcal{M}$ is locally isomorphic to $C := \cup_{j=1}^k N\times \{\C.e_j\}$
 where $e_j$ is the $j-$th vector in the canonical basis of $\C^k$. If an irreducible hyper-surface has its co-normal bundle contained in $C$, it has to be equal to the co-normal of one of the hyperplanes $\{z_j = 0\}$. This means that $S$ is contained in $\{\sigma_k = 0 \}$.\\
  But any hyper-surface contained in $\{\Delta = 0 \}$ is equal to $\{\Delta = 0 \}$. So the only possible irreducible components of the support of $Q$ are $\{\sigma_k = 0 \}$ or $\{\Delta = 0 \}$.$\hfill \blacksquare$\\
 
 We shall use the following immediate corollary of this proposition:
 
 \begin{cor}\label{co-torsion 2}
 Let $Q$ be a coherent holonomic quotient of $\mathcal{M}$ which is supported in a closed analytic subset $S$ in $N$ with empty interior in $N$. If $Q$ vanishes near the generic points of $\{\sigma_k = 0 \} \cup \{\Delta = 0 \}$, then $Q = \{0\}$.$\hfill \blacksquare$\\
 \end{cor}
 
 \subsubsection{A local chart}
 Let $k \geq 2$. We shall study the $D_N-$module $\mathcal{M}$  near the generic point of the hyper-surface $\{\Delta = 0\}$ in $N$.\\

Let $z_1^0, z_3^0, \dots, z_k^0$ be $(k-1)$ distinct points in  $\C$ and let $r > 0$ a real number small enough in order that the discs $D_1, D_3, \dots, D_k$ with respective centers $z_1^0, z_3^0, \dots, z_k^0$ and radius $r$  are two by two disjoint. Let $\mathcal{U}_0 := D_1\times D_1\times  \prod_{j=3}^k D_j$ and $\mathcal{V}$ (equal to $D_1\times D_2$ for $k = 2$) the image of $\mathcal{U}_0$ by the quotient map $q : \C^k \to \C^k\big/\mathfrak{S}_k = N$. Note $\mathcal{U} := q^{-1}(\mathcal{V})$. Then $q$ induces an isomorphism of $\mathcal{U}_0\big/\mathfrak{S}_2$ onto $\mathcal{V}$.\\
Remark that for each $\sigma \in \mathcal{V}$ we have exactly two roots $z_1(\sigma), z_2(\sigma)$ distinct or not which are in $D_1$ and for each $j \in [3, k]$ we have exactly one (simple) root $z_j(\sigma)$ in $D_j$.
We have the following holomorphic map on $\mathcal{V}$:
\begin{enumerate}
\item The map $\tau = (\tau_1, \tau_2) : \mathcal{V} \to \C^2$ given by $\tau_1(\sigma) := z_1(\sigma) + z_2(\sigma)$ and  \\
$\tau_2 := z_1(\sigma).z_2(\sigma)$ where $z_1(\sigma)$ and $z_2(\sigma)$ are the roots of $P_\sigma$ which are in $D_1$.
\item For each $j \in [3, k]$ the map $z_j : \mathcal{V} \to D_j$ given by the unique  (simple) root of $P_\sigma$ in $D_j$.
\end{enumerate}
To be completely clear, these holomorphic maps are defined on $\mathcal{V}$ by the following integral formulas:
\begin{align*}
&  \tau_1(\sigma) := \frac{1}{2i\pi} \int_{\partial D_1} \frac{\zeta.P'_\sigma(\zeta).d\zeta}{P_\sigma(\zeta)} \\
&  2\tau_2(\sigma) = \tau^2_1 - \nu_2(\sigma) \quad {\rm where} \quad \nu_2(\sigma) := \frac{1}{2i\pi} \int_{\partial D_1} \frac{\zeta^2.P'_\sigma(\zeta).d\zeta}{P_\sigma(\zeta)} \\
& {\rm and \ for} \ j \in [3, k] \quad  z_j(\sigma) := \frac{1}{2i\pi} \int_{\partial D_j} \frac{\zeta.P'_\sigma(\zeta).d\zeta}{P_\sigma(\zeta)} 
\end{align*}

The following lemma is obvious:

\begin{lemma}\label{facile}
The holomorphic map $\Phi : \mathcal{V} \to \C^2\times \prod_{j=3}^k D_j$ given by $(\tau_1, \tau_2, z_3, \dots, z_k)$, is an isomorphism of $\mathcal{V}$ onto the open set $\mathcal{V}_1 := D_{1,2}\times \prod_{j=3}^k D_j$ where $D_{1,2} := (D_1\times D_1)\big/\mathfrak{S}_2$ is the image of $D_1\times D_1$ by the quotient map by the action of $\mathfrak{S}_2$.$\hfill \blacksquare$\\
\end{lemma}

In the sequel, we shall use the coordinate system on $\mathcal{V}$  given by the holomorphic functions  $\tau_1, \tau_2, z_3, \dots, z_k$ on $\mathcal{V}$.

Define on $\mathcal{V} \simeq \mathcal{V}_1=   D_{1, 2}\times \prod_{j=3}^k D_j$ the following partial differential operators in the coordinate system described above:
\begin{itemize}
\item $T^2 := \partial_{\tau_1}^2 + \tau_1\partial_{\tau_1}\partial_{\tau_2} + \tau_2\partial_{\tau_2}^2 + \partial_{\tau_2}$ 
\item $ B_{i,j} := \partial^2\big/\partial z_i\partial z_j$ for $3 \leq i < j \leq k$
\item $ C_{1,j} := \partial^2\big/\partial \tau_1\partial z_j$ for $j \in [3, k]$ 
\item $ C_{2, j} := \partial^2\big/\partial \tau_2\partial z_j$ for $j \in [3, k]$ 
\item $ V_0 := \tau_1\partial_{\tau_1} + 2\tau_2\partial_{\tau_2} + \sum_{j=3}^k  z_j\partial_{z_j}  $
\item $V_{-1} := 2\partial_{\tau_1} + \tau_1.\partial_{\tau_2} + \sum_{j=3}^k \partial_{z_j} .$
\end{itemize}

\begin{prop}\label{the point} The isomorphism of change of coordinates $\Phi$ on $\mathcal{V}$  given by $\sigma \mapsto (\tau_1, \tau_2, z_3, \dots, z_k)$ has the following properties:
\begin{enumerate}[(i)]
\item The image of ideal $\mathcal{I}$ of $D_N$ restricted to $\mathcal{V}$ by the isomorphism $\Phi $ is  the left ideal generated by $T^2$, $B_{i,j}$, $C_{1,j}$ and $C_{2, j}$ in $D_{\mathcal{V}_1}$.
\item The vector field $U_0$ is sent to $V_0$ and the vector field $U_{-1}$ is sent to $V_{-1}$ by this isomorphism.
\end{enumerate}
\end{prop}

\parag{Proof} We shall use the local version of the theorem 5.1.1   in \cite{[B.19]}  given in the proposition \ref{charact. 0}.\\
For $\Omega \subset \mathcal{V}$ it is easy to see that the Frechet space of trace functions admits as a dense subset the finite $\C-$linear combinations of the Newton functions $\nu_m, m \in \mathbb{N}$ of $z_1(\sigma), z_2(\sigma)$ and of the functions $z_j^m(\sigma), m \in \mathbb{N}$ for each $j \in [3, k]$. From the case $k = 2$ for which the left ideal $\mathcal{I}$ is generated by $T^2$ and the fact that each $B_{i,j}, C_{1,j}$ and $C_{2,j}$ kill each $\nu_m$ and each $z_j^m$, we conclude that $\mathcal{I} $ contains the left ideal generated by $T^2$, the $B_{i,j}$, the $C_{1,j}$ and the $C_{2, j}$. \\
Conversely, if $P$ is in $\mathcal{I} $ it has to kill any $\nu_m$  and each $z_j^m, \forall j \in [3, k]$. So $P$ has no order $0$ term. Modulo the ideal generated by the $B_{i,j}$, the $C_{1,j}$ and the $C_{2, j}$ we may assume that we can write
$$ P = P_0 +\sum_{m=1}^N \sum_{j=3}^k g_{j,m}.\partial_{z_j}^m $$
where $P_0$ is a differential operator in $\tau_1, \tau_2$ with no order $0$ term, and with holomorphic dependence in $z_3, \dots,z_k$ (but no derivation in these variables) and where $g_{j, m}$ are holomorphic  functions on $\mathcal{V}$. Applying $P$ to $z_j^N$ gives that $N! g_{j,N} = 0$. So we see that all $g_{j,m}$ must vanish. Then $P = P_0$ and $P(\nu_m) = 0$ implies that $P_0$ is in the left ideal generated by $T^2$ in the $\mathcal{O}_{\mathcal{V}_1}-$algebra generated by $\frac{\partial}{\partial \tau_1}$ and $ \frac{\partial}{\partial \tau_2}$. Then $P_0$ and also $P$ is in our ideal and  $(i)$ is proved.\\
 The verification of $(ii)$ is easy and left to the reader.$\hfill \blacksquare$\\

\begin{lemma}\label{case k = 2}
For $k=2$ we have for each $n \in \mathbb{N}^*$
\begin{equation}
 (T^2 - 2n.\partial_2).\Delta^n = \Delta^n.T^2 + 2n.(2n + 1).\Delta^{n-1}. 
 \end{equation}
\end{lemma}

\parag{Proof} Recall that we have   $E := \sigma_1.\partial_1 + \sigma_2.\partial_2$ and $T^2 = \partial_1^2 + \partial_2.E$ and  that  \\
$\Delta = \sigma_1^2 - 4\sigma_2$. So
\begin{align*}
& \partial_1.\Delta = \Delta.\partial_1 + 2\sigma_1 \\
& [\sigma_1.\partial_1, \Delta] = 2\sigma_1^2 \\
& \partial_2.\Delta = \Delta.\partial_2 - 4 \\
& [\sigma_2.\partial_2, \Delta] = -4\sigma_2 \\
& [E, \Delta] = 2\sigma_1^2 - 4\sigma_2 \\
& [\partial_2.E, \Delta] = \partial_2.(\Delta.E + 2\sigma_1^2 - 4\sigma_2) - (\partial_2.\Delta +4).E = 2\Delta.\partial_2 - 4\sigma_1.\partial_1 - 4 \\
&  \partial_1^2.\Delta = \partial_1.(\Delta.\partial_1 + 2\sigma_1) =   ( \Delta.\partial_1 + 2\sigma_1).\partial_1 + 2\sigma_1.\partial_1 + 2 =  \Delta.\partial_1^2 + 4\sigma_1.\partial_1 + 2 \\
& [T^2,\Delta] = [\partial_1^2, \Delta] + [\partial_2.E, \Delta] \\
& T^2.\Delta = \Delta.T^2 + 4\sigma_1.\partial_1 + 2 + 2\Delta.\partial_2 - 4\sigma_1.\partial_1 - 4 = \Delta.T^2 + 2\partial_2.\Delta + 8 +2 - 4   \\
& (T^2 - 2\partial_2).\Delta = \Delta.T^2 + 6
\end{align*}
which proves $(16)$ for $n = 1$.\\
Assume now that we have proved the formula $(16)$ for $n \geq 1$.  Then we have, using that $\Delta^n.\partial_2 = \partial_2.\Delta^n + 4n.\Delta^{n-1}$
\begin{align*}
& (T^2 - 2n.\partial_2).\Delta^{n+1} = \Delta^n.T^2.\Delta + 2n.(2n + 1).\Delta^n \\
&  (T^2 - 2n.\partial_2).\Delta^{n+1} =  \Delta^{n+1}.T^2 + \Delta^n.2\partial_2.\Delta  + 6\Delta^n + 2n.(2n + 1).\Delta^n \\
& (T^2 - 2n.\partial_2).\Delta^{n+1} =  \Delta^{n+1}.T^2 + 2\partial_2.\Delta^{n+1} + 8n.\Delta^n + 6\Delta^n + 2n.(2n + 1).\Delta^n \\
& (T^2 - 2(n+1).\partial_2).\Delta^{n+1} = \Delta^{n+1}.T^2 + 2(n+1).(2n+3).\Delta^n
\end{align*}
because $2n.(2n + 1) + 8n + 6 = 2.(n+1).(2n+3)$.$\hfill \blacksquare$\\

\begin{thm}\label{co-torsion delta}
Let $Q$ be a coherent $D_N-$module which is a quotient of $\mathcal{M}$ and which is supported by $\{\Delta = 0\}$. Then $Q = \{0\}$. Moreover any holonomic quotient of $\mathcal{M}$ has no $\Delta-$torsion.
\end{thm}

\parag{proof} Thanks to the corollary \ref{co-torsion 2} it is enough to prove that  such a quotient $Q$ is zero near the generic points of $\{\Delta = 0\}$. So assume that $Q$ is such a non zero quotient. Using now the result of the proposition \ref{the point}, the $D-$module $Q$ is given near the generic points of $\{\Delta = 0\}$ by the quotient of $D$ by a left ideal $\mathcal{K}$ which contains $T^2$. Then there exists an integer $n > 0$ such that $\Delta^n$ belongs to $\mathcal{K}$. Then the lemma \ref{case k = 2} implies that $\mathcal{K}$ contains $\Delta^{n-1}$. Then by a descending induction on $n$ we obtain that $1$ is in $\mathcal{K}$ and this contradicts the non vanishing of $Q$.\\
The characteristic variety of a holonomic quotient of $\mathcal{M}$ which is supported in codimension $\geq 1$ in $N$  is contained in the characterisc variety of $\mathcal{M}$ so is contained in the union of  $N\times \{0\}$ with the co-normal to $\{\sigma_k = 0\}$ and $\{\Delta = 0\}$ thanks to the proposition \ref{co-torsion 1}. But the lemma \ref{case k = 2} implies that near the generic point of $\{\Delta = 0\}$ a torsion element in such a quotient vanishes. So the torsion submodule of a holonomic  quotient of $\mathcal{M}$ cannot have the co-normal of $\{\Delta = 0\}$ in its characteristic variety. Then such a quotient has no $\Delta-$torsion.$\hfill \blacksquare$\\

\subsection{Action of $sl_2(\C)$ on $\mathcal{M}$}

Let $B$ be the sub-$\C-$algebra of the Weyl algebra $\C[\sigma]\langle \partial \rangle$ generated by the vector fields $U_p, p \geq -1$, where $U_p$ is the vector field on $N$ defined as the image by the differential $T_q$ of the quotient map 
$$ q : M := \C^k \longrightarrow N := \C^k\big/\mathfrak{S}_k \simeq \C^k $$
of the vector field \   $\sum_{j = 1}^k z_j^{p+1}.\frac{\partial}{\partial z_j} $. 

\begin{thm}\label{action}
The right action of $B$ on $D_N$ induces a morphism of algebras between $B$ and the algebra of left $D_N-$linear endomorphisms of $\mathcal{M}$. Moreover, we have  the identities
\begin{equation}
U_p.U_q  - U_q.U_p = (q-p).U_{p+q}  \quad {\rm modulo} \quad \mathcal{I} 
\end{equation}
which imply that in the action of $B$ on $\mathcal{M}$ we have $[U_p, U_q] = (q-p).U_{p+q}$
\end{thm}

\parag{Proof} It will be enough to show that for each integer  $p \geq -1$ we have the inclusion $\mathcal{I}.U_p \subset \mathcal{I}$. If it is not difficult to prove such an inclusion for $p = -1$ or $p= 0$ by a direct computation of the commutators of $U_p$ with the generators of $\mathcal{I}$, it seems rather difficult to do it for $p$ large because the coordinates of $U_p$ in the $\C[\sigma]$ basis $\partial_1, \dots, \partial_k$ of the polynomial vector fields on $N$ seems more and more complicated. So we shall use the local version of the theorem 5.1.1  in \cite{[B.19]} given in the proposition \ref{charact. 0}.\\
Let $P \in \mathcal{I}$, $p \geq -1 $ an integer  and $m \in \mathbb{N}$. Using the  formula $U_p[N_m] = m.N_{m+p}$ which is easy to verify on $M$, we get
$$P[U_p[N_m]] = P[m.N_{m+p}] = 0 $$
when $P$ annihilates any Newton polynomial. Then $P.U_p$ also annihilates any Newton polynomial and thanks to the proposition \ref{charact. 0}  we conclude that $P.U_p$ belongs to $ \mathcal{I}$ proving the first assertion.\\
We shall prove the second assertion using the same result: to prove the formula $(17)$ it is enough to check that for each $m \in \mathbb{N}$ we have
$$ \big(U_p.U_q - U_q.U_p - (q-p).U_{p+q}\big)[N_m] = 0 .$$
But 
\begin{align*}
& U_p[U_q[N_m]] = U_p[m.N_{q+m}]= m.(m+q).N_{m+p+q} \\
& U_q[U_p[N_m]] = U_q[m.N_{p+m}] = m.(m+p).N_{m+p+q} \\
& U_{p+q}[N_m] = m.N_{m+p+q}
\end{align*}
concluding the proof.$\hfill \blacksquare$\\

The commutation relations $[U_0, U_{-1}] = -U_{-1}, [U_0, U_1] = U_1$ and $[U_1, U_{-1}] = 2U_0$ which are easy to check in $M$ show that the Lie algebra $\mathcal{L}$ generated by  the $U_p$ (with the commutators given by  the formula $(17)$) contains a sub-Lie algebra isomorphic to $ sl_2(\C)$. The formula $(17)$ shows that $\mathcal{L}$ acts on $\mathcal{M}$ and  then induces a structure of $sl_2(\C)-$module on $\mathcal{M}$.

\section{The $D_N-$modules $\mathcal{N}_\lambda$}

\subsection{Homothetie and translation}

\parag{Notations} Let $\lambda$ a complex number. We define the left ideal
 $$\mathcal{J}_\lambda := \mathcal{I} + D_N.(U_0 - \lambda)$$
 in $D_N$ and let $\mathcal{N}_\lambda$ be the quotient $D_N\big/\mathcal{J}_\lambda$. We shall denote by $q_\lambda : \mathcal{M} \to \mathcal{N}_\lambda$ the quotient map.\\
 We shall denote respectively by $\mathscr{H}_{\lambda}$ and $\mathscr{T}$ the endomorphisms of left $D_N-$modules on $\mathcal{M}$ induced  respectively by the right multiplications by $U_{0} - \lambda$ and $U_{-1} $ (see the theorem \ref{action} ).
They satisfy the commutation relation (see {\it loc. cit.})
  $$\mathscr{H}_{\lambda}.\mathscr{T}- \mathscr{T}.\mathscr{H}_{\lambda} = - \mathscr{T}$$
  for each $\lambda \in \C$ and  $\mathcal{N}_{\lambda}$ is, by definition, the co-kernel of $\mathscr{H}_{\lambda}$.\\
 As $\mathcal{I}.U_{-1} \subset \mathcal{I}$, writing this relation in the form $\mathscr{H}_{\lambda-1}.\mathscr{T} = \mathscr{T}.\mathscr{H}_{\lambda}$ we see that the right multiplication by $U_{-1}$ induces a left $D_N-$modules morphism $\mathscr{T}_\lambda : \mathcal{N}_{\lambda-1} \to \mathcal{N}_{\lambda}$ for each $\lambda$.

\begin{prop}\label{fond.1}
For each $\lambda \in \C$ we have an exact sequence of left $D_N-$modules on $N$
\begin{equation}
0 \to \mathcal{M} \overset{\mathscr{H}_{\lambda}}{\longrightarrow} \mathcal{M} \overset{q_{\lambda}}{\longrightarrow} \mathcal{N}_{\lambda} \to 0 
\end{equation}
where $q_{\lambda}$ is the obvious quotient map.
\end{prop}

\parag{Proof}  The  quotient map $q_{\lambda}$ is surjective by definition, so the point is to prove that the kernel of $q_{\lambda}$ is isomorphic to $\mathcal{M}$.\\ 
This kernel is obviously given by
\begin{equation}
\mathcal{J}_{\lambda}\big/ \mathcal{I} \simeq \big(\mathcal{I} + D_{N}.(U_0 - \lambda)\big)\big/\mathcal{I} \simeq  D_{N}.(U_0 - \lambda)\big/\mathcal{I}\cap D_{N}(U_{0}- \lambda)
\end{equation}
The proof  will be an easy consequence of  the following lemma.

\begin{lemma}\label{fond.0}
Let $P$ be a germ in $D_{N, \sigma}$ for some $\sigma \in N$ such that $P.(U_0 - \lambda) $ is in 
$\mathcal{I}_{\sigma}$. Then $P$ is in $\mathcal{I}_\sigma$. So \  $\mathcal{I}_\sigma \cap D_{N, \sigma}.(U_0 -\lambda) = \mathcal{I}_\sigma.(U_0 -\lambda)$.
\end{lemma}

\parag{Proof} Assume that the lemma is wrong. Then let $P_{0}$  in $ D_{N, \sigma}$ having minimal order among germs $P$ in $ D_{N, \sigma}$ satisfying the following properties
\begin{enumerate}
\item $P.(U_0 -\lambda)$ is in $\mathcal{I}_{\sigma}\cap D_{N, \sigma}(U_{0}- \lambda)$ ;
\item $P$ is not in $ \mathcal{I}_{\sigma}$.
\end{enumerate}
 Let $\pi$ be the symbol of $P_{0}$  and let $g$ be the symbol of $U_{0}$. We have $\pi.g \in p_* (I_Z)_\sigma$. But we know that $g$ does not vanish on any non empty open set of $Z$ because $\{g = 0\} \cap Z$ has pure co-dimension $1$ in $Z$ (see lemma \ref{remplace} above). Then $\pi$ vanishes on $(V \times \C^k) \cap Z$ where $V$ is a neighborhood  of $\sigma$ in $N$ and, as we have proved that $Z$ is reduced and is the characteristic cycle of $\mathcal{M}$, their exists a germ $P_{1}$ in $\mathcal{I}_{\sigma} $ with symbol  equal to $\pi$. Then $(P_{0}- P_{1}).(U_{0}- \lambda)$ satisfies again the properties 1. and  2. and is of order strictly less than the order of $P_{0}$. So  $P_0 - P_{1}$ is in $\mathcal{I}_{\sigma}$ and this contradicts the fact that we assumed that $P_{0}$ is not in $\mathcal{I}_{\sigma}$.$\hfill \blacksquare$

\parag{End of proof of \ref{fond.1}} The previous lemma shows that for each $\lambda \in \C$
$$ \mathcal{I} \cap D_{N}.(U_{0} - \lambda) = \mathcal{I}.(U_{0} - \lambda).$$

So the right multiplication by $U_{0}- \lambda$ induces an isomorphism of left $D_N-$modules
$$ \mathcal{M} \to D_{N}.(U_{0}- \lambda)\big/\mathcal{I}.(U_{0}- \lambda) $$
and the kernel of $q_{\lambda}$ is isomorphic to $\mathcal{M}$ by the inverse of this isomorphism.$\hfill \blacksquare$\\

\begin{defn}\label{etoiles}
Define the $D_N-$module $\tilde{\mathcal{M}}$ as the quotient $D_N\big/(\mathcal{I} + D_N.U_{-1})$. 
 For each $\lambda \in \C$ then define $\mathscr{T}_\lambda : \mathcal{N}_\lambda \to \mathcal{N}_{\lambda+1}$ as the $D_N-$linear map induced by $\mathscr{T}$.
\end{defn}

\begin{lemma}\label{11/06}
For each $\lambda \in \C$ the co-kernel of the $D_N-$linear map $\mathscr{T}_{\lambda+1}$
 is naturally isomorphic to  the co-kernel of  the $D_N-$linear  map $\tilde{\mathscr{H}}_{\lambda+1}  : \tilde{\mathcal{M}} \to \tilde{\mathcal{M}}$ induced by $\mathscr{H}_{\lambda+1} $ and 
there is also  a natural isomorphism of $D_N-$modules between the kernels of $ \tilde{\mathscr{H}}_{\lambda+1}$ and $\mathscr{T}_{\lambda+1}$.
\end{lemma}

\parag{proof} Consider the commutative diagram of left $D_N-$modules with exact lines and columns:
$$ \xymatrix{ \quad & \quad & \quad & \quad & 0 \ar[d] & \quad \\ \quad & \quad & \quad & \quad &Ker(\tilde{\mathscr{H}}_{\lambda+1}) \ar[d] & \quad \\ \quad & \quad & \mathcal{M} \ar[d]^{\mathscr{H}_\lambda} \ar[r]^{\mathscr{T}} & \mathcal{M} \ar[d]^{\mathscr{H}_{\lambda+1}} \ar[r] & \tilde{\mathcal{M}} \ar[d]^{\tilde{\mathscr{H}}_{\lambda+1}} \ar[r] & 0 \\ \quad & \quad & \mathcal{M} \ar[r]^{\mathscr{T}}  \ar[d]  & \mathcal{M} \ar[d] \ar[r] & \tilde{\mathcal{M}} \ar[r] \ar[d] & 0 \\ 
0 \ar[r] & Ker(\mathscr{T}_{\lambda+1}) \ar[r]  \ar@/^5pc/[rrruuu]  & \mathcal{N}_\lambda \ar[r]^{\mathscr{T}_{\lambda+1}} \ar[d] & \mathcal{N}_{\lambda+1} \ar[d] \ar[r] &  \mathcal{N}^\square_{\lambda+1} \ar[r] & 0 \\ \quad & \quad & 0 & 0 & \quad & \quad } $$
where $\mathcal{N}^\square_{\lambda+1}$ is, by definition the co-kernel of $\mathscr{T}_{\lambda+1} : \mathcal{N}_\lambda \to \mathcal{N}_{\lambda+1}$. By a  simple diagram chasing it is easy to see that $\mathcal{N}^\square_{\lambda+1}$  is  also the co-kernel of $\tilde{\mathscr{H}}_{\lambda+1} : \tilde{\mathcal{M}} \to \tilde{\mathcal{M}}$.\\
A elementary diagram chasing gives also the  isomorphism between  kernels of $ \tilde{\mathscr{H}}_{\lambda+1}$ and $\mathscr{T}_{\lambda+1}$. $\hfill \blacksquare$\\

We shall prove now that for $\lambda \not= 0, 1$ the map $\mathscr{T}_\lambda$ is an isomorphism of left $D_N-$modules. This implies $\mathcal{N}^\square_\lambda = \{0\}$ for $\lambda \not= 0, 1$.

\begin{lemma}\label{isom 2} Let $\mathscr{G}_\lambda : \mathcal{N}_{\lambda+1} \to \mathcal{N}_\lambda$ the left $D_N-$linear map given by right multiplication by $U_1$. Then we have for each $\lambda \in \C$
\begin{align*}
&  \mathscr{T}_\lambda\circ \mathscr{G}_{\lambda-1} = \lambda.(\lambda - 1) \quad {\rm on} \quad \mathcal{N}_{\lambda} \tag{A} \\
& \mathscr{G}_{\lambda-1}\circ \mathscr{T}_\lambda = \lambda.(\lambda - 1)  \quad {\rm on} \quad \mathcal{N}_{\lambda-1} \tag{B}
\end{align*}
So for $\lambda \not= 0, 1$ the  left $D_N-$linear map $\mathscr{T}_\lambda$ is an isomorphism.
\end{lemma}

\parag{proof} First recall that the computations made at the end of the proof of the theorem \ref{action}  give, for $p = \pm 1$ and $q = -p$, the equalities
\begin{equation}
 U_1.U_{-1} = U_0.(U_0 - 1) \ {\rm modulo} \ \mathcal{I} \quad  {\rm and} \quad U_{-1}.U_1 = U_0.(U_0 + 1)  \ {\rm modulo} \ \mathcal{I} 
 \end{equation}
 These give $(A)$ and $(B)$ and the conclusion follows.$\hfill \blacksquare$\\
 
 The following important result shows that adding to the ideal $\mathcal{A}$ the invariance by translation and the homogeneity $1$, that is to say considering the left ideal in $D_N$:
   $\mathcal{A} + D_N.U_{-1} + D_N.(U_0 - 1)$, we recover the ideal $\mathcal{J}_1 + D_N.U_{-1}$ and $D_N\big/(\mathcal{J}_1 + D_N.U_{-1})$   is the co-kernel of the map $\mathcal{T}_1 : \mathcal{N}_0 \to \mathcal{N}_1$. So, as a corollary,  we shall obtain the equality
 \begin{equation}
 \mathcal{N}^\square_1 =  D_N\big/ \mathcal{A} + D_N.U_{-1} + D_N.(U_0 - 1).
 \end{equation}
 \smallskip

\begin{prop}\label{magic}
For $h \in [2, k]$ we have the equality
\begin{equation*}
\partial_h.(U_0 - 1) + \partial_{h-1}.U_{-1} = k.\mathcal{T}^h +  \sum_{q=1}^{k-1} (k-q).\sigma_q.A_{h-1, q+1} \tag{$E_h$}
\end{equation*}
and for $h = 1$ the equality
\begin{equation*}
-\partial_1.(U_0 - 1) +  E.U_{-1} =  \sum_{q=1}^{k-1} (k-q).\sigma_{q}.\mathcal{T}^{q+1}  . \tag{$E_1$}
\end{equation*}
\end{prop}

\parag{Proof} Recall first that, if we put $E := \sum_{h=1}^k h.\sigma_h.\partial_h$ then for any $m \in [2, k]$ we have
 $$\mathcal{T}^m = \partial_1\partial_{m-1} + \partial_m.E = \partial_1\partial_{m-1} + E.\partial_m + \partial_m .$$
 For $h \in [2, k]$ we have
\begin{align*}
& \partial_h.U_0 + \partial_{h-1}.U_{-1} = \sum_{q=1}^k q.\sigma_q.\partial_q\partial_h + h.\partial_h + \sum_{q=0}^{k-1} (k-q).\sigma_q.\partial_{q+1}\partial_{h-1} + (k-h+1).\partial_h \\
& \quad =  \sum_{q=1}^k q.\sigma_q.\partial_q\partial_h + \sum_{q=1}^{k} (k-q).\sigma_{q}.A_{h-1, q+1}  + k.\partial_1\partial_{h-1} + \sum_{q=1}^{k} (k-q).\sigma_q.\partial_q\partial_h  + (k+1).\partial_h \\
& \quad = k.E.\partial_h + k.\partial_1\partial_{h-1} + k.\partial_h + \partial_h + \sum_{q=1}^{k-1} (k-q).\sigma_q.A_{h-1, q+1}\\
& \quad =  k.\mathcal{T}^h +  \partial_h + \sum_{q=1}^{k-1} (k-q).\sigma_{q}.A_{h-1, q+1}
\end{align*}
which is $(E_h)$.\\
For $h = 1$ let us compute $ \sum_{q=1}^{k-1} (k-q).\sigma_{q}.\mathcal{T}^{q+1} $:
\begin{align*}
&  \sum_{q=1}^{k-1} (k-q).\sigma_{q}.\mathcal{T}^{q+1}  =  \sum_{q=1}^{k-1} (k-q).\sigma_{q}.\big(\partial_1\partial_{q} + \partial_{q+1}.E \big) \\
&  \sum_{q=1}^{k-1} (k-q).\sigma_{q}.\mathcal{T}^{q+1}  =  \big(\sum_{q=1}^{k-1} (k-q).\sigma_q.\partial_q\big).\partial_1 +  \sum_{q=1}^{k-1} (k-q).\sigma_{q}.\partial_{q+1}.E    \\
&  \sum_{q=1}^{k-1} (k-q).\sigma_{q}.\mathcal{T}^{q+1}   = \big(\sum_{q=1}^{k-1} (k-q).\sigma_q.\partial_q\big).\partial_1 + (U_{-1} - k.\partial_1).E \\
& \quad \quad = k.(E - \sigma_k.\partial_k).\partial_1 - (U_0 - k.\sigma_k.\partial_k).\partial_1 + (U_{-1} - k.\partial_1).E \\
& \quad \quad =  k.E.\partial_1 - k.\partial_1.E - U_0.\partial_1 + U_{-1}.E \\
& \quad \quad  = E.U_{-1} - \partial_1.(U_0 - 1)
\end{align*}
using the commutation relations $[U_{-1}, E] = k.\partial_1$, $[E, \partial_1] = -\partial_1$ and $[U_0, \partial_1] = -\partial_1$. So we  obtain the equality $(E_1)$.$\hfill \blacksquare$\\

\parag{Remarks} 
\begin{enumerate}
\item An interesting way to look at these relations is to compare them with the minors of the $(k+1, 2)$ matrix\\ 

\begin{equation*}
\begin{pmatrix}  \mathcal{T} &  - (\mathscr{H}_0 - 1) \\  \partial_1 & -\mathscr{E} \\  \partial_2 & \partial_1 \\ . & .\\ . & .\\ . & . \\ \partial_k &  \partial_{k-1}  \end{pmatrix}
\end{equation*} 

where, by definition, $\mathscr{E}$ is the right product by $E$ in the Weyl algebra $\C[\sigma]\langle \eta \rangle$. The relations $(E_h), h \in [1, k]$ may also be seen as the fact  that\\

\begin{equation*}
 \begin{pmatrix} \partial_1 \\ \partial_2 \\ . \\ . \\ . \\  \partial_k \end{pmatrix}.(U_0 - 1) + \begin{pmatrix} -E \\  \partial_1 \\ . \\ . \\ . \\ \partial_{k-1} \end{pmatrix}.U_{-1} 
 \end{equation*}
 
is a global section of $\mathcal{I}^k \subset D_N^k$.
\item Let $g$ and $\gamma$ be the symbols of $U_0$ and $U_{-1}$ respectively. Looking at the symbols in formulas $(E_h), h \in [1, k]$,  we obtain (recall that $l_\sigma(\eta)$ is the symbol of $E$):
\begin{equation*}
 g(\sigma, \eta).\begin{pmatrix} \eta_1\\ \eta_2 \\ \dots \\ \dots \\ \eta_k \end{pmatrix} + \gamma(\sigma, \eta).\begin{pmatrix} -l_\sigma(\eta) \\ \eta_1 \\ \eta_2 \\ \dots \\ \eta_{k-1} \end{pmatrix} = 0 \quad {\rm on} \quad Z  \tag{F}
 \end{equation*}
  \item  For $\lambda \not= 0$ the sheaf of solutions\footnote{We mean here $Sol^0(\mathcal{N}_\lambda) := \underline{Hom}_{D_N}(\mathcal{N}_\lambda, \mathcal{O}_N)$.} of $\mathcal{N}_\lambda$ near the generic point in $N$ is the rank $k$ local system with basis $z_j^\lambda$. This is consequence of the the fact that any local trace function\footnote{See \cite{[B.19]}.} $F$ which satisfies $(U_0 - \lambda)[F] = 0$ is the trace of a homogeneous function of degree $\lambda$. 
\item Thanks to the lemma \ref{isom 2}, the map induced by  $\mathcal{T}_\lambda$ on solutions 
$$ Sol^0(\mathcal{N}_\lambda) \to Sol^0(\mathcal{N}_{\lambda-1}) $$
 sends $z_j^\lambda$ to $\lambda.z_j^{\lambda-1}$. This  is clearly an isomorphism for $\lambda \not= 0, 1$.
\end{enumerate}

\subsection{Characteristic variety}

\begin{prop}\label{Nouvelle version}
Let $\mathcal{I}$ be a coherent left ideal in $D_N$ such that its characteristic ideal $I_Z$ is the reduced ideal of an analytic subset $Z$ in $N \times \C^k$. Let $U \in \Gamma(N, D_N)$ be a differential operator of order $q$ such that its symbol $u$ does not vanish on any non empty open set in $Z$. Assume that $\mathcal{I}.U \subset \mathcal{I}$. Then the characteristic ideal of $\mathcal{I} + D_N.U$ is equal to $I_Z + \mathcal{O}_{N\times \C^k}.u$.\\
Moreover, for any $\sigma \in N$ and any germ at $\sigma$ of order $q+r$:  $Q  = P + B.U$ where $P \in \mathcal{I}_\sigma$ and $B \in D_{N, \sigma}$ there exists $P_1 \in \mathcal{I}_\sigma$ of order at most $q+r$ and $B_1 \in D_{N, \sigma}$ of order at most $r$ such that $Q = P_1 + B_1.U$.
\end{prop}

\parag{Proof} First assume that there exists $P \in \mathcal{I}_\sigma$ and $B \in D_{N, \sigma}$ such that the symbol of $Q := P + B.U$ is not in $I_Z + (u)$. Then consider such a couple $(P_0, B_0)$ with $B_0$ of order $b$ minimal among all such couples. Then $P_0$ and $B_0.U$ have the same order because when their orders are different we have $s(Q) = s(P_0)$ or $s(Q) = s(B_0).u$ contradicting the fact that $s(Q)$ is not in $I_Z + (u)$.\\
Also, if $P_0$ and $B_0.U$ have equal orders which is the order of $Q$, we have the equality  $s(Q) = s(P_0) + s(B_0).u$ contradicting our assumption.\\
So the only case left is when $P_0$ and $B_0.U$ have the same order $b_0 + q$ which is strictly bigger than the order of $Q$. In this case we have $ s(P_0) + s(B_0).u = 0$ which implies that $s(B_0).u$ vanishes on $Z$. But our hypothesis on $u$ implies then that $s(B_0)$ vanishes on $Z$. As $I_Z$ is reduced and is the characteristic ideal of $\mathcal{I}$ we may find a germ $B \in \mathcal{I}_\sigma$ such that $s(B) = s(B_0)$. Then write
$$ Q = P_0 + B_0.U = P_0 + B.U + (B_0 - B).U .$$
Since $B$ is in $\mathcal{I}_\sigma$ and $\mathcal{I}.U \subset \mathcal{I}$ we have $P_1 = P_0 + B.U$ in $\mathcal{I}_\sigma$ and the order of $B_0 - B_1$ is strictly less that $b$. This contradicts the minimality of $b$ and proves our first assertion.\\
Assume now that $Q = P + B.U$ has order $q+r$ and that $B$ has order $r + s$ with $s \geq 1$.  If the order of $P$  and $B.U$ are  not equal then either $P$ or $B.U$ is of order $q+r$ and  $P$ and $B.U$ have orders at most $q+r$ we are done. \\
So we may assume that $P$ and $B.U$ have the same order $q+r+s$ with $s \geq 1$. Then the previous considerations will produce $B' \in \mathcal{I}_\sigma$ with $s(B') = s(B)$ and then $P_1 := P + B'.U$ and $B_1 := B - B'$ give that $Q = P_1 + B_1.U$ with $P_1 \in \mathcal{I}_\sigma$ and $B_1 \in D_{N, \sigma}$ of order at most $q + r + s - 1$. By a descending induction on $s$  this completes our proof because when $B_1$ has order at most $r$ the order of $P_1$ is at most $q+r$ because we assume that $Q$ has order $q+r$.$\hfill \blacksquare$\\

The following  two corollaries  are immediate applications of the previous proposition, using the proposition \ref{charact. 1} and the theorem \ref{action} which allow to verify that needed hypotheses.

\begin{cor}\label{charact. 0}
The characteristic cycle of $\mathcal{N}_\lambda, \forall \lambda \in \C$, is the cycle associated to the ideal $I_Z + (g)$ in $\mathcal{O}_N[\eta]$ where $g$ is the symbol of $U_0$.\\
 Also the characteristic cycle of  $\tilde{\mathcal{M}}$ is the cycle associated to the ideal $I_Z + (\gamma)$ in $\mathcal{O}_N[\eta]$ where $\gamma$ is the symbol of $U_{-1}$.$\hfill \blacksquare$
\end{cor}

\begin{cor}\label{our sheaf case}
Let $\mathcal{I}$ the left ideal  in $D_N$ that we introduced in the definition \ref{ideal} and let $U := U_0 -\lambda$. Then for any non zero germ $Q \in \mathcal{I} + D_N.U$ of order $q+1$ there exist a germ $P \in \mathcal{I}$ of order at most $q+1$ and a germ $B \in D_N$ of order at most $q$ such that $Q = P + B.U$. $\hfill \blacksquare$
\end{cor}

\parag{Remark} 
 \begin{enumerate}
\item We shall be interested mainly by the special case of the corollary \ref{our sheaf case}. \\
Define for each $q \geq 0$
$$ \mathcal{J}_\lambda(q+1) = \mathcal{I}(q+1) + D_N(q).(U_0 - \lambda) .$$
Then this corollary gives, for each $\lambda \in \C$ and for each $q \in \mathbb{N}^*$ the equality $\mathcal{J}_{\lambda} \cap D_N(q) = \mathcal{J}_\lambda(q)$. This implies that the natural map
\begin{equation}
\mathcal{N}_\lambda(q) \to \mathcal{N}_\lambda 
\end{equation}
is injective 
\item Note that $\mathcal{J}_\lambda(0) := \mathcal{I}(0) = \{0\}$ as no non zero differential operator of order $0$ annihilates the Newton polynomials (in fact $N_0 := k$ is enough !)
\item Also the fact that $\mathcal{I}(1) = \{0 \}$ (see the  lemma \ref{remplace})  implies the equality  $\mathcal{J}_\lambda(1) = \mathcal{O}_N.(U_0 - \lambda)$.
\end{enumerate}

\parag{The irreducible component $X$}
Let $\Delta := \{\Delta(\sigma) = 0 \}$ in $N$. At the generic point $\sigma$ of this hyper-surface, the polynomial $P_\sigma$ has exactly one double root $\varphi(\sigma)$ and $\varphi : \Delta \dashrightarrow \C$ is a meromorphic function  which is locally bounded on $\Delta$. Then define the meromorphic map 
$$ \Phi : \Delta \dashrightarrow \mathbb{P}_{k-1} $$
by letting $\Phi_h(\sigma) = (- \varphi(\sigma))^{k-h}$ for $h \in [1, k]$ in homogeneous coordinates. Let $X \subset N \times \C^k$ be the $N-$relative cone over the graph of the meromorphic map $\Phi$. This is a $k-$dimensional irreducible subset in $\Delta \times \C^k$ and its fiber at the generic point in $\Delta$ is the line directed by the vector $\Phi_h(\sigma), h \in [1, k]$. \\

We shall consider the following sub-spaces in $N \times \mathbb{P}_{k-1}$ (where $s(P)$ is the symbol of $P$)

  $$\mathbb{P}(Z) := \{(\sigma, \eta) \in N \times \mathbb{P}_{k-1} \ / \ s(P)(\sigma, \eta) = 0 \quad \forall P \in \mathcal{I}\setminus \{0\} \} $$ 
  $$ \mathbb{P}(X) := \{(\sigma, \eta) \in \mathbb{P}(Z) \  / \   \gamma(\sigma, \eta) :=  \sum_{h=0}^{k-1} (k-h).\sigma_h.\eta_{h+1} = 0 \} $$
   $$  \mathbb{P}(Y) := \{(\sigma, \eta) \in \mathbb{P}(Z) \  / \  g(\sigma, \eta) :=  \sum_{h=1}^k h.\sigma_h.\eta_h = 0 \} .$$

The next proposition will justify our notations in proving that   $\mathbb{P}(X)$ is the graph of the meromorphic map $\Phi$ !

\begin{prop}\label{description}
The subspace $\mathbb{P}(Z)$ is a complex sub-manifold of dimension $k$ which is a $k-$branched covering of $N$ via the natural projection $N \times \mathbb{P}_{k-1} \to N$. The sub-space $\mathbb{P}(X)$ is reduced and equal to  the irreducible component of 
$$\mathbb{P}(Z) \cap \big(\{ \Delta = 0\}\times \mathbb{P}_{k-1}\big) $$ 
which is the graph of the meromorphic map $\Phi : \{ \Delta = 0 \} \to \mathbb{P}_{k-1}$ defined above, and $\mathbb{P}(Y)$ is the sum (as a cycle) of \, $\mathbb{P}(X)$ with the reduced hyper-surface in $\mathbb{P}(Z)$ defined by the (reduced) divisor  $\{\eta_{k-1} = 0\}$ in $\mathbb{P}(Z)$.
\end{prop}

\parag{proof} First remark that if $(\sigma, \eta)$ is in $Z$ and satisfies $\eta_k = 0$ then we have $\eta = 0$. So \, $\mathbb{P}(Z)$ is contained in the open set $ \Omega_k := \{\eta_k \not= 0\}$ and, on this open set which is  isomorphic to $N \times \C^{k-1}$, we may use the coordinates $\sigma_1, \dots, \sigma_k, \eta_1/\eta_k, \dots, \eta_{k-2}/\eta_k$ and $z := - \eta_{k-1}/\eta_k$.

\begin{lemma}\label{chart k}
We have an isomorphism
$$\varphi_k :  \mathbb{P}(Z) \to \C^k \quad {\rm given \ by} \quad (\sigma, \eta) \mapsto (\sigma_1, \dots, \sigma_{k-1}, z) .$$
\end{lemma}

\parag{proof} Remark first that the vanishing of the $(2, 2)$ minors which give the equations of $Z$ (see formula $(10)$ before proposition \ref{charact. 1}) implies, assuming $\eta_k \not= 0$, that:
$$ \eta_h/\eta_k = (-z)^{k-h} \quad {\rm for}  \  h \in [1, k-1] $$
 and also, as the symbol of $\mathcal{T}^k$ is equal to $ \eta_1.\eta_{k-1} + \eta_k.l_\sigma(\eta)$, that $l_\sigma(\eta)/\eta_k  = -(-z)^k$. But  then,
  $$l_\sigma(\eta)/\eta_k = \sum_{h=1}^k (-1)^{k-h}.\sigma_h.z^{k-h} = (-1)^k.(P_\sigma(z) - z^k)$$
   shows that $P_\sigma(z) = 0$ on $\mathbb{P}(Z)$.\\
Let us show that the holomorphic map $ \psi_k : \C^k \to  \mathbb{P}(Z)$   given by
$$ \eta_h = (-z)^{k-h} \quad {\rm for}  \  h \in [1, k] \quad {\rm and} \quad \sigma_k = -\sum_{h=0}^{k-1} (-1)^{k-h}.\sigma_h.z^{k-h} $$
with the convention $\sigma_0 \equiv 1$ gives an inverse to $\varphi_k$. \\
First, we shall verify that $\psi_k$ takes its values in $\mathbb{P}(Z)$. Note that the definition of $\sigma_k$  implies $P_\sigma(z) = 0$. We have for $(\sigma, \eta) = \psi_k(\sigma', z)$ the equality:
$$ l_\sigma(\eta) = \sum_{h=1}^k \sigma_h.\eta_h = \sum_{h=1}^k (-1)^{k-h} \sigma_h.z^{k-h} = -(-z)^k$$
Then we have to verify that the vectors $(\eta_1, \dots, \eta_{k-1}, 1)$ and $((-z)^k, \eta_1, \dots, \eta_{k-1})$ are co-linear. This is clear as the second one is $\eta_{k-1} = (-z).\eta_k-$times the first one (see again formula $(10)$).\\
To complete the proof, it is enough to check that $\psi_k \circ \varphi_k $ and $\varphi_k \circ \psi_k$ are the identity maps. This is  easy verification is left to the reader.$\hfill \blacksquare$\\

\parag{End of proof of \ref{description}} In this chart we have 
\begin{align*} 
& g(\sigma, \eta)/\eta_k = \sum_{h=1}^k h.\sigma_h.\eta_h/\eta_k = \sum_{h=1}^k (-1)^{k-h}.h.\sigma_h.z^{k-h} \\
& g(\sigma, \eta)/\eta_k = (-1)^{k+1}\Big(\sum_{h=1}^k (-1)^h.(k-h).\sigma_h.z^{k-h} - k.\sum_{h=1}^k (-1)^h.\sigma_h.z^{k-h}\Big)
\end{align*}
 and this gives
$$ g(\sigma, \eta)/\eta_k = (-1)^{k+1}\Big(z.P'_\sigma(z) -k.z^k - k.(P_\sigma(z) - z^k)\Big) = (-1)^{k+1}.z.P'_\sigma(z) .$$
We have also 
$$\gamma(\sigma, \eta)/\eta_k = \sum_{h=0}^{k-1} (-1)^{k-h-1}.(k-h).\sigma_h.z^{k-h-1} = (-1)^{k+1}.P'_\sigma(z) .$$
So $g = z.\gamma$ in this chart\footnote{compare withe the formula $(F)$ at the end of the paragraph 3.1.}, and the ideal generated by $g$ and $\gamma$ in $\mathbb{P}(Z)$  is generated by $\gamma$ which defined the  hyper-surface $\mathbb{P}(X)$.\\
But on this hyper-surface we have $P_\sigma(z) = 0$ and $P'_\sigma(z) = 0$, so $z$  is a double root of $P_\sigma$. This implies that $\Delta(\sigma) = 0$ for $(\sigma, \eta)$ in the analytic subset $ \vert \mathbb{P}(X)\vert$. \\
On a Zariski dense open set in $\Delta$ the unique double root of $P_\sigma$ is equal to  $\varphi(\sigma) $ which is given by $z = -\eta_{k-1}/\eta_k$ when $(\sigma, \eta) \in \vert \mathbb{P}(X)\vert$. So $\vert \mathbb{P}(X)\vert$ contains the graph of the meromorphic map $\Phi$. Moreover, as  the projection $\mathbb{P}(Z) \to N$ is clearly a branched covering (of degree $k$) and over the generic point in $\{\Delta(\sigma) = 0$ there exists an unique root of multiplicity  $2$ for $P_{\sigma}$, $\mathbb{P}(X)$ has generic degree $1$ over $\{\Delta = 0 \}$. And because $P''_\sigma(z)$ does not vanish at the generic point in $\mathbb{P}(X)$ (which has to be over the generic point of $\{\Delta = 0 \}$) implies that the hyper-surface $\mathbb{P}(X)$ of $\mathbb{P}(Z)$ is reduced. This is enough to conclude that $\mathbb{P}(X)$ is equal to the graph of $\Phi$.\\
The previous computation shows also that $\mathbb{P}(Y)$ is the sum of $\mathbb{P}(X)$ with the divisor define by $\{z = 0\}$ in $\mathbb{P}(Z)$ which is a smooth and reduced hyper-surface
given by the equation  $\eta_{k-1} = 0$ in $\mathbb{P}(Z)$.$\hfill \blacksquare$\\

The determination of the characteristic cycles of the holonomic $D_N-$modules $\tilde{\mathcal{M}}$ and $\mathcal{N}_\lambda$ is an immediate consequence of the previous proposition thanks  to  \ref{charact. 0}.

\begin{cor}\label{caract. cycles}
The characteristic cycle of the $D_N-$module $\tilde{\mathcal{M}}$ is equal to $\mathbb{P}(X)$.
For each complex number $\lambda$ the characteristic cycle of the $D_N-$module $\mathcal{N}_\lambda$ is equal to $ \mathbb{P}(Y) = \mathbb{P}(X) + \big(\mathbb{P}(Z) \cap \{ \eta_{k-1} = 0 \}\big) $.$\hfill \blacksquare$
\end{cor}

\parag{Remarks}
\begin{enumerate}
\item The intersection $ \mathbb{P}(Z) \cap \{ \eta_{k-1} = 0 \} $ is equal to $N \times [v]$ where $v$ is the point $ (0, \dots, 0, 1) \in \mathbb{P}_{k-1}$ and this intersection is the projectivization of the co-normal to the hyper-surface $\{\sigma_k = 0\}$.
\item At the set-theoretical level we have
$$ Z \cap \{\gamma = 0 \} = X \cup (N\times \{0\}) \quad {\rm and}$$
$$Z \cap \{g = 0\} = X \cup\Big( \{\sigma_k = 0\} \times\{\eta_1= \eta_2 = \dots = \eta_{k-1} = 0 \} \Big) \cup (N\times \{0\}) .$$
\item Despite the previous results, $g(\sigma, \eta)$ does not belongs to the ideal of $\C[\sigma, \eta]$ generated by $I_Z$ and $\gamma(\sigma, \eta)$ at the generic point in $N\times \{0\}$. This is consequence of the fact that $I_Z$ does not contain a non zero element in $\C[\sigma, \eta]$ which is homogeneous of degree $1$ in $\eta$, using the corollary \ref{our sheaf case}.
\end{enumerate}

The following  lemma will be useful later on

\begin{lemma}\label{useful 2}
Assume that $f.\partial_k^n.U_{-1}$ is in $\mathcal{J}_{\lambda, \sigma}$ for some $f \in \mathcal{O}_{N, \sigma}$, some integer $n \geq 1$ and some $\lambda \in \C$. Then $f$ is in $\sigma_k.\mathcal{O}_{N, \sigma}$.
\end{lemma}

\parag{proof} The fact that $f.\partial_k^n.U_{-1}$  is in $\mathcal{J}_{\lambda, \sigma}$  implies that $f.\eta_k^n.\gamma$ vanishes on the characteristic variety of the $D_N-$module $\mathcal{N}_\lambda$. So thanks to the corollary \ref{caract. cycles}  $f.\eta_k^n.\gamma$ vanishes on $C$, the co-normal bundle  of the hyper-surface $\{\sigma_k = 0 \}$.
But  $\eta_k$ and  $\gamma$ do not vanish on any non empty open set in  $C$:\\
 this is clear for $\eta_k$ and the restriction of $\gamma$ to $C$ is equal to $\sigma_{k-1}.\eta_k$ and $\sigma_{k-1}$ also does not vanish   on any non empty open set in $C$. So $f \in \mathcal{O}_{N, \sigma}$ has to vanish on $C$ and we conclude that $f$ is in $\sigma_k.\mathcal{O}_{N, \sigma}$. $\hfill \blacksquare$\\

\subsection{The case $\lambda \not\in \mathbb{N}$}

\parag{Notation} For each $\lambda \in \C$ and each $q \geq 0$ we shall note
 $$\mathcal{J}_\lambda(q+1) := \mathcal{I}(q+1) + D_N(q).(U_0 - \lambda)$$ 
  and
 $$\mathcal{N}_\lambda(q+1) := D_N(q+1)\big/\mathcal{J}_\lambda(q+1).$$
 For $q = 0$ we note $\mathcal{J}_\lambda(0) := \mathcal{I}(0)$ and $\mathcal{N}_\lambda(0) := \mathcal{O}_N\big/\mathcal{J}_\lambda(0)$.\\

The goal of this paragraph is to prove the following theorem.

\begin{thm}\label{no torsion}
For $\lambda \in \C \setminus \mathbb{N}^*$ the $D_N-$module $\mathcal{N}_\lambda$ has no $\mathcal{O}_N-$torsion. 
\end{thm}

\parag{Proof} This result is a direct consequence of  the  proposition \ref{non torsion}, thanks to the injectivity for each $q \geq 0$  of the natural map $\mathcal{N}_\lambda(q) \to \mathcal{N}_\lambda$ (see remark 1 following the corollary \ref{our sheaf case}).$\hfill \blacksquare$\\

The case $\lambda = 0$ will be seen separately at the end of this paragraph.

\begin{defn}\label{theta}
For any  $\lambda \in \C \setminus \mathbb{N}^*$, for any integer  $q \geq 2$ and for any integer $r \in [q, k.(q-1)]$ define the following elements in $\mathcal{W}_q$ (see formulas $(12)$ and $(13)$ in lemma \ref{computation})
\begin{align}
& \theta_{q, r} := (r - \lambda).[\partial^\alpha.\mathcal{T}^m]  - (q - 1).[\partial^\beta.(U_0 - \lambda)] \\
&{\rm so} \quad   \theta_{q, r} = \sum_{h=0}^k  (r - \lambda - (q-1).h).\sigma_h.y_{q, r+h}.
\end{align}
where in formula $(23)$ we assume that  $\alpha \in \mathbb{N}^k$ and $m \in [2, k]$  satisfy $\vert \alpha \vert = q-2$ and $w(\alpha) = r-m$, and that $\beta \in \mathbb{N}^k$ satisfies $\vert \beta \vert = q-1$ and $w(\beta) = r$.
\end{defn}

\begin{cor}\label{charact. 2}
 For any integer $q \geq 1$ the kernel of the quotient map
  $$ l_q : \mathcal{W}(q) \to \mathcal{N}_\lambda(q)$$
  is equal to the sub$-\mathcal{O}_N-$module generated by $U_0 - \lambda \in \mathcal{W}(1)$ and the elements $\theta_{p, r}, \forall p \in [2, q],  \forall r \in [p, (k-1).p],$. 
\end{cor}

\parag{Proof} We have to prove that if a non zero differential operator $P$ of order $p \leq q$ is in $\mathcal{J}_\lambda$ then it may be written as $Q + B.(U_0 -\lambda)$ with $Q \in \mathcal{I}$ of order at most $p$ (or $Q = 0$) and $B$ of order at most $p-1$ (or $B = 0$). When $p \geq 2$ this is precisely the statement proved in the corollary \ref{our sheaf case}. For $p \leq 1$ the only $P$ which are in $\mathcal{J}_\lambda(1)$ are in $\mathcal{O}_N.(U_0 - \lambda)$ thanks to the  remarks 2 and 3 following the corollary \ref{our sheaf case}. $\hfill \blacksquare$\\ 

\begin{lemma}\label{basis}
Let  $\lambda$ be in  $ \C \setminus \mathbb{N}^*$;  for each integer $q \geq 2$ the elements $\theta_{q, r}$ and $y_{q, s}$, with $r \in [q, k.(q-1)]$ and $ s \in [k.(q-1)+1, k.q]$ form a $\mathcal{O}_N-$basis of $\mathcal{W}_q$.
\end{lemma}

\parag{proof} Let $\mathcal{W}_{q, p}$ be the $\mathcal{O}_N-$module of $\mathcal{W}_q$ with basis the $y_{q, r}$ for $r \geq p+1$. Then we have for $r \in [q, k.(q-1)]$
$$ \theta_{q, r} \in  (r- \lambda).y_{q, r} + \mathcal{W}_{q, r+1} $$
so the determinant of the $k.(q-1) - q + 1 + k  = k.q - q + 1$ vectors $\theta_{q, r}, y_{q, s}$ in the basis $\big(y_{q, r}, r \in [q, k.q]\big)$ of $ \mathcal{W}_q $ is  upper triangular and is equal to $\prod_{r = q}^{k.(q-1)} (r - \lambda) $ which is in $\C^*$ as soon as $\lambda$ is not in the subset $[q, k.(q-1)]$ of $\mathbb{N}^*$.$\hfill \blacksquare$\\

\begin{prop}\label{non torsion}
Let $q \geq 1$ be an integer and assume that $\lambda $ is not an integer in $[0, k.(q-1)]$. Let $L_q : \mathcal{W}_q \to \mathcal{N}_\lambda(q)$ be the restriction to $\mathcal{W}_q$ of  quotient map $l_q$. This $\mathcal{O}_N-$linear map is surjective and its kernel  is the sub-module of $\mathcal{W}_q$ with basis the $\theta^q_r$ for  $r \in [q, k.(q-1)]$. So $\mathcal{N}_\lambda(q)$ is a free $\mathcal{O}_N-$module of  rank $k $.
\end{prop}

\parag{Proof} Remark first that for $q= 1$ the result is clear as for $\lambda \not= 0$ we have $$\mathcal{N}_\lambda(1) = \oplus_{h=1}^k \mathcal{O}_N.\partial_h $$ 
thanks to the remark 3 following the corollary \ref{our sheaf case}  and $\mathcal{W}_1 = \oplus_{h=1}^k \mathcal{O}_N.y_{1, h}$ with $L_1(y_{1, h}) = [\partial_h]$. So we may assume that $q \geq 2$.\\
We shall prove first that $\mathcal{N}_\lambda(q)$ is equal to the image of $L_q$ by induction on $q \geq 2$.\\
Assume that $q = 2$. Then the image of 
$$\partial_j.(U_0 - \lambda) - (j-\lambda).\partial_j = \sum_{h=1}^k h.\sigma_h.y_{2, h+j} \in \mathcal{W}_2$$
 by $L_2$ is the class of $ -(j-\lambda).\partial_j$ in $\mathcal{N}_\lambda$. So the image of $L_2$ contains the classes of  $\partial_1, \dots, \partial_k$ as $\lambda$ is not in $[1, k]$ and also contains the class of $1$ as we assume $\lambda \not= 0$ and as  the equality $\lambda = \sum_{h=1}^k h.\sigma_h.\partial_h $ holds in $\mathcal{N}_\lambda$. But the image of $L_2$ contains obviously the classes of $\partial^\alpha$ for any multi-index $\alpha \in \mathbb{N}^k, \vert \alpha \vert = 2$. So our assertion is proved for $q = 2$.\\
Assume now that $q \geq 3$ and that our assertion is proved for $q-1$. Remark that the image of $L_q$ contains obviously the classes of $\partial^\alpha$ for each $\alpha \in \mathbb{N}^k, \vert \alpha \vert = q$. We shall use now the following easy formula:
\begin{itemize}
\item For any $r \in [q-1, k.(q-1)]$ and any $j \in [1, k]$ we have in $\mathcal{N}_\lambda$ the equality
$$ \partial_j L_{q-1}(y_{q-1, r}) = L_q(\partial_j y_{q-1,r}) =  L_q(y_{q, r+j}) .$$
\end{itemize}
For any $\beta \in \mathbb{N}^k \setminus \{0\}$ with $\vert \beta\vert \leq q-1$ we may find $j \in [1, k]$ and $\gamma \in \mathbb{N}^k$ such that $\partial^\beta = \partial_j\partial^\gamma$. By our inductive assumption there exists $x \in \mathcal{W}_{q-1}$ such that $L_{q-1}(x) = \partial^\gamma$. Then $\partial_jx$ is in $\mathcal{W}_q$   and thanks to the formula above we have
 $$L_q(\partial_jx) = \partial_jL_{q-1}(x) = \partial_j.\partial^\gamma = \partial^\beta \quad {\rm in } \quad \mathcal{N}_\lambda .$$
 Again we conclude that the class of $1$ in $\mathcal{N}_\lambda(q)$ is in the image of $L_q$ using $\lambda \not= 0$ and the equality $\lambda = \sum_{h=1}^k h.\sigma_h.\partial_h $ which holds in $\mathcal{N}_\lambda(q)$.
 This complete the proof of our first statement.\\
 But it is clear that $\theta_{q,r}$ for $r \in [q, k.(q-1)]$ are in the kernel of $L_q$. So the  $\mathcal{O}_N-$free rank $k$ module with basis $\big(y_{q, r}, r \in [k.(q-1)+1, k.q]\big)$ is surjective via $L_q$ onto $\mathcal{N}_\lambda(q)$.  The next lemma completes the proof, as we already know that $\mathcal{N}_\lambda(1)$ is a $\mathcal{O}_N-$free rank $k$ sub-module of $\mathcal{N}_\lambda(q)$ with basis $\partial_1, \dots, \partial_k$. $\hfill \blacksquare$\\
 
 \begin{lemma}\label{algebre}
 Let $A$ be an integral commutative ring and let $M$ be a $A-$module. Assume that there exists a surjective $A-$linear map $p : A^k \to M$ and an injective $A-$linear map $i : A^k \to M$. Then $p$ is an isomorphism.
 \end{lemma}
 
 \parag{proof} Let $j : A^k \to A^k$ be a $A-$linear map such that $j\circ p = i$. So $j$ is injective and the  co-kernel $C$ of $j$ is a torsion module. Let $q : A^k \to C$ be the quotient map and let $K$ be the kernel of $p$. The restriction of $q$ to $K$ is injective because if $x \in K$ satisfies $q(x) = 0$ then $x = j(y)$ for some $y \in A^k$ and then $i(y) = p(j(y)) = p(x) = 0$, which implies $y = 0$ and $x=0$. So $K$ is a sub-module of $C$ and then $K$ is a $A-$torsion module. But as $K \subset A^k$ we have $K = 0$ and so $p$ is an isomorphism.$\hfill \blacksquare$

$$\xymatrix{ 0 & & & & & \\  & C \ar[lu] & & & & \\ 0 \ar[r] & K \ar[u] \ar[r] & A^k \ar[r]^p  \ar[ul]_q& M \ar[r] & 0 \\  & 0 \ar[u] & & A^k \ar[u]_{i} \ar[ul]_j &  \\  & & & 0 \ar[u] & 0 \ar[ul] } $$

\begin{lemma}\label{connection 1}
For $\lambda \not\in \mathbb{N}$ we have $\sigma_k.\Delta(\sigma).\mathcal{N}_\lambda(2) \subset \mathcal{N}_\lambda(1)$.
\end{lemma}

The proof will be a simple consequence of the following lemma.

\begin{lemma}\label{determinant}
Let $y := (y_{2}, \dots, y_{2k})$ be  in $\C[\sigma]^{2k-1}$ and consider the $\C[\sigma]-$linear system $(2k-1, 2k-1)$ on $\C[\sigma]^{2k-1}$ given by the following $\C[\sigma]-$linear forms:
\begin{align*}
& L_{q}(y) := \sum_{h=0}^{k} h.\sigma_{h}.y_{q+h} \quad {\rm for } \quad q \in [1, k] \\
& \Lambda_{r}(y) := \sum_{h=0}^{k} \sigma_{h}.y_{r+h} \quad {\rm for} \quad r \in [2, k]
\end{align*}
Then the determinant of this linear system is equal to $\sigma_{k}.\Delta(\sigma)$ where $\Delta(\sigma)$ is the discriminant of the polynomial $P_{\sigma}(z) := \sum_{h=0}^{k} (-1)^{h}.\sigma_{h}.z^{k-h}$ with the convention $\sigma_{0}\equiv 1$.
\end{lemma}

\parag{Proof} Remark first that $\Delta(\sigma)$ is also the discriminant of the polynomial (see the computation below):
 $$\tilde{P}_{\sigma}(z) := \sum_{h=0}^{k} \sigma_{h}.z^{k-h}.$$
Then remark also that the resultant of the polynomials $\tilde{P}_{\sigma}(z)$ and $k.\tilde{P}_{\sigma}(z) - z.(\tilde{P}_{\sigma})'(z)$ coincides with the determinant of the $(2k-1, 2k-1)$ $\C[\sigma]-$linear system  defined in the statement of the lemma. So it is enough to compute this resultant. It is given by
$$ R(\sigma) = \prod_{\tilde{P}_{\sigma}(z_{j}) = 0} \big(k.\tilde{P}_{\sigma}(z_{j}) - z_{j}.(\tilde{P}_{\sigma})'(z_{j})\big) = \sigma_{k}.\prod_{P_{\sigma}(-z_{j}) = 0} (-1)^{k-1}.P'_{\sigma}(-z_{j}) = \sigma_{k}.\Delta(\sigma)$$ 
as $\tilde{P}_\sigma(-z) = (-1)^{k}.P_\sigma(z) $ implies $\tilde{P}'_\sigma(-z) = (-1)^{k-1}.P'_\sigma(z)$. This conclude the proof.$\hfill \blacksquare$\\

\parag{proof of \ref{connection 1}} It is enough to prove that for each $(p, q) \in [1, k]^2$ there exist polynomials $a_{h, q}^p(\lambda)$ in $\C[\sigma, \lambda]$ (in fact affine in $\lambda$) such that
$$ \sigma_k.\Delta(\sigma).\partial_p\partial_q - \sum_{h=1}^k  a_{h, q}^p(\lambda).\partial_h \in \mathcal{J}_\lambda .$$
For $m \in [2, k] $ we have
$$ \mathcal{T}^m = y_{2, m} + \sum_{h=1}^k \sigma_h.y_{2, m+h} + y_{1, m} \in \mathcal{I} \subset \mathcal{J}_\lambda $$
and for $q \in [1, k]$:
$$ \partial_q.(U_0 - \lambda) = \sum_{h=1}^k h.\sigma_h.y_{2, q+h} + (q-\lambda).y_{1, q} \in \mathcal{J}_\lambda .$$
This gives $(2k-1)$ $\C[\sigma]-$linear relations between the basis elements $y_{2, r}, r \in [2, 2k]$ of $\mathcal{W}_2 $ modulo $L_2^{-1}(\mathcal{N}_\lambda(1))$. But the determinant of these $2k-1$ vectors in the basis $y^2_r$ of $\mathcal{W}_2$ is equal to $\sigma_k.\Delta(\sigma)$ thanks to the previous lemma. The conclusion follows, as we know that $L_2 : \mathcal{W}_2 \to \mathcal{N}_\lambda(2)$ is surjective for $\lambda \not\in \mathbb{N}$.$\hfill \blacksquare$\\

\begin{lemma}\label{petit}
Assume that on a Stein open set $U$ in $N$  the equality $\mathcal{N}_\lambda(2)_{\vert U} = \mathcal{N}_\lambda(1)_{\vert U}$ is true for some $\lambda \in \C$. Then we have
$$ (\mathcal{N}_\lambda)_{\vert U} = \mathcal{N}_\lambda(1)_{\vert U} .$$
\end{lemma}

\parag{proof} It is enough to prove the equality $\mathcal{N}_\lambda(q)_{\vert U} = \mathcal{N}_\lambda(1)_{\vert U}$ for any $q \geq 2$ because we know that $\mathcal{N}_\lambda = \cup_{q \geq 0}\  \mathcal{N}_\lambda(q)$. As this is true for $q = 2$ by assumption, we shall prove this equality by induction on $q \geq 2$. So assume that this equality is proved for some $q \geq 2$ and we shall prove it for $q+1$.\\
Let $\alpha \in \mathbb{N}^k$ such that $\vert \alpha\vert = q+1$ and write $\partial^\alpha = \partial_p\partial^\beta$ for some $p \in [1, k]$ and some $\beta \in \mathbb{N}^k$ with $\vert \beta\vert = q$. By the inductive assumption we may write $\partial^\beta = \sum_{h=1}^k b_h.\partial_h$ in $\mathcal{N}_\lambda(q)$ with $b_h \in \mathcal{O}(U)$ because we know that $\mathcal{N}_\lambda(1) = \oplus_{h=1}^k \mathcal{O}_N.\partial_h$ on $N$. Then we obtain that $\partial_p\partial^\beta $ is in $\mathcal{N}_\lambda(2)_{\vert U} = \mathcal{N}_\lambda(1)_{\vert U}$, concluding our induction.$\hfill \blacksquare$\\

\begin{cor}\label{connection 2}
For each $\lambda \in \C \setminus \mathbb{N}$ there exists a meromorphic integrable connection $\nabla_\lambda : \mathcal{O}_N^k \to \frac{1}{\sigma_k.\Delta}.\mathcal{O}_N^k\otimes \Omega^1_N$ with a simple pole on the reduced hyper-surface $\{\sigma_k.\Delta(\sigma) = 0 \} \subset N$ such that the restriction of  $\mathcal{N}_\lambda$ to the Stein (in fact affine)  open set  \ $U := \{\sigma_k.\Delta(\sigma) \not= 0 \}$ is isomorphic to the $D_U-$module defined by $(\mathcal{O}_N^k, \nabla_\lambda)$. Moreover, this isomorphism is the restriction of an injective $D_N-$linear map $$\mathcal{N}_\lambda \to \big(\mathcal{O}_N^k(*\sigma_k.\Delta(\sigma)), \nabla_\lambda\big).$$
\end{cor}

\parag{proof} This is an easy consequence of the $\mathcal{O}_N$ isomorphism $\mathcal{N}_\lambda(1) \to \oplus_{h=1}^k \mathcal{O}.\partial_h$ and the previous  lemmas \ref{connection 1} and \ref{petit}.$\hfill \blacksquare$\\

The following results will complete the theorem \ref{no torsion} for $\lambda = 0$ (see lemma below) and gives in this case the analog of the corollary \ref{connection 2}.

\parag{Notations} Let $\varphi_0 : \mathcal{N}_0 \to \mathcal{O}_N$ be the $D_N-$linear map defined by  $\varphi_0(1) = 1$. Let $\mathcal{N}_0^*$ its kernel. Then we have a $D_N-$linear direct sum decomposition 
$$ \mathcal{N}_0 = \mathcal{N}_0^* \oplus \mathcal{O}_N $$
given by $\varphi_0$ and the map $P \mapsto P - P[1]$ from $\mathcal{N}_0$ to $\mathcal{N}_0^*$.\\
Define $\mathcal{N}_0^*(1) := \mathcal{N}_0^* \cap \mathcal{N}_0(1)$. This $\mathcal{O}_N-$module is isomorphic to the  coherent $\mathcal{O}_N-$module
 $$\mathcal{F} := \oplus_{h=1}^k \mathcal{O}_N.\partial_h\big/ \mathcal{O}_N.(\sum_{h=1}^k h.\sigma_h.\partial_h)$$
 which is locally free of rank $(k-1)$  on $N\setminus \{0\}$. It has no torsion.

\begin{lemma}\label{complement}
The $D_N-$module $\mathcal{N}_0$ has no torsion.
\end{lemma}

\parag{proof} We shall apply the proposition \ref{no torsion} to the $D_N-$module $\mathcal{N}_0^*$. In this case the map $L_1 : W_1 \to \mathcal{N}_0^*(1) = \mathcal{F}$ is clearly surjective and the induction for the surjectivity of $L_q$ for $q \geq 2$ stays valid.$\hfill \blacksquare$\\

\begin{cor}\label{connection 3}
There exists a meromorphic (regular) integrable connection 
 $$\nabla_0 : \mathcal{F} \to \frac{1}{\sigma_k.\Delta}.\mathcal{F}\otimes \Omega^1_N$$
 with a simple pole on the reduced hyper-surface $\{\sigma_k.\Delta(\sigma) = 0 \} \subset N$ such that the restriction of  $\mathcal{N}_0^*$ to the open set $U = \{\sigma_k.\Delta(\sigma) \not= 0 \}$ is isomorphic to the $D_U-$module associated to $(\mathcal{F}, \nabla_0)$. Moreover, this isomorphism is the restriction of an injective $D_N-$linear map $\mathcal{N}_0^* \to \big(\mathcal{F}(*\sigma_k.\Delta(\sigma)), \nabla_0\big)$.
\end{cor}

\parag{proof} It is enough to remark that the key point is to extend the lemma \ref{connection 1} to $\mathcal{N}_0^*$. For this purpose it is enough to extend the proposition \ref{no torsion} to $\mathcal{N}_0^*$. This is done in the previous lemma \ref{complement}.$\hfill \blacksquare$

\parag{Remark} A horizontal basis of the regular holonomic connection $(\mathcal{F}, \nabla_0)$ is given by the local branches $Log\, z_1 - Log\, z_j, j \in [2, k]$ where $z_h, h \in [1, k]$ are the local branches of the roots of $P_\sigma(z) = 0$. This comes from the fact  that the local branches $Log\, z_h$ are trace functions so are killed by $\mathcal{I}$ and also because  
$$U_0(Log\, z_h) = (\sum_{j=1}^k z_j.\frac{\partial}{\partial z_j}(Log\, z_h)) = 1$$
 so that $U_0$ kills the $Log\, z_1 - Log\, z_j$ for $j \in [2, k]$.

\bigskip

We shall conclude this section by the following theorem.

\begin{thm}\label{minimal}
Let $\lambda \in \C \setminus \mathbb{Z}$. Then $\mathcal{N}_\lambda$ is the minimal extension of the meromorphic connection given by $(\mathcal{N}_\lambda(1), \nabla_\lambda)$. So $\mathcal{N}_\lambda$ is a simple $D_N-$module.
\end{thm}

\parag{Proof} To see that $\mathcal{N}_\lambda$ is the minimal extension of the simple pole meromorphic  connection $\big(\mathcal{N}_\lambda(1), \nabla_\lambda\big)$ it is enough to prove that $\mathcal{N}_\lambda$ has no torsion, and this is given by the proposition \ref{non torsion},  and no "co-torsion", that is to say that there is no non trivial coherent left ideal $\mathcal{K}$ in $D_N$  containing $\mathcal{J}_\lambda$ and generically equal to $\mathcal{J}_\lambda$ on $N$. Such an ideal defines a  holonomic quotient $Q$ of $\mathcal{N}_\lambda$ which is supported in a closed analytic subset $S$ of $N$ with empty interior in $N$. As $\mathcal{N}_\lambda$ is a quotient of $\mathcal{M}$, we may apply the corollary \ref{co-torsion 2} and so it is enough to show that near the generic points of $\{\sigma_k.\Delta(\sigma) = 0\}$ such an ideal $\mathcal{K}$ is equal to $\mathcal{J}_\lambda$ or to $D_N$.\\
Near the generic point of $\{\sigma_k = 0\}$ we have $\Delta \not= 0$ and we may use a local isomorphism of $N$ given by a holomorphic section of the quotient map
 $$q : M = \C^k \to \C^k\big/\mathfrak{S}_k = N.$$
 Via such an isomorphism  $\mathcal{N}_\lambda$ is the quotient of $D_{\C^k}$ by the left ideal generated by the $\frac{\partial^2}{\partial z_i \partial z_j}$ for $i \not= j \in [1, k]$ and $\sum_{j=1}^k z_j\frac{\partial}{\partial z_j} - \lambda$. The lemma below allows to conclude this case.
For the other case, that is to say near  the generic point of $ \{ \Delta = 0 \}$, the theorem  \ref{co-torsion delta}  completes the proof.\\
The fact that $\mathcal{N}_\lambda$ is a simple $D_N-$module is then consequence of the irreductibility of the monodromy representation of its associated meromorphic connection.  $\hfill \blacksquare$\\

\begin{lemma}\label{no co-torsion 1}
Let $\mathcal{J}_\lambda$ for $\lambda \not\in -\mathbb{N}^*$ be the ideal in $D_{\C^k}$ generated by the differential operators  $\frac{\partial^2}{\partial z_i\partial z_j}$ for $1 \leq i < j \leq k$ and $\sum_{h=1}^k z_h.\frac{\partial}{\partial z_h} - \lambda$.
Let assume that $Q$ is a quotient of the $D_{\C^k}-$module $\mathcal{N}_{\lambda} := D_{\C^k}\big/\mathcal{J}_\lambda$  in a neighborhood $U$ of the point $(z^0_1, \dots, z^0_k)$ in $\C^k$ where $z_1^0 = 0$ and $z_i \not= z_j$ for $1 \leq i < j \leq k$, with support in $\{z_1 = 0 \}$. Then $Q = 0$.
\end{lemma}

\parag{proof} Assume that $Q \not= 0$ Then $Q = D_U\big/\mathcal{K}$ where $\mathcal{K}$ is a left ideal in $D_U$ such that $\mathcal{J}_\lambda \subsetneq \mathcal{K} \subsetneq D$. Then restricting the open neighborhood $U$ of $z^0$ if necessary, there exists a positive integer $n$ such that $z_1^n $ belongs to $\mathcal{K}$\footnote{The class of $1$ in $Q$ is of $z_1-$torsion !}. Then we have
\begin{align*}
& \frac{\partial}{\partial z_1}.z_1^n = n.z_1^{n-1} + z_1^n.\frac{\partial}{\partial z_1} \in \mathcal{K} \quad {\rm so \ writing\ this \ as} \\
& n.z_1^{n-1} + z_1^{n-1}.(\sum_{h=1}^k z_h.\frac{\partial}{\partial z_h} - \lambda) + \lambda.z_1^{n-1} - z_1^{n-1}.(\sum_{h=2}^k z_h.\frac{\partial}{\partial z_h}) \in \mathcal{K} \quad {\rm and \ then}\\
& (n + \lambda).z_1^{n-1} - \sum_{h=2}^k z_h.z_1^{n-1}.\frac{\partial}{\partial z_h} \in \mathcal{K} \tag{@}
\end{align*}
as \  $\sum_{h=1}^k z_h.\frac{\partial}{\partial z_h} - \lambda \in \mathcal{J}_\lambda \subset \mathcal{K}$ on $U$. But $z_1^n \in \mathcal{K}$ implies also, for each $j \in [2, k]$:
\begin{align*} 
& \frac{\partial^2}{\partial z_1\partial z_j}.z_1^n = n.z_1^{n-1}.\frac{\partial}{\partial z_j} + z_1^n. \frac{\partial^2}{\partial z_1\partial z_j} \in \mathcal{K} \quad {\rm which \ implies}\\
& n.z_1^{n-1}.\frac{\partial}{\partial z_j} \in \mathcal{K} \quad \forall j \in [2, k] \tag{@@}
\end{align*}
again as \  $\mathcal{J}_\lambda \subset \mathcal{K}$. Combining $(@)$ and $(@@)$  we conclude that $z_1^{n-1}$ belongs to $\mathcal{K}$, as we assume $n > 0$ and $\lambda \not\in \mathbb{N}^*$.\\
By a descending induction on $n$ we conclude that $1$ belongs to $\mathcal{K}$ which contradicts our assumption that $Q$ is not $0$.$\hfill \blacksquare$\\

\parag{Remark} Note that the $D_N-$linear map 
 $$\varphi_{-1} : \mathcal{N}_{-1} \to \underline{H}_{[\sigma_k = 0]}^1(\mathcal{O}_N) := \mathcal{O}_N[\sigma_k^{-1}]\big/\mathcal{O}_N$$
 defined by $\varphi_{-1}(1) := \sigma_{k-1}/\sigma_k$  is surjective because $\varphi_{-1}(\partial_{k-1}) = 1/\sigma_k$. This shows that for $p= -1$ the sheaf $\mathcal{N}_{-1}$ has a non zero quotient supported by $\{\sigma_k = 0\}$. Then, using the isomorphism $\mathcal{T}_\lambda : \mathcal{N}_{\lambda-1} \to \mathcal{N}_\lambda$ for $\lambda \in -\mathbb{N}^*$ in order to deduce the case $\lambda-1$ from the case $\lambda$ for each $\lambda \in - \mathbb{N}^*$, we see that the sheaf $\mathcal{N}_{-p}$ has a non zero quotient supported by $\{\sigma_k = 0\}$ for any $p \in \mathbb{N}^*$.

\section{The $D_N-$modules $\mathcal{N}_p, p \in \mathbb{Z}$}

\subsection{Structure of  $\mathcal{N}_p, p \geq 1$}

The first important remark is that, thanks to the lemma \ref{isom 2}, it is enough to determine the structure of $\mathcal{N}_1$ as for each $p \geq 2$ the $D_N-$module $\mathcal{N}_p$ is isomorphic to $\mathcal{N}_1$ via the right multiplication by $U_1^{p-1}$.

\subsubsection{Minimality of $\mathcal{N}_1^\square$}

 Recall that $\mathcal{N}_1^\square$ is the co-kernel of the left $D_N-$linear map $\mathcal{T}_1 : \mathcal{N}_0 \to \mathcal{N}_1$ defined by the right multiplication by $U_{-1}$. \\
Thanks to formulas $E_h, h \in [2,k]$ (see the  proposition \ref{magic}) we obtain that $\mathcal{N}_1^\square$ is the quotient of $D_N$ by the left ideal $\mathcal{A} + D_N.(U_0 - 1) + D_N.U_{-1}$ because these formulas imply that the partial differential operators $\mathcal{T}^m, m \in [2, k]$ are contained in $ \mathcal{A} + D_N.(U_0 - 1) + D_N.U_{-1}$ and we have $\mathcal{J}_1 = \mathcal{I} + D_N.(U_0 - 1)$ by definition (see the formula $(1)$ for the definition of $\mathcal{I}$ and the  beginning of the paragraph 3.1 for the definition of $\mathcal{J}_\lambda$).\\
We shall note $\mathcal{N}_1^\square(q) := D_N(q)\big/\big(\mathcal{J}_1 \cap D_N(q)\big)$ for each integer $q \geq 0$.

\begin{prop}\label{injective}
For each $q$ the natural map $\mathcal{N}_1^\square(q) \to \mathcal{N}_1^\square$ is injective.
\end{prop}

\parag{Proof} The proof will use the proposition \ref{Nouvelle version} two times : the first time for the left ideal $\mathcal{A}$ and with $U := U_{-1}$ and the second time for the left  ideal $\mathcal{A} + D_N.U_{-1}$ and with $U := U_0 - 1$. This will give the equalities
\begin{align*}
&  \big(\mathcal{A} + D_N.U_{-1}\big) \cap D_N(q) = \mathcal{A}(q) + D_N(q-1).U_{-1} \quad {\rm and} \\
&  \big(\mathcal{A} + D_N.U_{-1} + D_N.(U_0 - 1) \big) \cap D_N(q) = \mathcal{A}(q) + D_N(q-1).U_{-1} + D_N(q-1).(U_0 - 1).
\end{align*}
This will conclude the proof.\\
In order to apply the proposition  \ref{Nouvelle version} we have to show that the following properties
\begin{enumerate}[i)]
\item The coherence of $\mathcal{A}$ and of $\mathcal{A} + D_N.U_{-1}$.
 \item The fact that the characteristic ideals of $\mathcal{A}$ and of $\mathcal{A} + D_N.U_{-1}$ are reduced.
 \item The inclusions $\mathcal{A}.U_{-1} \subset \mathcal{A}$ and $\big(\mathcal{A} + D_N.U_{-1}).(U_0 - 1)\big) \subset \mathcal{A} + D_N.U_{-1}$.
 \item The symbol of $U_{-1}$ does not vanish on any non empty  open set of the characteristic variety of $D_N\big/\mathcal{A}$.
 \item The symbol of $U_0 - 1$ does not vanish on any non empty  open set of the characteristic variety of $D_N\big/\mathcal{A} + D_N.U_{-1}$.
 \end{enumerate}
 The point $i)$  is clear.\\
 The characteristic ideal of $\mathcal{A}$ is the pull-back by the projection $p_2: N \times \C^k \to \C^k$ of the ideal of the reduced ideal $IS(k)$ of the surface $S(k)$ (see the corollary \ref{irreduct.} in the appendix). \\
 The point $ii)$ is completed by the following lemma:
 
 \begin{lemma}\label{point 2}
 Let $\gamma(\sigma, \eta) := \sum_{h=0} (k-h).\sigma_h.\eta_{h+1}$ and $g(\sigma, \eta) := \sum_{h=1}^k h.\sigma_h.\eta_h$. Then defined  the following ideals in $\mathcal{O}_N[\eta]$, where $I_1 := (p_2)^*(IS(k))$:
 $$ I_2 := I_1 +  (\gamma) \quad {\rm and} \quad I_3 := I_2 + ( g) . $$
 Then $I_2$ is  reduced   and $g$ does not vanish on any non empty  open set of the analytic subset  $(N \times S(k))\cap \{\gamma = 0\}$.
 \end{lemma}
 
 \parag{Proof} To see that $I_2$ is reduced, as $N \times S(k)$ is normal, it is enough to prove that $\{\gamma = 0\} $ defined a reduced and irreducible hyper-surface in $N \times S(k)$. Looking at the chart on the dense open set $\eta_k\not= 0$ of $N \times S(k)$ which is given by the map  $(\sigma, \eta) \mapsto (\sigma, -\eta_{k-1}/\eta_k, \eta_k) \in N \times \C \times \C^*$ (see the paragraph  3.2) we find that $\gamma$ is given in this chart by
 $$ \gamma(\sigma, \eta) =  (-1)^{k-1}.P'_\sigma(z).\eta_k \quad {\rm where} \quad  z := -\eta_{k-1}/\eta_k     $$
using the fact that $\eta_h = (-z)^{k-h}.\eta_k$ in this chart. This gives the fact that  $\{\gamma = 0\} $ is reduced and irreducible in $N \times S(k)$.\\
The computation of $g$ in the same chart gives that 
$$ g(\sigma, \eta) = (-1)^k.z.P'_\sigma(z).\eta_k - (-1)^k.k.P_\sigma(z).\eta_k .$$
and this proves that $g$ does not vanishes identically on any non zero open set in $(N \times S(k)) \cap \{\gamma = 0\}$ because
 $$(N \times S(k)) \cap \{\gamma = 0\} \cap \{ g = 0\} \subset Z \cap \{\gamma = 0\} $$
which has dimension $k$, so  co-dimension $2$ in $N \times S(k)$.$\hfill \blacksquare$

\parag{End of proof of \ref{injective}} The point $iii)$ is consequence of the following easy  formulas:
\begin{align*}
& A_{p, q}.U_{-1} =  U_{-1}.A_{p, q} - (k-p-1).A_{p+1, q} - (k-q).A_{p, q+1} \\  
& A_{p, q}.U_0 = U_0.A_{p, q} - (p+q).A_{p, q} \\
& U_{-1}.(U_0 - 1) = U_0.U_{-1}.
\end{align*}
The points $iv)$ and $v)$ are obvious because a non zero germ of section of $\mathcal{O}_N[\eta]$ which is homogeneous of degree $1$ in $\eta$ does not vanishes of $N \times S(k)$.$\hfill \blacksquare$\\

Recall that in $\mathcal{W}(q) := \oplus_{p=0}^q \mathcal{W}_p$ we have, for each $\beta \in \mathbb{N}^k$ with $\vert \beta\vert = q-1$ and $w(\beta) = r-1$ (compare with formulas $(13)$ and $(14)$, but here $w(\beta) = r-1$)

\begin{equation}
[\partial^\beta.(U_0 - 1)] = \sum_{h=1}^k  h.\sigma_h.y_{q, r+h-1} + (r - 2).y_{q-1, r -1}.
\end{equation}
and
\begin{equation}
[\partial^\beta U_{-1} ]= \sum_{h=0}^k (k-h).\sigma_h.y_{q, r+h}  + (k.(q-1) - r+1).y_{q-1, r}
\end{equation}

Now note $\beta^+ $ a multi-index with $\vert \beta^+\vert = q-1$ and $w(\beta^+) = r$, when $r \not= k.(q-1)+1$ and $\beta^+ = 0$ for $r = k.(q-1)+1$.\\
Then for $r \not= k.(q-1) + 1$ we have 
\begin{equation*}
\partial^{\beta^+}.(U_0 -1) = \sum_{h=0}^k h.\sigma_h.y_{q, r+h} + (r-1).y_{q-1, r} 
\end{equation*}
with the convention $\sigma_0 \equiv 1$ and $\partial^{\beta^+}.(U_0 -1) = 0$ for $r = k.(q-1)+1$.

Then define for $q \geq 2$ and  $r \in [q, k.(q-1)]$  the following elements in $\mathcal{W}_q$: 
\begin{equation}
\tilde{\theta}_{q, r} := (r-1).\partial^{\beta} U_{-1}- (k.(q-1) - r+1).\partial^{\beta^+}.(U_0 - 1)
\end{equation}

This gives
\begin{equation}
\tilde{\theta}_{q, r}= k.\sum_{h=0}^k \big( (r-1) - h.(q-1) \big).\sigma_h.y_{q, r+h} 
\end{equation}
Remark that for $r = k.(q-1)+1$ and $h = k$ the vector $y_{q, r+k}$ is not defined in $\mathcal{W}_q$ and we cannot use the formula $(27)$ to define $\tilde{\theta}_{q, k.(q-1)+1}$.  But with our convention $\partial^{\beta^+}.(U_0 -1) = 0$ for $r = k.(q-1)+1$, we  define the vector   \begin{equation}
\tilde{\theta}_{q, k.(q-1)+1} := k.(q-1).\partial^{q-1}_k.U_{-1} = k.(q-1).\sum_{h= 0}^{k-1} (k-h).\sigma_h.y_{q, k.(q-1) + 1 + h}
\end{equation}
 which is  in $\mathcal{W}_q$.\\
Then,  for $q \geq 2$,  let $\tilde{\Theta}_q \subset \mathcal{W}_q$ be the sub$-\mathcal{O}_N-$module generated by the $\tilde{\theta}_{q, r}, r \in [q, k.(q-1]+1]$. Of course $\tilde{\Theta}_q$ is in the kernel of the $\mathcal{O}_N-$linear map
$$ L_q : \mathcal{W}_q \to \mathcal{N}_1^\square(q)$$
induced by the quotient map $\mathcal{W} \to \mathcal{N}^\square_1$.\\
For $q = 1$ define $\tilde{\Theta}_1 := \mathcal{O}_N.U_{-1}$ and $V_1 := \oplus_{h=2}^k \mathcal{O}_N.[\partial_h]$ (where $y_{1, h} := [\partial_h]$ in $\mathcal{W}_1$).

\begin{lemma}\label{direct}
For each  $q \geq 1$ we have a direct sum decomposition $\mathcal{W}_q = \tilde{\Theta}_q \oplus V_q$ where $V_q$ is the $\mathcal{O}_N-$sub-module with basis $y_{q, r}$ with $r \in [k.(q-1)+2, k.q]$.
\end{lemma}

\parag{proof}  For $q = 1$ our assertion is clear. For $q \geq 2$ ( and so $r \geq 2$) the difference  $\tilde{\theta}_{q, r} - k.(r-1).y_{q, r} $ is a $\C[\sigma]-$linear combination of the $y_{q, s}$ for $s \geq r+1$, so the matrix of the vectors $\tilde{\theta}_{q, r}$ for $r \in [q, k.(q-1)+1]$ and $y_{q, s}, s \in [k.(q-1)+2, k.q]$ is triangular in the basis $y_{q, t}, t \in [q, k.q]$ of $\mathcal{W}_q$ with determinant  equal to
 $$k^{(k-1)(q-1)+1}. \prod_{r = q}^{k.(q-1)+1} (r-1) = k^{(k-1)(q-1)+1}.\frac{(k.(q-1))!}{(q-2)!} $$ 
which is a positive integer.  $\hfill \blacksquare $\\

\begin{lemma}\label{surject.}
For each $q \geq 1$ the map $l_q : V_q \to \mathcal{N}^\square_1(q)$ induced by $L_q$ is bijective.
\end{lemma}

\parag{proof}
We shall prove this lemma by induction on $q \geq 1$. First remark that  the map $l_1 : V_1 \to \mathcal{N}^\square_1(1)$ is surjective (in fact an isomorphism of free rank  $(k-1)$ $\mathcal{O}_N-$modules)  because $1 = \sum_{h=1}^k h.\sigma_h.\partial_h$ and $k.\partial_1 = - \sum_{h=1}^{k-1} (k-h).\sigma_h.\partial_{h+1}$ in $\mathcal{N}^\square_1$.\\
So let $q \geq 2$ and assume that $l_{q-1} :V_{q-1} \to \mathcal{N}^\square_1(q-1)$ is surjective. Then $\mathcal{N}^\square_1(q-1)$ is contained in the image of $L_q$ because for $r \in [k.(q-2)+2, k.(q-1)]$ the relation $(26)$ shows that the image of $y_{q-1, r}$ by $l_{q-1}$ is in the image of $L_q$.\\
Then remark that $L_q$ induces a surjective  map onto  the quotient $\mathcal{N}_1^\square(q)\big/\mathcal{N}_1^\square(q-1)$  and that $\tilde{\Theta}_q$ is in  the kernel of $L_q$. So $l_q$ is surjective on $\mathcal{N}^\square_1(q)$. So we have a surjective map $l_q$ of the rank $k-1$ free $\mathcal{O}_N-$module $V_q$ onto  $\mathcal{N}^\square_1(q)$ and an injective map of the rank $k-1$ free $\mathcal{O}_N-$module $\mathcal{N}_1^\square(1)$ into $\mathcal{N}_1^\square(q)$. The lemma \ref{algebre} gives that $l_q$ is bijective.$\hfill \blacksquare$\\

\begin{thm}\label{minimality}
The restriction of $\mathcal{N}^\square_1(1)$ to the Zariski open set $\{\Delta(\sigma) \not= 0\}$ is a rank $(k-1)$ free\footnote{isomorphic to $\oplus_{h=2}^k \mathcal{O}_N.\partial_h$.}  $\mathcal{O}_N-$module  with a simple pole meromorphic connection along $\{\Delta = 0\}$ given by the inclusion $\Delta(\sigma).\mathcal{N}^\square_1(2) \subset \mathcal{N}^\square_1(1)$ (see the lemma \ref{manquant} below). Its sheaf of horizontal sections is locally generated by $z_i - z_j$ where $z_h,  h \in [1, k],$ are local branches of the multivalued function $z(\sigma)$ defined by  $P_{\sigma}(z(\sigma)) \equiv 0$. The $D_N-$module $\mathcal{N}^\square_1$ is the minimal extension on $N$ of this vector bundle with its integrable regular meromorphic connection. So it is a simple $D_N-$module.
\end{thm}

\parag{Proof} 
The lemma \ref{surject.} gives that the map $l_q : V_q \to \mathcal{N}^\square_1(q)$ is a isomorphism of $\mathcal{O}_N-$modules for each $q \geq 1$ and the proposition \ref{injective} implies that $\mathcal{N}^\square_1$ is the union of the sheaves $\mathcal{N}^\square_1(q), q \geq 1$. So the $D_N-$module $\mathcal{N}^\square_1$ has no $\mathcal{O}_N-$torsion.\\
Thanks to lemma \ref{manquant}  below we have the inclusion $\Delta.\mathcal{N}^\square_1(2) \subset \mathcal{N}^\square_1(1)$. This implies that $\mathcal{N}^\square_1(1) \simeq \mathcal{O}_N^{k-1}$ has an integrable meromorphic connection $\nabla_1$ with a simple pole along $\{\Delta(\sigma) = 0\}$ on $N$. The fact that $\mathcal{K}_1 =  \mathcal{J}_1 + D_N.U_{-1}$ implies that the horizontal sections of $\mathcal{N}^\square_1(1)$ are trace functions (see \cite{[B.19]}) which are homogeneous of degree 1 and killed by $U_{-1}$. So they are $\C-$linear combinations of $z_1(\sigma), \dots, z_k(\sigma)$, the local branches of the multivalued function $z(\sigma)$ on $N$ defined by $P_\sigma(z(\sigma)) = 0$.\\
The condition for $\sum_{h=1}^k a_h.z_h(\sigma), a_h \in \C$ to be killed by $U_{-1} \simeq \sum_{h=1}^k \frac{\partial}{\partial z_h} $ is given by $\sum_{h=1}^k a_h = 0$ and then the horizontal sections are linear combinations of the differences  $z_i - z_j, \ i, j \in [1, k]$. A basis of horizontal sections is given, for instance,  by $z_2(\sigma) - \sigma_1/k, \dots, z_k(\sigma) - \sigma_1/k$ (note that $\sum_{j=1}^k (z_j(\sigma) - \sigma_1/k) \equiv 0$).\\
The $D_N-$module  $\mathcal{N}^\square_1$ has neither $\mathcal{O}_N-$torsion nor $\mathcal{O}_N-$co-torsion because its characteristic variety  is the union of $N \times \{0\}$ and $X$ and,  thanks to the theorem \ref{co-torsion delta},  it has neither $\Delta-$torsion nor $\Delta-$co-torsion as a quotient of $\mathcal{M}$. So $\mathcal{N}^\square_1$ is the minimal extension of the meromorphic connection $\big(\mathcal{N}^\square_1(1), \nabla_1\big)$ and it is a simple $D_N-$module because the monodromy representation of the local system of horizontal sections of $(\mathcal{N}^\square_1(1), \nabla_1)$ is irreducible. $\hfill \blacksquare$\\

\begin{lemma}\label{manquant}
We have $\Delta.\mathcal{N}^\square_1(2) \subset \mathcal{N}^\square_1(1)$.
\end{lemma}

\parag{proof} In $\mathcal{W}_2\big/\mathcal{W}_1$ the $2k-1$ vectors induced by $\partial_j.(U_0 - 1), j \in [2, k]$ and $\partial_h.U_{-1}, h \in [1, k]$ are given in the basis $y_{2, r}, r \in [2, 2k]$ of this free $\mathcal{O}_N-$module by the relations
\begin{align*}
&A _j := \partial_j.(U_0 - 1) = \sum_{p=1}^k p.\sigma_p.y_{2, j+p} \\
& B_h := \partial_h.U_{-1} = \sum_{p=0}^{k-1} (k-p).\sigma_p.y_{2, h+p+1}
\end{align*}
with the convention $\sigma_0 = 1$. \\
Put $\tilde{P}_\sigma(z) := \sum_{p=0}^k \sigma_p.z^{k-p}$ and $y_{2, k+p} = z^{k-p}$.\\
 Then $B_k = \tilde{P}_\sigma'(z)$ and $A_k = z.\tilde{P}_\sigma'(z) - k.\tilde{P}_\sigma(z)$. So the resultant of $A_k$ and $B_k$ is equal to $(-k)^{k-1}.\Delta(\sigma)$.
 The  determinant of the vectors $A_j, j \in [2, k]$ and $B_h, h \in [1, k]$ in the basis $y_{2, r}, r\in [2, 2.k]$ of $\mathcal{W}_2 \simeq \mathcal{W}_2\big/\mathcal{W}_1 $ is then equal to $(-k)^{k-1}.\Delta(\sigma)$ (compare with  lemma \ref{determinant}).$\hfill \blacksquare$\\

\subsubsection{The structure theorem for $\mathcal{N}_1$}

We first examine the case $p = 1$. As already explained in the beginning of this section this will be enough to describe the structure of $\mathcal{N}_p$ for any $p \in \mathbb{N}^*$.\\
The torsion sub-module of $\mathcal{N}_1$ is described by the following result. Remark that we already know from the previous theorem \ref{minimality} that the torsion sub-module of $\mathcal{N}_1$ is contained in the image of $\mathcal{T}_1 : \mathcal{N}_0 \to \mathcal{N}_1$ as $\mathcal{N}^\square_1 = \mathcal{N}_1\big/Im(\mathcal{T}_1)$ has no torsion.

\begin{prop}\label{torsion 1.1}
There exists a injective $D_N-$linear map $\chi : \underline{H}^1_{[\sigma_k = 0]}(\mathcal{O}_N) \to \mathcal{N}_1$ which sends the class $[1/\sigma_k]$ in $\underline{H}^1_{[\sigma_k = 0]}(\mathcal{O}_N)$ to the class $[\partial_k.U_{-1}]$ in $\mathcal{N}_1$. Its image is the torsion sub-module of $\mathcal{N}_1$.
\end{prop}

\parag{Proof} Note first that $\underline{H}^1_{[\sigma_k = 0]}(\mathcal{O}_N)$ is given by $D_N\big/\big(\sum_{h=1}^{k-1} D_N.\partial_h + D_N.\sigma_k\big)$ as the annihilator of $[1/\sigma_k]$ is generated by $\partial_h, h \in [1, k-1]$ and $\sigma_k$. To show that $\chi$ exists it is enough to show that $\partial_h, h \in [1, k-1]$ and $\sigma_k$ annihilate the class  $[\partial_k.U_{-1}]$ in $\mathcal{N}_1$.
The fact that $\partial_h.[\partial_k.U_{-1}] = 0$ in $\mathcal{N}_1$ for $h \in [1, k-1]$ is a direct consequence of the formulas $(E_h), h \in [2, k]$ which give $[\partial_h.U_{-1}] = 0$ in $\mathcal{N}_1$. Then the formula $(E_1)$ gives the vanishing of  the class of $E.U_{-1} = \sum_{h=1}^k \sigma_h.\partial_h.U_{-1}$ in $\mathcal{N}_1$. So we obtain that $\sigma_k.[\partial_k.U_{-1}]$ vanishes in $\mathcal{N}_1$ and $\chi$ is well defined. Moreover, as $\underline{H}^1_{[\sigma_k = 0]}(\mathcal{O})$ is a $D_N-$module with support in $\{\sigma_k = 0 \}$, its image is contained in the torsion sub-module in $\mathcal{N}_1$.\\
Note that we know that the torsion in $\mathcal{N}_1$ is only $\sigma_k-$torsion thanks to the corollary \ref{caract. cycles} and the theorem \ref{co-torsion delta}.\\
To prove the injectivity of $\chi$, assume that the kernel of $\chi$  is not $0$ and consider an element  $K := \sum_{p=1}^m f_p.\partial_k^p[1/\sigma_k]$ in this kernel with $f_p \in \mathcal{O}_N\big/(\sigma_k)$ and with $m$ minimal. Then we have $0 = \chi(K) = [\sum_{p=1}^m f_p.\partial_k^{p+1}.U_{-1}]$ in $\mathcal{N}_1$. So $f_m.\eta_k^{m+1}.\gamma$ is the symbol of an element in $\mathcal{J}_1$. The lemma \ref{useful 2} then implies that $f_m$ is in $\sigma_k.\mathcal{O}_N$ contradicting the minimality of $m$. So $\chi$ is injective.\\
To complete the proof we have to show that if $P$ induces a torsion class in $\mathcal{N}_1$ then there exists $Q \in D_N$ such that $P - Q.\partial_k.U_{-1} $ is in $\mathcal{J}_1$. As we already know ( because $\mathcal{N}^\square_1$ has no torsion) that there exists $P_1 \in D_N$ such that $\mathcal{T}_1(P_1) = [P_1.U_{-1}] = [P]$ in $\mathcal{N}_1$ and as we know that $\partial_h.U_{-1} = 0$ for each $h \in [1, k-1]$ we may assume that $P_1$ is in $\mathcal{O}_N[\partial_k]$. But $\partial_k^n.U_{-1}$ is torsion in $\mathcal{N}_1$ for $n \geq 1$ because $\partial_k.U_{-1}$ is torsion (see above). So the only point to prove is that if  $f.U_{-1}$ is torsion in $\mathcal{N}_1$ for some $f \in \mathcal{O}_N$  then $f = 0$. This a consequence of the  following lemma. $\hfill \blacksquare$

\begin{lemma}\label{useful 0}
Le class of $U_{-1}$ is not  in the $\sigma_k-$torsion of $\mathcal{N}_1$.
\end{lemma}

\parag{Proof} Assume that $\sigma_k^n.U_{-1}$  is in $\mathcal{J}_1$ for some $n \in \mathbb{N}$. Then choose $n$ minimal with this property and compute
$$ \partial_k\sigma_k^n.U_{-1} = n.\sigma_k^{n-1}.U_{-1} + \sigma_k^n.\partial_k.U_{-1} \in \mathcal{J}_1. $$
 As $\sigma_k.\partial_k.U_{-1} $  is in $\mathcal{J}_1$ (see above) we obtain that $n = 0$ by minimality of $n$. But $U_{-1}$ is not in $\mathcal{J}_1$ because its symbol $\gamma(\sigma, \eta)$ restricted to the co-normal $C$ to the hyper-surface  $\{\sigma_k = 0 \}$ is equal to $\sigma_{k-1}$ which does not vanish identically on $C$. And $C$ is a component of the characteristic variety of $\mathcal{N}_1$ (see paragraph 3.2). This concludes the proof.$\hfill \blacksquare$\\

\begin{thm}\label{fin 1}
 The diagram below describes the structure of $\mathcal{N}_1$,  where $\mathscr{T}$ is the torsion sub-module of $\mathcal{N}_1$, where $\varphi_1: \mathcal{N}_1 \to \mathcal{O}_N $ is the $D_N-$linear map defined by $\varphi_1(1) = \sigma_1$ and where the isomorphism   $\chi :\mathscr{T} \simeq \underline{H}^1_{[\sigma_k = 0]}(\mathcal{O})$ is  defined by sending  $1/\sigma_k $ to $[\partial_k.U_{-1}] $.\\
  The $D_N-$modules $\mathscr{T} \simeq \underline{H}^1_{[\sigma_k = 0]}(\mathcal{O}_N)$ and $\mathcal{N}^\square_1$ are simple $D_N-$modules.\\
 Moreover we have the direct  sum decomposition of left $D_N-$modules:
 $$ \mathcal{N}_1\big/\mathscr{T} = Im(\mathcal{T}_1)\big/\mathscr{T} \oplus \mathcal{N}_1^*\big/\mathscr{T} = \mathcal{O}_N.[U_{-1}] \oplus \mathcal{N}^\square_1. $$
The following commutative diagram of left $D_N-$modules has exact lines and columns where the maps $i$ and $e$ are defined  by $i([U_{-1}]) = [U_{-1}]$ and $e([U_{-1}]) = 1/k$:
 $$ \xymatrix{\quad & 0 \ar[d] & 0 \ar[d] & \quad & \quad \\ 0 \ar[r] & \mathscr{T} \ar[d] \ar[r] & Im(\mathcal{T}_1) \ar[d] \ar[r]^{i} & \mathcal{O}_N.[U_{-1}] \ar[d]^{e}_\simeq \ar[r] & 0 \\
 0 \ar[r] & \mathcal{N}_1^* \ar[r] \ar[d] & \mathcal{N}_1 \ar[d] \ar[r]^{\varphi_1} & \mathcal{O}_N \ar[r] & 0 \\ \quad & \mathcal{N}_1^*\big/\mathscr{T} \ar[d] \ar[r]^\theta & \mathcal{N}^\square_1 \ar[d] & \quad  \\
  \quad &  0 & 0 & \quad & \quad} $$
   \end{thm}

\parag{Proof} Note first that the quotient by the torsion sub-module $\mathscr{T}$ (which is the image of $D_N.\partial_k.U_{-1}$ in $\mathcal{N}_1$ ; see proposition \ref{torsion 1.1}) of the image by $\mathcal{T}_1$  of  $D_N.U_{-1}$ in $\mathcal{N}_1$ is isomorphic to $\mathcal{O}_N$ because its generator $[U_{-1}]$ is killed by $\partial_h, \forall h \in [1, k]$ (see formulas $(E_h), h \in [1, k]$)  and this quotient has no torsion because $\mathscr{T}$ is also the torsion sub-module of the image of $D_N.\partial_k.U_{-1}$ in $\mathcal{N}_1$. This gives the exactness of the first line. The exactness of the second line  and of the columns are clear.\\
Note also that $\varphi_1(U_{-1}) = k$ so the upper right square commutes. The commutations of the other squares are obvious.\\
To show that the map $\theta$ is well defined and is an isomorphism is a simple exercice in diagram chasing which is left to the reader.\\
The direct sum decomposition of $\mathcal{N}_1\big/\mathscr{T}$  is given by the left $D_N-$linear map
 $$r :  \mathcal{N}_1\big/\mathscr{T} \to Im(\mathcal{T}_1)\big/\mathscr{T}$$
 constructed as follows:\\
Note first that $\varphi_1(U_1) = k$. For $[P] \in \mathcal{N}_1$ let $f := \varphi_1([P])$. Then we define
 $$r([P]) := [(f/k).U_{-1}] \in Im(\mathcal{T}_1)\big/\mathscr{T}.$$
 As $\mathscr{T}$ is in the kernel of $\varphi_1$, this map is well defined on $\mathcal{N}_1\big/\mathscr{T}$ and  $[P] - r([P])$  is in $ker(\varphi_1) = \mathcal{N}_1^*$ and defines a class in $\mathcal{N}_1^*\big/\mathscr{T}$. Remark that the lemma \ref{useful 0} shows that the kernel  of $r$ is equal to $\mathcal{N}_1^* \big/\mathscr{T}$ because $\varphi_1$ is injective on $\mathcal{O}_N.[U_{-1}]$. This gives the desired splitting, as $r$ induces the identity on $Im(\mathcal{T}_1)\big/\mathscr{T}$.$\hfill \blacksquare$\\

\subsection{The structure of $\mathcal{N}_0$}
Define $\mathcal{N}_0^*$ as the kernel of the $D_N-$linear map $\varphi_0 : \mathcal{N}_0 \to \mathcal{O}_N$ given by $\varphi_0(1) = 1$. The sub-module $\mathcal{N}_0^*$ is generated by $\partial_1, \dots, \partial_k$. We shall show that it contains a copy of $\mathcal{O}_N$.

\begin{prop}\label{Kernel}
The kernel of the $D_N-$linear map $\mathcal{T}_1 : \mathcal{N}_0 \to \mathcal{N}_1$ given by the right multiplication by $U_{-1}$ is $D_N.U_1$ which is contained in $\mathcal{N}_0^*$ and the quotient $\mathcal{N}_0^*\big/D_N.U_1$ is isomorphic to the $D_N-$module $\underline{H}_{[\sigma_k = 0 ]}^1(\mathcal{O}_N)$.
\end{prop}

The proof of this proposition will used the following results from \cite{[B.19]} proposition 5.2.1.

\begin{prop}\label{citation}
For each $m \in \mathbb{Z}, m \geq -k+1$ and for each $\sigma \in N$ such that $\Delta(\sigma) \not= 0$ define
\begin{equation}
DN_{m}(\sigma) := \sum_{P_{\sigma}(x_{j}) = 0} \frac{x_{j}^{m+k-1}}{P'_{\sigma}(x_{j})} 
\end{equation}
Each $DN_{m}$ is the restriction to the open set $\{\Delta(\sigma) \not= 0 \}$ of a polynomial of (pure) weight $m$  in $\C[\sigma_{1}, \dots, \sigma_{k}]$ and the following properties are satisfied:
\begin{enumerate}[i)]
\item For $m \in [-k+1, -1], DN_{m} = 0 $ and $DN_0 = 1$.
\item For each $m \geq 1, \quad \sum_{h=0}^{k} (-1)^{h}.\sigma_{h}.DN_{m-h} = 0 $ with the convention $\sigma_{0} \equiv 1$.
\item For each $h \in [1, k]$ and  each $m \geq 0$ we have
\begin{equation*}
\partial_h N_{m} = (-1)^{h-1}.m.DN_{m-h} .\qquad \qquad \qquad  \qquad \qquad \qquad \qquad \qquad \qquad  \hfill \blacksquare
\end{equation*}
\end{enumerate}
\end{prop}

We shall use also the following lemma.

\begin{lemma}\label{25/7}
For any $h \in [2, k]$ we have
\begin{equation*}
 \partial_hU_1 + \partial_{h-1}.(U_0+1)  \in \mathcal{I} .\tag{$F_h$}
 \end{equation*}
 Moreover we have also
\begin{equation*}
 \partial_1.U_1 - E.(U_0+1) \in \mathcal{I} . \tag{$F_1$}
 \end{equation*}
So we have $\partial_h.U_1 = -\partial_h $ for $h \in [2, k]$ and $\partial_1U_1 = E$ in $\mathcal{N}_0^*$ and  $\partial_h.U_1 = 0$ in $\mathcal{N}_{-1}$ for any $h \in [1, k]$.
\end{lemma}

\parag{proof} Thanks to the theorem which characterizes trace functions in \cite{[B.19]} it is enough to prove that for each $m \in \mathbb{N}$ we have $\partial_h.U_1[N_m] = -\partial_{h-1}.(U_0 + 1)[N_m]$ for $h \in [2, k]$ and $\partial_1U_1(N_m) = E.(U_0+1)[N_m]$ for all $m \in \mathbb{N}$. This is consequence of the following formulas which use the results of \cite{[B.19]} recalled in the proposition \ref{citation} and the formula $U_p[N_m] = m.N_{m+p} $ which is valid for each $m \in \mathbb{N}$ and each integer $p \geq -1$ because $U_p$ is the image by the quotient map $q : M = \C^k \to \C^k\big/\mathfrak{S}_k = N$ of the vector field $\sum_{j=1}^k z_j^{p+1}.\frac{\partial}{\partial z_j}$:
\begin{align*}
& \partial_h.U_1[N_m] = \partial_h[m.N_{m+1}] = (-1)^{h-1}.m.(m+1).DN_{m+1-h} \quad \forall h \in [1, k] \quad \forall m \in \mathbb{N} \\
& \partial_{h-1}.(U_0 + 1)[N_m] = \partial_{h-1}[(m+1).N_m] = (-1)^h.m.(m+1).DN_{m-h+1} \quad \forall h \in [2, k]
\end{align*}
proving the formulas $(F_h)$ for $h \in [2, k]$.\\
Now  $\sigma_h.\partial_h[N_m] = (-1)^{h-1}.m.\sigma_h.DN_{m-h}$ gives
$$ E.(U_0+1)[N_m] = m.(m+1).\sum_{h=1}^k (-1)^{h-1}.\sigma_h.DN_{m-h} = m.(m+1).DN_m $$
because for $m \geq 1$ we have $\sum_{h=0}^k (-1)^h.\sigma_h.DN_{m-h} = 0$ (see  $ii)$ above) and also $E[1] = 0$ for $m = 0$.
But for $h = 1$ we have  
$$\partial_1.U_1[N_m] = \partial_1[m.N_{m+1}] = m.(m+1).DN_{m} \quad \forall m \in \mathbb{N}.$$
This  gives the formula $(F_1)$ .$\hfill \blacksquare$\\

\parag{Proof of the proposition \ref{Kernel}}Remark first that $U_1$ is in $\mathcal{N}_0^*$ and, thanks to the previous lemma, that $D_N.U_1$ contains $\partial_1, \dots, \partial_{k-1}$ and $\sigma_k.\partial_k$.  Define the sub-$D_N-$module $S := \sum_{h=1}^{k-1} D_N.\partial_h$. Then we have a natural surjective $D_N-$linear map $ \alpha : S + D_N\partial_k\big/(S + D_N.\sigma_k.\partial_k) \to \mathcal{N}_0^*\big/D_N.U_1$. But we have
$$ S + D_N\partial_k\big/(S + D_N.\sigma_k.\partial_k) \simeq D_N.\partial_k\big/(S \cap D_N.\sigma_k.\partial_k) \simeq D_N\big/(S + D_N.\sigma_k) $$
thanks to the equality $S \cap D_N.\sigma_k.\partial_k = (S + D_N.\sigma_k).\partial_k$. Moreover, the $D_N-$module
$$D_N\big/(S + D_N.\sigma_k) \simeq \underline{H}_{[\sigma_k = 0 ]}^1(\mathcal{O}_N)$$
is simple, so $\alpha$ must be an isomorphism.$\hfill \blacksquare$\\

\begin{thm}\label{structure 0}
Define $\mathcal{N}_0^{\square} := D_N.U_1 \subset \mathcal{N}_0^*$. Then $\mathcal{N}_0^{\square}$ is simple and isomorphic to $\mathcal{N}^\square_1$ via the map induced by the map  $\square.U_1 : \mathcal{N}^\square_1 \to \mathcal{N}_0^*$, and the quotient $\mathcal{N}_0^*\big/\mathcal{N}_0^{\square}$ is isomorphic to $\underline{H}^1_{[\sigma_k = 0]}(\mathcal{O}_N)$.
\end{thm}

\parag{Proof}  The only point which is not already proved above is the fact that right multiplication by $U_1$,  $\square.U_1: \mathcal{N}_1 \to \mathcal{N}_0$ has its image in $\mathcal{N}_0^*$ and induces an isomorphism of 
$\mathcal{N}^\square_1$ to $\mathcal{N}^\square_0 = D_N.U_1 \subset \mathcal{N}_0^*$. But $U_1$ is in $\mathcal{N}_0^*$ so the first assertion is clear. This map vanishes on  the image of $\mathcal{T}_1$ because $U_{-1}.U_1 = (U_0+1).U_0 \quad {\rm modulo} \ \mathcal{I}$ (and also $\mathcal{I}.U_1 \subset  \mathcal{I}$ see the theorem \ref{action}), so that right multiplication by $U_1$  induces a map which is clearly surjective. As $\mathcal{N}^\square_1$ is simple, this surjective map is an isomorphism.\\
We have the following commutative diagram with exact lines and columns describing the structure of the $D_N-$module  $\mathcal{N}_0$:

$$ \xymatrix{\quad &\quad &0 \ar[d] &0 \ar[d] &\quad & \quad &\quad \\ 
0 \ar[r] &\mathcal{N}_0^{\square} \ar[r] & \mathcal{N}_0^* \ar[d] \ar[r] & \underline{H}_{[\sigma_k = 0 ]}^1(\mathcal{O}_N) \ar[d] \ar[r] & 0 & \quad  &\quad \\
0 \ar[r] & \mathcal{N}_0^*  \ar[r]  \ar[ur]^= & \mathcal{N}_0 \ar[d]^{\varphi_0} \ar[r]^{\mathcal{T}_1} & Im(\mathcal{T}_1) \ar[d] \ar[r] & \mathcal{N}_1 \ar[r] & \mathcal{N}^\square_1 \ar[r] \ar@/_11pc/[llllu]_{\simeq \quad \square.U_1} & 0 \\
\quad & \quad & \mathcal{O}_N \ar[d] \ar[r]^{\square.U_{-1}} & \mathcal{O}_N.[U_{-1}] \ar[d] & \quad & \quad & \quad \\
\quad & \quad & 0 & 0 & \quad & \quad & \quad } $$

\bigskip

where the map $\square.U_1 : \mathcal{N}_1 \to \mathcal{N}_0$ given by right multiplication by $U_1$  induces an isomorphism $\mathcal{N}^\square_1 \to \mathcal{N}_0^{\square}$, showing that $\mathcal{N}_0^{\square}$ is a simple $D_N-$module.$\hfill \blacksquare$ \\

Note that the local horizontal basis of $\mathcal{N}_0^{\square}$ on the open set $\{\Delta(\sigma).\sigma_k \not= 0\}$ is  (locally) generated by the $(Log \, z_i - Log\, z_j)$ and their images by the isomorphism induced by $U_1$ are the $(z_i - z_j)$ which generates a local horizontal basis of $\mathcal{N}^\square_1$ on the open set   $\{\Delta(\sigma).\sigma_k \not= 0\}$.

\subsection{The structure of $\mathcal{N}_p$ for $p \in -\mathbb{N}^*$}

Again it is enough, thanks to the lemma \ref{isom 2}, to describe the structure of $\mathcal{N}_{-1}$. Define $\mathcal{N}_{-1}^*$ as the kernel of the $D_N-$linear map
$$\varphi_{-1} : \mathcal{N}_{-1} \to \mathcal{O}_N[\sigma_k^{-1}] $$
which is given by $\varphi_{-1}(1) = \sigma_{k-1}\big/\sigma_k $. This map is well defined because the meromorphic function $\sigma_{k-1}\big/\sigma_k = \sum_{j=1}^k 1\big/z_j$ is a local trace function of the open set $\{\sigma_k \not= 0\}$ and so  it is killed by $\mathcal{I}$ everywhere as the $D_N-$module $\mathcal{O}_N[\sigma_k^{-1}] $ has no torsion. Moreover we have, still on the open set $\{\sigma_k \not= 0\}$:
 $$U_0(\sigma_{k-1}\big/\sigma_k) =( k-1)\sigma_{k-1}\big/\sigma_k -k\sigma_{k-1}\big/\sigma_k = -\sigma_{k-1}\big/\sigma_k.$$
  So $U_0 +1$ is also in the annihilator of $\sigma_{k-1}\big/\sigma_k $ in $\mathcal{O}_N[\sigma_k^{-1}] $. Therefore the map $\varphi_{-1}$ is well defined. It is surjective because $\varphi_{-1}(\partial_{k-1}) = 1/\sigma_k$.\\

\begin{lemma}\label{symb. U_1}
The symbol of $U_1$ does not vanish identically on $X$ for each integer $k \geq 2$.
\end{lemma}

\parag{Proof} We have $\eta_h = (-z)^{k-h}.\eta_k$ on $X$ where $z = -\eta_{k-1}/\eta_k$ is the double root of $P_\sigma$ at the generic point of $\Delta$ (recall that $X$ is the closure of the graph of the meromorphic map $\Delta \dashrightarrow \C$ given by the double root of $P_\sigma$ at the generic point of $\Delta$). As $U_1 = \sum_{h=1}^k (\sigma_1\sigma_h - (h+1).\sigma_{h+1}).\partial_h $ (with the convention $\sigma_{k+1} = 0$), we obtain
\begin{align*}
& s(U_1) = (-1)^k\eta_k.\sum_{h=1}^k (\sigma_1\sigma_h - (h+1).\sigma_{h+1}).(-1)^h.z^{k-h} \\
& \quad  = (-1)^k\eta_k.\sum_{h=1}^k (-1)^h.\sigma_1\sigma_h.z^{k-h} + (-1)^k\eta_k.\sum_{p=2}^k (-1)^{p-1}.(k-p).\sigma_p.z^{k-p+1} + \\
& \qquad \qquad - (-1)^kk.\eta_k.\sum_{p=2}^k (-1)^{p-1}.\sigma_p.z^{k-p+1} 
\end{align*}
\begin{align*}
& s(U_1) = (-1)^k\eta_k.\sigma_1.\big(P_\sigma(z) - z^k\big) + (-1)^k\eta_k.z^2.\big(P'_\sigma(z) - k.z^{k-1} + (k-1).\sigma_1.z^{k-2}\big) \\ 
&  \qquad \qquad    - (-1)^kk.\eta_k.z.\big(P_\sigma(z) - z^k + \sigma_1.z^{k-1}\big) \\
& \quad  =  -(-1)^k\eta_k.\sigma_1.z^k - (-1)^kk.\eta_k.z^{k+1} + (-1)^k(k-1).\eta_k.\sigma_1.z^k + \\
& \qquad \qquad  (-1)^kk.\eta_k.z^{k+1} - (-1)^kk.\eta_k.\sigma_1.z^k \\
& \quad = -2(-1)^k\eta_k.\sigma_1.z^k
\end{align*}
as $P_\sigma(z) = P'_\sigma(z) = 0$ on $X$.$\hfill \blacksquare$\\

So $U_1$ is not zero in any $\mathcal{N}_\lambda$ for any $\lambda \in \C$.

\begin{prop}\label{serieux}
The kernel of the $D_N-$linear map $\mathcal{T}_0 : \mathcal{N}_{-1} \to \mathcal{N}_0$ given by right multiplication by $U_{-1}$ is equal to $D_N.U_1 = \mathcal{O}_N.U_1$ in $\mathcal{N}_{-1}$.
\end{prop}

\parag{Proof} Recall that the relations $(F_h), h \in [1, k]$ show the equality $D_N.U_1 = \mathcal{O}_N.U_1$ in $\mathcal{N}_{-1}$ (see the lemma \ref{25/7}). Moreover we know that $\mathcal{N}_{-1}$ has no $\mathcal{O}_N-$torsion, thanks to the theorem \ref{no torsion},  so $D_N.U_1$ is a sub-module of $\mathcal{N}_1$ which is isomorphic to $\mathcal{O}_N$ as $U_1$ is not zero in $\mathcal{N}_{-1}$ because its symbol does not vanish on $X$ (see the lemma \ref{symb. U_1} above).\\
The rest of the proof of the proposition will use the following lemmas:

\begin{lemma}\label{useful 1}
Let $a$ and $b$ be holomorphic function on an open set $U$ in $N$ such that the function $a.\gamma - b.g$ is a section on $U \times  \C^k$ which vanishes on $Z \cap (U \times  \C^k)$. Then $a$ and $b$ vanishes identically on $U$
\end{lemma}

\parag{proof} The first remark is that we have $a.\gamma = b.g$ on $U \times \C^k$ because the sheaf $p_*(I_Z)$ has no non zero section which homogeneous of degree $1$ in $\eta_1, \dots, \eta_k$ (see lemma \ref{remplace}). Then looking at the coefficients of $\eta_1$ and $\eta_2$ in the equality $a.\gamma = b.g$ gives
$$ k.a = \sigma_1.b \quad {\rm and} \quad (k-1).\sigma_1.a = 2\sigma_2.b \quad {\rm and \ so} \quad (k-1)\sigma_1^2.b = 2k.\sigma_2.b $$
which implies $b \equiv 0$ and then $a \equiv 0$ on $U$.$\hfill \blacksquare$

\begin{lemma}\label{idee}
Let $P\in D_N$  such that $P.U_{-1} = A.U_0 + Q$ with $A \in D_N$ and  $Q \in \mathcal{I}$. Then $A$ is unique modulo $\mathcal{I}  $.
\end{lemma}

\parag{proof} We have to show that $A.U_0 \in \mathcal{I}$ implies that $A$ is in $\mathcal{I} $. \\
If this is not true, let $A \in D_N \setminus \mathcal{I}$ be of minimal order such that $A.U_0 \in \mathcal{I}$. We have $s(A).g \in I_Z$ and, as $g$ is generically $\not= 0$ on $Z$  and $I_Z$ is prime (so reduced), there exists $A_1 \in \mathcal{I}$ with $A- A_1$ of order strictly less than the order of $A$. Then $(A - A_1).U_0$ is in $\mathcal{I}$.\\
 This contradicts the minimality of $A$ since $A - A_1$ cannot be in $\mathcal{I}$.$\hfill\blacksquare$

\begin{lemma}\label{idee suite}
There exists a natural $D_N-$linear map $\psi : Ker(\mathcal{T}_0) \to \mathcal{N}_1$ given by $ \psi(P) = [A] $ when $P.U_{-1} = A.U_0$ modulo $\mathcal{I}$.
\end{lemma}

\parag{Proof}Firs recall that the right multiplication by $U_{-1}$ induces a $D_N-$linear map $\mathcal{T}_0 : \mathcal{N}_{-1} \to \mathcal{N}_0$ because we have $\mathcal{I}.U_{-1} \subset \mathcal{I}$ and the relation
$ (U_0 +1).U_{-1} = U_{-1}.U_0 $.  If $P \in D_N$ induces a germ of section of $Ker(\mathcal{T}_0)$ then the previous lemma shows that if we write $P.U_{-1} = A.U_0 + Q$ with $Q \in \mathcal{I}$, the image of the germ $A$ in $D_N\big/\mathcal{I}$ is well defined. Then we have a  $D_N-$linear map $Ker(\mathcal{T}_0) \to D_N\big/\mathcal{I} = \mathcal{M}$ and after composition by the quotient maps $D_N\big/\mathcal{I} \to \mathcal{N}_1$ we obtain the desired map.$\hfill \blacksquare$

\parag{End of proof of the proposition \ref{serieux}} First remark that $U_1$ is sent to $0$ in $\mathcal{N}_1$ by the  map $\psi$ because of the relation $ U_1.U_{-1} = (U_0 -1).U_0 \quad {\rm modulo} \quad \mathcal{I}$. \\
So $\psi$ composed with the quotient map $\mathcal{N}_1 \to \mathcal{N}^\square_1$  induces a map 
$$\tilde{\psi} : Ker(\mathcal{T}_0)\big/\mathcal{O}_N.U_1 \to \mathcal{N}^\square_1.$$
 We shall prove that this map is injective and not surjective. As $\mathcal{N}^\square_1$ is simple, this will prove that $Ker(\mathcal{T}_0) = \mathcal{O}_N.U_1$ completing the proof of the proposition.\\
We shall first prove the injectivity of $\tilde{\psi}$, so the fact  that if $[P] \in Ker(\mathcal{T}_0)$ satisfies  $\psi(P) = [A]$ with $[A] = 0$  in $\mathcal{N}^\square_1$ then $[P]$ is a germ of section of the sub-sheaf $D_N.U_1$ of $Ker(\mathcal{T}_0)$.\\
Let $P \in D_N$ of minimal order such that  the class of $[P]$ in $Ker(\mathcal{T}_0)\big/D_N.U_1$ is not zero and satisfies  $\tilde{\psi}([P]) = 0$. Then we have 
$$ (P + X.(U_0 + 1) + Q_0).U_{-1} = A_0.U_0 + Q_1 \quad {\rm with} \quad Q_0, Q_1 \in \mathcal{I} \quad {\rm and} \ X \in D_N .$$
Then, thanks to the relations $(U_0 +1).U_{-1} = U_{-1}.U_0 \ {\rm modulo} \quad \mathcal{I} $ \  we obtain
$$P.U_{-1} = A.U_0 + Q_2  \quad {\rm with}  \ Q_2 := Q_1 - Q_0.U_{-1} \in \mathcal{I}  \quad {\rm and} \ A = A_0 -X.U_{-1} .$$
Then our hypothesis implies that  $A = R.(U_0 -1) + S.U_{-1} +Q_3$ \ with $Q_3 \in \mathcal{I}$. So
$$ (P - S.(U_0 +1) - R.U_1).U_{-1} = Q_3.U_0 \quad{\rm modulo} \ \mathcal{I} $$
 and $Q_3.U_0$ is again in $\mathcal{I}$.  So looking at the  symbols we find 
  $$s(P - S.(U_0 +1) - R.U_1).\gamma \in  I_Z.$$ 
  As $\gamma$  is generically $\not= 0$ on $Z$  and $I_Z$ is prime (then reduced)  we conclude that there exists $P_1 \in \mathcal{I}$ with symbol $s(P_1) = s(P - S.(U_0 +1) - R.U_1)$. So the order of $P - P_1- S.(U_0 +1) - R.U_1$ is strictly less than the order of $P$ but the class of $P - P_1 -S.(U_0+1) - R.U_1$ in $Ker(\mathcal{T}_0)\big/D_N.U_1$ is the same than the class induced by $P$. This contradict the minimality of the order of $P$;  so $Ker(\mathcal{T}_0) = D_N.U_1$ and the map $\tilde{\psi}$ is injective.\\
  To conclude it is now enough to prove that $\tilde{\psi}$ is not surjective, as explained above. So  assume that there exists  $P \in D_N$ with $P.U_{-1} = A.U_0 + Q$ with $Q \in \mathcal{I}$ and $[A  - 1] = 0 $ in $\mathcal{N}^\square_1$. This would implies that $P.U_{-1} = (1 + T.(U_0 - 1) + Y.U_{-1}).U_0 $ modulo $\mathcal{I}$ and so we obtain the equality
  $$(P -Y.(U_0+1) -T.U_1).U_{-1} = U_0 \quad {\rm  modulo} \quad \mathcal{I}.$$
  So looking at the symbols restricted to $Z$ this gives:
  $$s\big(P -Y.(U_0+1) -T.U_1\big).\gamma = g $$
  in $\mathcal{O}_Z$. By homogeneity in $\eta$ this implies that $f := s\big(P -Y.(U_0+1) -T.U_1\big) $ is the pull-back of a holomorphic function on an open set in $N$ and this gives a contradiction thanks to the lemma \ref{useful 1}. $\hfill \blacksquare$\\
  
  \begin{prop}\label{the simple -1}
  Let $\mathcal{N}_{-1}^{\square}$ be the sub-$D_N-$module of $\mathcal{N}_{-1}^*$ which is generated by $\partial_1, \dots, \partial_{k-2}$. Then $\mathcal{T}_0$ sends $\mathcal{N}_{-1}^{\square}$ onto $\mathcal{N}_0^{\square}$ and induces an isomorphism between theses two simple $D_N-$modules.
  \end{prop}
  
  \parag{Proof} As we know that $\mathcal{N}_0^{\square}$ is equal to $D_N.U_1 \subset \mathcal{N}_0^*$, we first check that the generators of $\mathcal{N}_{-1}$ have their images by $\mathcal{T}_0$ in $D_N.U_1$.\\
  For $h \in [1, k-2]$ the formulas $(E_{h+2})$ and $(F_{h+1})$ implies
  $$ \partial_{h+2}.U_1 + \partial_{h+1}.(U_0 +1) + \partial_h.U_{-1} + \partial_{h+1}.(U_0- 1) \in \mathcal{I} $$
  which implies $ \partial_h.U_{-1}  = - \partial_{h+2}.U_1$ in $\mathcal{N}_0$.\\
  Note that, as $\mathcal{N}^\square_1$ is obviously generated\footnote{In fact knowing that it is simple, it is generated by any non zero element in it.} by $\partial_1, \dots, \partial_{k-2}$ its image  by the right multiplication by $U_1$ in $\mathcal{N}_0^*$ is generated by $\partial_h.U_1, h \in [1, k-2]$ giving a direct proof of the surjectivity of  $\mathcal{T}_0 :\mathcal{N}_{-1}^{\square}
  \to D_N.U_1 = \mathcal{N}_0^{\square}$.\\
   The injectivity of this map is clear thanks to the proposition \ref{serieux} and  the fact that $\varphi_{-1}(U_1) = -k$ which implies that the sub-modules  $D_N.U_1 = \mathcal{O}_N.U_1 $ and $ Ker(\mathcal{T}_0)$ of $ \mathcal{N}_{-1}$ has an intersection reduced to $\{ 0 \}$ .$\hfill \blacksquare$\\

  \begin{prop}\label{egalite}
  The sub-module $\mathcal{N}_{-1}^{\square}$ is equal to $\mathcal{N}_{-1}^*$.
  \end{prop}
  
  \parag{Proof} By definition $\mathcal{N}_{-1}^*$ is the kernel of the map $\varphi_{-1} : \mathcal{N}_{-1} \to \mathcal{O}_N[\sigma_k^{-1}]$ which sends $[1]$ to $\sigma_{k-1}/\sigma_k$. Then $\mathcal{N}_{-1}^*$ is 
  generated by  the annihilator of $\sigma_{k-1}/\sigma_k$ in  $\mathcal{O}_N[\sigma_k^{-1}]$. So $\mathcal{N}_{-1}^*$ is generated by the class of
   $$\partial_1, \dots, \partial_{k-2}, \partial_{k-1}^2, \sigma_{k-1}.\partial_{k-1} - 1, \sigma_k.\partial_k +1.$$
   We already know that $\partial_1, \dots, \partial_{k-2}$ are in $\mathcal{N}_{-1}^{\square}$ thanks to formulas $(E_{h+2})$ and $(F_{h+1})$ for $ h \in [1, k-2]$, see proposition \ref{the simple -1}. Moreover we have $\partial_{k-1}^2 = \partial_k.\partial_{k-2} \quad {\rm modulo} \ \mathcal{I}$. So it is enough to prove that $a := \sigma_{k-1}.\partial_{k-1} - 1$ and $b =  \sigma_k.\partial_k +1$ are in $\mathcal{N}_{-1}^{\square}$. \\
   The formula $(E_1)$ gives $E.U_{-1} = -\partial_1 $ in $\mathcal{N}_0$ and the formula $(F_1)$ gives $\partial_2.U_1 = -\partial_1$ in $\mathcal{N}_0$. This implies that $\mathcal{T}_0(E) = \partial_1.U_1 \in \mathcal{N}^\square_0 $. This implies that $E$  is in $\mathcal{N}_{-1}^{\square} + ker(\mathcal{T}_0)$.\\
 So write $ E = e + f.U_1$ with $e \in \mathcal{N}_{-1}^{\square} $ and $f \in \mathcal{O}_N$ using the proposition \ref{serieux}.\\
 Now
  $$\varphi_{-1}(E) = E[\sigma_{k-1}/\sigma_k] = 0 \quad {\rm  and}  \quad \mathcal{N}^\square_{-1} \subset \mathcal{N}_1^* = ker(\varphi_{-1}).$$
 So $\varphi_{-1}(f.U_1) = f.\varphi_{-1}(U_1) = -k.f = 0$. This implies $E = e$  is in $\mathcal{N}^\square_{-1}$.  But  $a + b = E \quad {\rm modulo} \ \mathcal{N}_{-1}^{\square}$. So $a + b $ belongs to $\mathcal{N}_{-1}^{\square}$.
 We have also in $\mathcal{N}_{-1}$ :
  $$ 0 = U_0 + 1 = (k-1).(\sigma_{k-1}.\partial_{k-1} - 1) + k.(\sigma_k.\partial_k + 1) \quad {\rm modulo} \ \mathcal{N}_{-1}^{\square} $$
  and this gives $(k-1).a + k.b \in \mathcal{N}_{-1}^{\square}$, concluding the proof.$\hfill \blacksquare$\\
  
  \begin{thm}\label{structure -1}
 We have the following commutative diagram of $D_N-$module with exact lines and columns, where the $D_N-$linear map  $\varphi_{-1} : \mathcal{N}_{-1} \to \mathcal{O}_N[\sigma_k^{-1}]$ is defined by $\varphi_{-1}(1) = \sigma_{k-1}\big/\sigma_k$:
 $$\xymatrix{ \quad & 0 \ar[d] & 0 \ar[d] & 0 \ar[d] & \quad \\ 0 \ar[r] & \mathcal{N}_{-1}^* \ar[d]^{=} \ar[r] & \mathcal{N}_{-1}^*\oplus \mathcal{O}_N.U_1 \ar[d] \ar[r]^{\varphi_{-1}} & \mathcal{O}_N \ar[d] \ar[r] & 0 \\
 0 \ar[r] &  \mathcal{N}_{-1}^* \ar[r] \ar[d] & \mathcal{N}_{-1} \ar[d] \ar[r]^{\varphi_{-1}} & \mathcal{O}_N[\sigma_k^{-1}] \ar[d] \ar[r] & 0 \\ \quad & 0 \ar[r] & \mathcal{Q} \ar[r]^{\chi} \ar[d] & \underline{H}^1_{[\sigma_k = 0]}(\mathcal{O}_N) \ar[d] \ar[r] & 0 \\ \quad & \quad & 0 & 0 & \quad \\} $$
 So $\chi$ is an isomorphism. Moreover the map $\mathcal{T}_0$ induces an isomorphism of $\mathcal{N}_{-1}^{*} $ onto the simple $D_N-$module $\mathcal{N}_0^{\square} = D_N.U_1 \subset \mathcal{N}_0^*$.
 \end{thm}
 
 \parag{proof} The exactness of the first line is consequence of  the  equality  $\varphi_{-1}(U_1) = -k$. The exactness of the second line is consequence of the surjectivity of $\varphi_{-1}$ which is consequence of  the equality $\varphi_{-1}(\partial_{k-1}) = 1/\sigma_k$.\\
As  $\mathcal{Q}$ is the obvious quotient the injectivity of the induced map $\chi$ is easily obtained by a diagram chasing. $\hfill \blacksquare$\\

 The local  solutions of $\mathcal{N}_{-1}$ are the $1/z_j, j \in [1, k]$ and  the local solutions of $\mathcal{N}_{-1}^*$ are the $1/z_j - 1/z_h$ which generate the linear combinations of the $1/z_j$ which are killed by $U_1 = \sum_{j=1}^k z_j^2.\frac{\partial}{\partial z_j}$.\\
 
 \parag{Conclusion} For  each integer  $p \geq 2$ define
  $$\mathcal{N}_p^{\square} := D_N\big/\mathcal{I} + D_N.(U_0 - p) + D_N.U_{-1}^{p}$$
  and  $ \mathcal{N}_{-p}^* := ker(\varphi_{-p})$, where $\varphi_{-p} : \mathcal{N}_{-p} \to \mathcal{O}_N[*\sigma_k]$ is given by $\varphi_{-p}(1) = U_{-1}^{p-1}[\sigma_{k-1}/\sigma_k]$.
   Then we have the chain of isomorphisms:
 $$\xymatrix{\dots  \mathcal{N}_{-p-1}^* \ar@/_/ [r]_{\mathcal{T}_{-p}} & \mathcal{N}_{-p}^* \ar@/_/ [l]_{\square.U_1}\dots \ar@/_/[r]_{\mathcal{T}_{-1}} \ar@/_/[l]_{\square.U_1} &\mathcal{N}_{-1}^* \ar@/_/[r]_{\mathcal{T}_0}  \ar@/_/[l]_{\square.U_1} & \mathcal{N}_0^{\square} & \mathcal{N}_1^{\square} \ar@/_/[l]_{\square.U_1} \ar@/_/ [r]_{\mathcal{T}_{2}} & \ar@/_/[l]_{\square.U_1} \mathcal{N}_2^\square \dots  \ar@/_/ [r]_{\mathcal{T}_{p-1}} &\ar@/_/ [l]_{\square.U_1} \mathcal{N}_p^{\square}  \ar@/_/ [r]_{\mathcal{T}_{p}} & \mathcal{N}_{p+1}^{\square} \ar@/_/ [l]_{\square.U_1}  \dots }$$
 where $\mathcal{T}_p := \square.U_{-1}$ is given by right multiplication by $U_{-1}$.

 \subsection{Some higher order solutions of $\mathcal{N}_p$ for $p \in \mathbb{N}$}

Let $N := \C^{k}$ with coordinates $\sigma_{1}, \dots, \sigma_{k}$ and note $D_{N}$ the sheaf of (holomorphic) differential operators on $N$ and $\mathcal{D}b_{N}^{p, q}$ the sheaf of $(p, q)-$currents on $N$.\\
Recall that $\mathcal{D}b_{N}^{p, q}$ is a left  $D_{N}-$module and that we have the following theorem due to M. Kashiwara (see  \cite{[K.86]})

\begin{thm}\label{K-vanishing}
For any regular holonomic $D_{N}-$module $\mathcal{N}$ and any integer $j \geq 1$ we have
$$ Ext_{D_{N}}^{j}(\mathcal{N}, \mathcal{D}b_{N}^{0, p}) = 0 .$$
\end{thm} 

Note that the case $p \geq 1$ is an obvious consequence of the case $p = 0$ as $\mathcal{D}b_{N}^{0, p}$ is the direct sum of $C_{k}^{p}$ copies of $\mathcal{D}b_{N}^{0, 0}$ as a left $D_{N}-$module.

\begin{cor}\label{acyclique}
For any regular holonomic $D_{N}-$module $\mathcal{N}$ and any integer $j \geq 0$ we have a natural isomorphism of sheaves of $\C-$vector spaces
$$  Sol^{j}(\mathcal{N}) :=  Ext_{D_{N}}^{j}(\mathcal{N}, \mathcal{O}_{N}) \simeq H^{j}\big((Hom_{D_{N}}(\mathcal{N}, \mathcal{D}b_{N}^{0, \bullet}), \bar\partial^{\bullet})\big) .$$
\end{cor}

For instance, if $\mathcal{N} := D_{N}\big/\mathcal{J}$ is a regular holonomic system (where $\mathcal{J}$ is a coherent left ideal in $D_{N}$), we have a natural isomorphism of sheaves of complex vector spaces, for each $j$:
$$ Sol^{j}(D_{N}\big/\mathcal{J}) \simeq \{T \in \mathcal{D}b_{N}^{0, j} / \ /  \mathcal{J}.T = 0 , \ \bar\partial T = 0 \} \Big/ \bar\partial\big( \{T \in \mathcal{D}b_{N}^{0, j-1} / \ /   \mathcal{J}.T = 0 \}\big) .$$

\parag{Proof} As the Dolbeault-Grothendieck complex $ (\mathcal{D}b_{N}^{0, \bullet}, \bar\partial^{\bullet})$ is a resolution of $\mathcal{O}_{N}$ by $D_{N}-$modules for which the functor 
$$ \mathcal{N} \mapsto Hom_{D_{N}}(\mathcal{N}, -) $$
is exact, thanks to the previous theorem, the conclusion follows by degeneracy of the spectral sequence.$\hfill \blacksquare$\\

\begin{prop}\label{non vanishing}
Let $\sigma^{0}$ be a point the hypersurface $\{\sigma_{k} = 0 \}$  in $N$ and let $d$ be the multiplicity of the root $0$ in $P_{\sigma^{0}}$. Let $U$ be a small open neighborhood of $\sigma^{0}$ in $N$ on which there exists a holomorphic map $f : U \to Sym^{d}(\C)$ whose value at $\sigma \in U$ is the $d-$tuple of roots of $P_{\sigma}$ which are near by $0$\footnote{To be more precise, let $D$ be an open disc with center $0$ in $\C$ such that $\bar D$ contains only the root $0$ of $P_{\sigma^{0}}$. Then choose $U$ small enough such that for all $\sigma \in U$ the polynomial $P_{\sigma}$ has exactly $d$ roots in $D$.}.\\
Then define for $q \in \mathbb{N}$ the distribution  on $U$ (given by a locally integrable function)
\begin{equation}
T_{q}(\sigma) = \sum_{j=1}^{d} \  z_{j}^{q}.Log\vert z_{j}\vert^{2} \quad {\rm where} \  [z_{1}, \dots, z_{d}] = f(\sigma) 
\end{equation}
Then the current $\bar\partial T_{q}$ defines a section on $U$ of the sheaf $Sol^{1}(\mathcal{N}_{q})$ such that its germ at a point $\sigma^{0}$ is non zero in $Sol^{1}(\mathcal{N}_{q})_{\sigma^{0}}$.\\
\end{prop}

\parag{Proof} Let  $pr : H \to \C$ and $\pi : H \to N$ are the projections, where
 $$H := \{(\sigma, z) \in N\times \C \ / \ P_\sigma(z) = 0\}.$$
  We may assume that the open set $pr(\pi^{-1}(U))$ is the disjoint union of $D$ with an open set $\Omega$ in $\C$. Then if we define the locally integrable function $f : D \cup \Omega$ as $f(z) = z^{q}.Log\vert z\vert^{2}$ on $D$ and $f \equiv 0$ on $\Omega$ we have $T_{q}(\sigma) = \pi_{*}(f)(\sigma) = \sum_{j=1}^{k} f(z_{j})$ where $z_{1}, \dots, z_{k}$ are the roots of $P_{\sigma}$. It is then easy to verify that $\mathcal{I}.T_{q} = 0$ and that $(U_{0} - q).T_{q} = N_{q}(\sigma)$ the $q-$th Newton function of the $d-$tuple $d(\sigma)$ of roots of $P_\sigma$ which are in $D$. So it is holomorphic on $U$. Then the $(0, 1)-$current $\bar\partial T_{q}$ is $\bar\partial-$closed and is killed by $\mathcal{J}_{q}$. Then, thanks to corollary \ref{acyclique} it induces a section on $U$ of the sheaf $Sol^{1}(\mathcal{N}_{q})$.\\
Fix now $\tau \in U$ and assume that the germ at $\tau$ of the previous section vanishes. Then there exists on an open polydisc $V$ with center  $\tau$ in $U$ and a $(0, 0)-$current $S$ satisfying
$$ \mathcal{I}.S = 0, \quad (U_{0} - q).S = 0 \quad {\rm and} \quad \bar\partial S = \bar\partial T .$$
Then we may write $S = T- F$ where $F$  is holomorphic on $V$. But then $F$ satisfies also $\mathcal{I}.F = 0$ and $(U_{0}- q).F(\sigma) = N_{q}(d(\sigma))$ for all $\sigma \in V$. The first equation implies that $F$ is a global trace  function on $V$ (up to shrink $V$ around $\tau$ if necessary)  and using  lemma 3.1.2 in \cite{[B.19]} we se that, up to a locally constant function on $D\cup \Omega$, $(U_{0} - q).F$ is the trace of a holomorphic function $h$ define by $h(z) = z^{q}$ on $D$ and $0$ on $\Omega$.  But if $F = Trace(g)$ where $g$ is holomorphic on $D \cup \Omega$ this implies
$$ z.\frac{\partial g}{\partial z}(z) = h(z) + k(z) $$
where $k$ is constant equal to $k$ on $D$. So, on $D$ the meromorphic  function $G := g\big/z^{q}$ satisfies
$$G'(z) = 1/z + k/z^{q+1} .$$
This is clearly impossible for $q \geq 1$. For $q= 0$ this gives that  $G = g $ is constant and so is $F = Trace(g)$. But then $U_{0}.F = d$ is impossible for $d \geq 1$. This shows that at each point $\sigma^{0}$ of the hyper-suface $\{ \sigma_{k} = 0 \}$ in $N$ the germ induced by $\bar\partial T_q$ in $ Sol^{1}(\mathcal{N}_{q})_{\sigma^{0}}$ is not zero. So the support of the sheaf $Sol^{1}(\mathcal{N}_{q})$ contains this hyper-surface for each $q \in \mathbb{N}$.$\hfill \blacksquare$\\

\parag{Remark}  The exact sequence
 $$0 \to \mathcal{M} \overset{(U_{0}-q)}{\to} \mathcal{M} \to \mathcal{N}_{q} \to 0 $$
   gives a long exact sequence
$$ 0 \to Sol^{0}(\mathcal{N}_{q}) \to Sol^{0}(\mathcal{M}) \overset{U_{0}-q}{\to} Sol^{0}(\mathcal{M}) \overset{\partial}{\to}  Sol^{1}(\mathcal{N}_{q}) \to Sol^{1}(\mathcal{M}) \to \dots $$
and it is clear that  the germ at the origin of the Newton polynomial $N_{q}$ in $ Sol^{0}(\mathcal{M})$ is not in the image of $U_{0}- q$. Our computation above shows that the image of the germ of $N_{q}$  at the point $0 $ is mapped by the connector $\partial$ to the germ in $Sol^{1}(\mathcal{N}_{q})_{0}$ which is constructed above.

\section{An application}

We shall consider now the universal monic degree $k$ equation near the point $\sigma^0$ defined by $\sigma^0_1 = \sigma^0_2 = \dots = \sigma^0_{k-1} = 0$ and $\sigma^0_k = -1$. We shall denote by $z(\sigma)$ the root of $P_{\sigma^0 + \sigma}(z) = 0$ which is  near the (simple) root $-1$ of the equation $P_{\sigma^0 + \sigma}(z) = z^k -(-1)^k = 0$ when $\sigma$ is small enough. Define 
\begin{equation}
F(\sigma^0 + \sigma) := z(\sigma) - \sigma_1/k := \sum_{\alpha \in \mathbb{N}^k} C_\alpha.\frac{\sigma^\alpha}{\alpha!} 
\end{equation}
 the Taylor expansion at the point $\sigma^0$  of  $z(\sigma) - \sigma_1/k$ which a solution near $\sigma^0$ of the $D_N-$module $\mathcal{N}^\square_1$ (see the theorem \ref{minimality}).\\
An easy consequence of the results in the paragraph 4.1  is the following theorem.

\begin{thm}\label{7/12 1}
The following differential operators annihilate the function $F $ in a neighborhood of $\sigma^0$, where we note $\partial_h$ for the partial derivative relative to $\sigma_h$.
\begin{enumerate}
\item $A_{p, q} := \frac{\partial^2}{\partial p\partial q} - \frac{\partial^2}{\partial_{p+1}\partial_{q-1}} \quad \forall p \in [1, k-1] \ {\rm and} \ \forall q \in [2, k]$.
\item $\hat{U}_0 - 1 := \sum_{h=1}^k h.\sigma_h.\partial_h - k.\partial_k - 1 .$
\item $ U_{-1} := \sum_{h=0}^{k-1} (k-h).\sigma_h.\partial_{h+1}$ \  with the convention $\sigma_0 \equiv 1$
\end{enumerate}
\end{thm}

\parag{Proof} This is consequence of the fact that $F$ is a  solution in an open neighborhood of $\sigma^0$ of  the regular holonomic system
 $\mathcal{N}_1^{\square} \simeq D_N\big/\mathcal{A} + D_N.(U_0 - 1) + D_N.(U_{-1}) $. Remark that the operator $\hat{U}_0$ is the expression of $U_0$ in the coordinates centered at $\sigma^0$. The other operators have
 in these coordinates the same expression than in the usual coordinates centered at the origin.$\hfill \blacksquare$\\
 
 \begin{cor}\label{7/12 2}
 The coefficients $C_\alpha$ is the expansion $(1)$ only depend on the integers  $q := \vert \alpha\vert = \sum_{h=1}^k \alpha_h$ and  $w(\alpha) := \sum_{h=1}^k h.\alpha_h$ so we may rewrite the expansion $(1)$ with the convention $C_{q, r} = 0$ when $r \not\in [q, k.q]$:
 \begin{equation}
 F(\sigma^0 + \sigma) = \sum_{q, r} C_{q, r}.m_{q, r}(\sigma)  
 \end{equation}
 where for $q \in \mathbb{N}$ and $r \in [q, k.q]$ we define the polynomial $m_{q, r} \in \C[\sigma]$ by the formula
 $$m_{q, r}(\sigma) := \sum_{\vert \alpha\vert = q, w(\alpha) = r}  \frac{\sigma^\alpha}{\alpha!}$$
 \end{cor}
 
 \parag{Proof} This is obvious consequence of the description of the holomorphic functions which are annihilated by the differential operators $A_{p, q}$ for all $p \in [1, k-1]$ and $q \in [2, k]$ (see the paragraph 2.1) which generate the left ideal $\mathcal{A}$ in $D_N$.$\hfill \blacksquare$\\

\begin{prop}\label{7/12 3}
We have the following formulas, with the conventions $m_{q, r} = 0$ for $r \not \in [q, k.q] $ (in particular for $q < 0$ or $r < 0$):
\begin{enumerate}
\item $(\hat{U}_0 -1)(m_{q, r}) = (r-1).m_{q, r} - k.m_{q-1, r-k}$   
\item $U_{-1}(m_{q, r}) = (k.q - r + 1).m_{q, r-1} +  k.m_{q-1, r-1}$.
\end{enumerate}
\end{prop}

\parag{Proof} The first formula is a direct consequence of the easy formulas
  $$U_0(\sigma^\alpha/\alpha!) = w(\alpha).\sigma^\alpha/\alpha! \quad {\rm and} \quad \partial_k(\sigma^\alpha/\alpha!)= \sigma^\beta/\beta!$$
   when $\alpha_k \geq 1$, with $\beta + 1_k = \alpha$ and $$\partial_k(\sigma^\alpha/\alpha!)= 0 \quad {\rm when} \quad \alpha_k = 0.$$
The second formula is little more tricky:\\
For $h \in [2, k-1]$ we have 
$$\sigma_h.\partial_{h+1}(\sigma^\alpha/\alpha!) = \beta_h.\sigma^\beta/\beta! \quad {\rm when} \quad \alpha_{h+1} \geq 1$$
 with $\beta + 1_{h+1} = \alpha + 1_h$, and
 $$\sigma_h.\partial_{h+1}(\sigma^\alpha/\alpha!) = 0 \quad{\rm when}  \quad \alpha_{h+1} = 0.$$
  More over, for any $\beta$ with $\vert \beta \vert = q-1$ and $w(\beta) = r-1$ and for each  $h \in [2, k]$ there exists exactly one $\alpha$ if 
 $\beta_h \not= 0$ with\quad $\sigma_h.\partial_{h+1}(\sigma^\alpha/\alpha!) =  \beta_h.\sigma^\beta/\beta!$, and it satisfies $\vert \alpha \vert = q$ and $w(\alpha) = r$, and no such $\alpha$ exists if $\beta_h = 0$. This means that
 that $\sigma_h.\partial_{h+1}(m_{q, r}) $ contains $\sigma^\beta/\beta!$ with the coefficient $\beta_h$.\\
 For $h = 1$ the situation is simpler: $ \partial_1(\sigma^\alpha/\alpha!) = \sigma^\beta/\beta!$ \  when $\alpha_1 \geq 1$ with $\beta + 1_1 = \alpha$, and $\partial_1(\sigma^\alpha/\alpha!) = 0$ when $\alpha_1 = 0$. \\
 Then for each $\beta$ with $\vert \beta \vert = q-1$ and $w(\beta) = r-1$ there exists a unique $\alpha$ such that $ \partial_1(\sigma^\alpha/\alpha!) = \sigma^\beta/\beta!$\ and it satisfies $\vert \alpha \vert = q$ and $w(\alpha) = r$.
 So we conclude that 
 \begin{align*}
 &  U_{-1}(m_{q, r}) = \sum_{\vert \alpha\vert = q, w(\alpha) = r} \sum_{h=0}^k (k-h).\sigma_h.\partial_{h+1}(\frac{\sigma^\alpha}{\alpha!}) \\
 & \qquad \quad  = k.m_{q-1, r-1} + \sum_{\vert\beta\vert = q, w(\beta) = r-1}\sum_{h=1}^k (k-h).\beta_h.\frac{\sigma^\beta}{\beta!}\\
 & \qquad \quad  = k.m_{q-1, r-1} + \big( k.(q-1) - (r-1)\big).m_{q, r-1} 
 \end{align*}
 concluding the proof.$\hfill \blacksquare$\\
 
 Taking in account the equations 2. and 3. of the theorem \ref{7/12 1} (the equations 1. are used already in the corollary \ref{7/12 2}) we obtain:
 
 \begin{cor}\label{7/12 4}
The coefficients $C_{q, r}$ of the Taylor expansion $(1)$ satisfies the relations:
\begin{align*}
& (r-1).C_{q, r} - k.C_{q+1, r+k} = 0  \quad \forall q \geq 1, \quad  \forall r \in [q, k.q]  \tag{A} \\
& (k.q - r + 1).C_{q, r} + k.C_{q+1, r} = 0 \quad \forall q \geq 1, \quad \forall r \in [q+1, k.q] . \tag{B}
\end{align*}
The  formula $(B)$ gives, for each $r \geq 2$ and each $s \in \mathbb{N}$ such that  $0 \leq s \leq \frac{(k-1).r}{k}$ 
\begin{equation*}
C_{r-s, r} = \frac{(-1)^s.C_{r,r}}{\prod_{j=1}^s (r - j - (r-1)/k) } \tag{$B^*$}
\end{equation*}
The formula $(A)$ gives for each $r \geq 1$:
\begin{equation*}
C_{r+k, r+k} = (-1)^{k-1} \frac{r-1}{k}.\prod_{p=0}^{k-2} \big(r+p-(r-1)/k\big).C_{r,r} \tag{$A^*$}
\end{equation*}
Moreover we have
\begin{equation*}
C_{q, r} = 0 \quad \forall q \geq 2 \ {\rm and} \quad \forall r \equiv 1 \ {\rm modulo} \ k, \ r \in [q, k.q] \tag{$C$}
\end{equation*}
\end{cor}

\parag{Proof}  Looking at the coefficient of $m_{q,r}$ for $q \geq 1$ and $r \in [q, k.q]$ in  the equality $(\hat{U}_0 - 1)[F] \equiv 0$  gives the  gives $(A)$. Looking at the coefficient of $m_{q, r-1}$ for $q \geq 1$ and $r \in [q+1, k.q]$ in 
$U_{-1}[F] \equiv 0$ gives $(B)$.\\
The formula $(B^*)$ is a direct consequence of the formula $(B)$ with $q := r-s$  by an easy induction on $s \in [0, \frac{(k-1).r}{k}]$.\\
Using formula $(B^*)$ with $r' = r+k$ and $s = k-1$ we obtain
\begin{equation*}
C_{r+1, r+k} = \frac{(-1)^{k-1}.C_{r+k, r+k}}{\prod_{p=0}^{k-2} (r+p - (r-1)/k)} \tag{@}
\end{equation*}
Combining this formula with the formula $(A)$ with $q = r \geq 1$ which gives $$C_{r+1, r+k} = \frac{r-1}{k}.C_{r, r}  $$
and we obtain the formula $(A^*)$.\\
Formula $(C)$ is a direct consequence of $(A^*)$ for $r = 1$ with an induction on $a \geq 1$ when $r = 1+ a.k$.$\hfill \blacksquare$\\

We shall see below that the vanishing of  $C_{1+a.k, 1+ a.k}$ is also valid for $a = 0$ giving $C_{q, 1+ a.k} = 0$ for any $q \geq 1$ and any $a$ such that $q \leq 1+ a.k \leq k.q$.

\parag{Remark} It is enough to compute $C_{0, 0}$ and $C_{1, h}$ for each $h \in [1, k]$ to determine all coefficients $C_{q, r}$ in $(1)$ with $r \in [q, k.q]$:\\
The formula $(B)$ determines $C_{h,h}, h \in [2, k]$ from $C_{1, h}$ with $r = h $  and $s = h-1$. Then the formulas $(A^*)$ and $(C)$ gives $C_{r, r}$ for any $r \geq 1$.
Then the formula $(B)$ completes the computation of $C_{q, r}$ for any $q \geq 0$ and any $r \in [q, k.q]$.

\begin{lemma}\label{final}
We have the following values:
\begin{align*}
& C_{0, 0} = -1  \\
& C_{1, 1} = 0 \\
& C_{1, h} = 1/k   \quad           {\rm for} \ h \in [2, k].
\end{align*}
\end{lemma}

\parag{proof} The value of $C_{0,0}$ is $F(\sigma^0)$ which is $ -1 $ by definition of $F$. The values of $C_{1, h}$ is the derivative $\partial_h(F + \sigma_1/k)(\sigma^0)$ because we have $m_{1, h} = \sigma_h$ for $h \in [1, k]$. So it is enough to make an order $1$ expansion of $F$ at $\sigma^0$ to compute the values of the $C_{1, h}, h \in [1, k]$. This is given by the following computation at the first order of $P_{\sigma^0 + \sigma}(z(\sigma)) \equiv 0$,  where we define
$$z(\sigma) = -1 +  \sum_{j=1}^k c_j.\sigma_j + o(\vert\vert \sigma\vert\vert) $$
 which gives $c_1 = C_{1, 1} + 1/k$ and  $c_j := C_{1, j}$ for $j \in [2, k]$ and then:
\begin{align*}
& \big(-1 + \sum_{h=1}^k c_h.\sigma_h\big)^k + \sum_{j=1}^k (-1)^j.\sigma_j.\big(-1 + \sum_{h=1}^k c_h.\sigma_h\big)^{k-j} - (-1)^k = o(\vert\vert \sigma\vert\vert) \\
& (-1)^{k-1}.k.(\sum_{h=1}^k c_h.\sigma_h) + \sum_{j=1}^k (-1)^{k}.\sigma_j =  o(\vert\vert \sigma\vert\vert)  \quad {\rm and \ so} \\
& c_j = 1/k \quad \forall j \in [1, k] \\
\end{align*}
Then  $C_{1, 1} = 0$ and  $C_{1, j} = c_j  = 1/k $ for each $j \in [2, k]$.$\hfill \blacksquare$\\

Then, thanks to the formula $(B^*)$ for $r = h \in [2, k]$ and $s = h-1$ we find
$$ C_{h, h} = (-1)^{k-h} \frac{\prod_{j=1}^{h-1} (h-j-(h-1)/k)}{\prod_{p=0}^{k-2} (h+p-(h-1)/k)} \quad \forall h \in [2, k].$$

\newpage

\section{Appendix: The surface $S(k)$}

For $k  \geq 2$ an integer and $\alpha \in \mathbb{N}^k$ define
\begin{itemize}
\item the {\bf length} of $\alpha$ given by  $\vert \alpha\vert := \sum_{h=1}^k \alpha_h $; 
\item the {\bf  weight} of $\alpha$ given by  $w(\alpha) := \sum_{h=1}^k  h.\alpha_h $.
\end{itemize}
We shall say that $\alpha$ and $\beta$ are {\bf equivalent}, noted by  $\alpha \sharp \beta$, when $\vert \alpha \vert = \vert \beta \vert$ and $w(\alpha) = w(\alpha)$. \\
Remark that for any $\gamma \in \mathbb{N}^k$ we have $(\alpha + \gamma) \sharp (\beta + \gamma )$ if and only if $\alpha \sharp \beta$.\\

Let $A$ be a $\C-$algebra which is commutative, unitary and integral. In the algebra  $A[x_1, \dots, x_k]$ let $IS(k)$  be  the ideal generated by the polynomials $x_p.x_q - x_{p+1}.x_{q-1}$ for all $p \in [1, k-1] $ and $q \in [2, k]$.\\
We shall say that the two monomials $x^\alpha$ and $x^\beta$ in  $A[x_1, \dots, x_k]$ are equivalent when $\alpha$ and $\beta$ are equivalent. In this case we shall also write $x^\alpha \sharp x^\beta$. \\
Remark that for any $p \in [1, k-1]$ and any $q \in [2, k]$ $x_p.x_q$ is equivalent to $x_{p+1}.x_{q-1}$.\\

For a monomial $m := x^\alpha$ we define its length by $l(m) := \vert \alpha\vert$ and its weight $w(m) := w(\alpha)$.\\

Our first result is the following characterization of the elements in $IS(k)$.

\begin{prop}\label{le retour}
 Two monomials $x^\alpha$ and $x^\beta$ in $A[x_1, \dots, x_k]$ are equivalent if and only if $x^\alpha - x^\beta $ is in $IS(k)$.
 \end{prop}
 
The proof of this proposition will need a preliminary lemma and the next definition.

\begin{defn}\label{minimal}
We shall say that a monomial $m$ in $A[x_1, \dots, x_k]$ is {\bf minimal} when it has one of the following forms:
\begin{enumerate}
\item there exists $p, q$ in $\mathbb{N}$ such that $m =  x_1^p.x_k^q$ ; 
\item there exists $p, q$ in $\mathbb{N}$ and $j \in [2, k-1] $ such that $m =  x_1^p.x_j.x_k^q$.
\end{enumerate}
\end{defn}

\parag{Remark} Any monomial (minimal or not) is not in the ideal $IS(k)$ because  the point $x_1 = x_2 = \dots = x_k = 1$ is not in $\vert S(k)\vert$ the common set of zeros in $A^k$ of the generators of $IS(k)$ and a monomial does not vanish at this point.

\begin{lemma}\label{reduction}
For each $\alpha \in \mathbb{N}^k$ there exists an unique minimal monomial $ x^{\mu(\alpha)}$ such that $x^\alpha \sharp x^{\mu(\alpha)}$. Moreover, for each $\alpha$  $x^\alpha - x^{\mu(\alpha)}$ is in $IS(k)$
\end{lemma}

\parag{proof} Let us begin by proving the uniqueness assertion.\\
We have to show that two minimal monomials which are equivalent are equal. If both are in case 1. (so $m := x_1^p.x_k^q$)  this is obvious as the length is equal to $l(m) = p + q$ and the weight is $w(m) = p + k.q$ and then $(k-1).q = w(m) - l(m)$ proving the uniqueness of $q$ and then of $p$.\\
If both are in case 2. let $m := x_1^p.x_j.x_k^q$ and $m' = x_1^{p'}.x_{j'}.x_k^{q'}$ then we have 
\begin{align*}
& l(m) = p + 1 + q = l(m') = p' + 1 + q' \quad {\rm and} \\
&  w(m) = p + j + k.q = w(m') = p' + j' + k.q' \quad {\rm which \ imply} \\
& j- j' = (k-1).(q' - q) \quad {\rm with} \quad \vert j - j'\vert \in [0, k-3].\\
& {\rm So}  \quad j = j' \quad {\rm and \ then} \quad  q = q' \quad {\rm and} \quad p = p'.
\end{align*}

If $m = x_1^p.x_j.x_k^q$  and $m' =  x_1^{p'}.x_k^{q'}$ we have
 \begin{align*}
 & (m) = p + 1 + q = l(m') =  l(m') = p' + q' \quad {\rm and} \\
& w(m) = p + j + k.q = w(m') = p' + k.q' \quad {\rm which \ imply} \\
 & j - 1 = (k-1).(q' - q) \quad {\rm with} \quad  j \in [2, k-1] \quad {\rm  and \ this \ is\ impossible}.
 \end{align*}
The assertion of existence  is clear for $\vert \alpha \vert = 0, 1$. We shall prove the existence of $\mu(\alpha)$  by an induction on the length  $\vert \alpha\vert$ of $\alpha$.\\
Assume that the lemma is proved for all $\beta \in \mathbb{N}^k$ with length $1 \leq \vert \beta \vert < \vert \alpha\vert$. Then write $x^\alpha = x_r.x^\beta$ for some $r \in [1, k]$. By the induction hypothesis we know that there exists a minimal monomial $x^{\mu(\beta)}$ with $x^\beta \sharp x^{\mu(\beta)}$. Then we obtain that $x^\alpha \sharp x_r.x^{\mu(\beta)}$. If $x^{\mu(\beta)} =  x_1^p.x_k^q$, then $x_r. x_1^p.x_k^q$ is minimal for any choice of $r \in [1, k]$. If $ x^{\mu(\beta)} =  x_1^p.x_j.x_k^q$ then remark that we have $x_r.x_j \sharp x_1.x_{r+j-1}$ for $r+j-1 \leq k$ and $x_r.x_j \sharp x_k.x_{r+j-k}$ for $r+j \geq k+1$ and this allows to conclude the induction.\\
Remark that if, in the induction above, we assume that $x^\beta - x^{\mu(\beta)}$ belongs to $IS(k)$ we obtain that $x^\alpha - x^{\mu(\alpha)}$ is also in $IS(k)$;  for instance in the case $x^{\mu(\beta)} =  x_1^p.x_j.x_k^q$ 
\begin{align*}
&  x^\alpha - x^{\mu(\alpha)} = x_r.(x^\beta - x^{\mu(\beta)}) +  (x_r.x_j -  x_1.x_{r+j-1}).x_1^p.x_k^q \quad {\rm for} \quad  r+j \leq k+1, \\
x&  x^\alpha - x^{\mu(\alpha)}  = x_r.(x^\beta - x^{\mu(\beta)}) +  (x_r.x_j -  x_k.x_{r+j-k}).x_1^p.x_k^q  \quad {\rm for} \quad r+j  \geq k+2. 
\end{align*}
The other cases are analogous.$\hfill \blacksquare$

\parag{Proof of the proposition \ref{le retour}}The previous lemma gives that  $x^\alpha \sharp x^\beta$ implies  $x^\alpha - x^{\mu(\alpha)}$ and $x^\beta - x^{\mu(\alpha)}$ are in $IS(k)$, so also $x^\alpha - x^\beta$. Conversely, assume that $x^\alpha - x^\beta $ is in $IS(k)$. As the ideal $IS(k)$ is homogeneous (in the sense of length) if $\l(\alpha) \not= l(\beta)$ we conclude that both $x^\alpha$ and $x^\beta$ are in $IS(k)$. This contradicts the remark following the definition \ref{minimal}.\\
In a similar way the ideal $IS(k)$ is quasi-homogeneous in the sense of the weight $w$. So if $w(\alpha) \not= w(\beta)$ then  $x^\alpha$ and $x^\beta$ are in $IS(k)$ which is again impossible. So $x^\alpha - x^\beta $ is in $IS(k)$ implies that $\alpha \sharp \beta$.$\hfill \blacksquare$

\begin{cor}\label{quotient}
For any $q \in \mathbb{N}$ and any $r \in [q, k.q]$ there exists a minimal monomial $\mu_{q, r}$ (necessarily unique) such that $\vert \mu_{q, r}\vert = q$ and $w(\mu_{q, r}) = r$.
\end{cor}

\parag{proof} The assertion is clear for $q = 0, 1$. So let us prove it by induction on $q$. So let $q \geq 2$ and let $r \in [q, k.q]$, and assume that we know that $\mu_{q', r'}$ exists for any $q' \leq q-1$ and any $r' \in [q', k.q']$. If $r$ is in $[q, k.(q-1)+ 1]$, then $r-1$ is in $[q-1, k.(q-1)]$ and $\mu_{q-1, r-1}$ exists. So $\mu_{q, r} := x_1.m_{q-1, r-1}$ is the solution. \\
If $r $ is in $[k.(q-1)+1, k.q]$ then $r-k$ is in $[q-1, k.(q-1)]$and,  because for $q \geq 2$ we have $k.q -2k +1 \geq q-1$ and also $r-k \leq k.q-k \leq k.(q-1)$, $\mu_{q-1, r-k}$ is defined and $\mu_{q, r} := x_k.\mu_{q-1, r-k}$ is the solution.$\hfill \blacksquare$\\

\begin{prop}\label{local description}
Let $L_1 := \{\eta_1 = 0\} \cap S(k)$ and $L_k := \{\eta_k = 0 \} \cap S(k)$. Then $L_1$ is the line directed by the vector $(0, \dots, 0, 1)$ and $L_k$ the line directed by the vector $(1, 0, \dots, 0)$.
The maps $\varphi_1 : S(k)\setminus L_1 \to \C^*\times \C$ and $\varphi_k : S(k) \setminus L_k \to \C^*\times \C$ which are defined by the formulas
\begin{equation}
\varphi_1(\eta) := (\eta_1, -\eta_2/\eta_1) \quad {\rm and} \quad \varphi_k(\eta) := (\eta_k, -\eta_{k-1}/\eta_k)
\end{equation}
are isomorphisms. So $S(k)\setminus \{0\}$ is smooth and connected surface.
\end{prop}

\parag{Proof of the proposition \ref{local description}} Consider the holomorphic map
$$ \psi_1 : \C^2 \to S(k) \quad (\zeta_0, \zeta_1) \mapsto  x_h := \zeta_0.(-\zeta_1)^{h-1}  \ \forall h \in [1,k] .$$
It induces the inverse to the map $\varphi_1$ on $\zeta_0 \not= 0$ and the map $\psi_k$ defined by \\  $x_h = (-\zeta_0)^{k-h}.\zeta_1 \quad \forall h \in [1, k]$ gives the inverse of $\varphi_k$
on $\zeta_1 \not= 0$.$\hfill \blacksquare$\\

\begin{cor}\label{irreduct.}
The ideal $IS(k)$  of $\C[x]$  is prime. Moreover $(S_k)$ is a normal surface.
\end{cor}

\parag{Proof of the corollary \ref{irreduct.}} The only point which is not a direct consequence of the previous proposition is the normality of $S(k)$. But as the blow-up of the maximal ideal at the origin in $S(k)$ gives
a desingularization of $S(k)$ with the rational curve\footnote{ the image of the map $\psi_0 : \mathbb{P}_1 \to \mathbb{P}_{k-1}$ by $ \xi_h = \zeta_1^{h-1}.\zeta_0^{k-h}$.} over the origin in $S(k)$.So is a rational singular point and $S(k)$ is normal.$\hfill \blacksquare$\\

\end{document}